
\input amstex
\vcorrection{-0.5in}

\documentstyle{amsppt}

\define\OmYE{\Omega^1_Y(\log E)}
\define\OmU{\Omega^1_U}
\define\EE{\Cal E}
\define\EX{\Cal E_X}
\define\FF{\Cal F}
\define\FY{\Cal F_Y}
\define\KK{\Cal K}
\define\Kn{\Cal K_n}
\define\MM{\Cal M}
\define\tMM{\tilde{\Cal M}}

\define\Ax{A_x}
\define\Axi{A_{x,i}}

\define\Byi{B_{y,i}}

\define\Hom{\Cal Hom}
\define\End{\Cal End}
\define\ZZ{\Bbb Z}
\define\hZZ{\hat{\Bbb Z}}
\define\QQ{\Bbb Q}
\define\Zm{\ZZ/m}

\define\Ql{\QQ_\ell}
\define\PP{\Bbb P^1}
\define\CC{\Bbb C}
\define\Gm{\Bbb G_m}
\define\aaa{\bold a}
\define\an{^{an}}
\define\ab{^{ab}}
\define\et{^{et}}
\define\tame{^{tame}}
\define\llll{\lambda}
\define\kk{\kappa}
\define\isom{\overset \sim\to\longrightarrow}

\define\PA#1{\pi_1(#1)^{\roman{ab}}}
\define\PAT#1{\pi_1(#1)^{\roman{ab,tame}}}
\define\PI#1#2{\pi_1(#1,#2)}
\define\QZ{\Bbb Q/\Bbb Z}
\define\kxx{\kappa(x)^\times}
\define\kx{\kappa(x)}
\define\kxt{{\kappa(\tilde x)}}
\define\kyx{\kappa(y)^\times}
\define\ky{\kappa(y)}
\define\Card{\roman{Card}\ }

\define\Ker{\roman{Ker}\ }
\define\Coker{\roman{Coker}\ }
\define\rank{\roman{rank}\ }
\define\chara{\roman{char}\ }
\define\res{\roman{res}}
\define\ssum{\roman{sum}}
\define\ord{\roman{ord}}

\define\Spec{\roman{Spec}\ }
\define\Gal{\roman{Gal}}
\define\Tr{\roman{Tr}}

\mag=\magstep1 \document
\topmatter
\title
Determinant of period integrals
\endtitle
\author
Takeshi Saito, Tomohide Terasoma
\endauthor
\affil
Department of Mathematical Sciences,
University of Tokyo,
7-3-1 Hongo Bunkyo, 113 Tokyo, Japan\\
Department of Mathematical Sciences,
University of Tokyo,
3-8-1 Komaba Meguro, 153 Tokyo,Japan
\endaffil
\endtopmatter

We prove a formula for the determinant of
period integrals.
Period integrals are comparison between
Betti cohomologies and
de Rham cohomologies.
Our formula
Theorem 1 in Section 4 expresses
the determinant of period integrals
as the product of the periods
evaluated at the relative canonical
cycles and of product of special values
of the $\Gamma$-function.
The formula is a Hodge version of
Theorem 1 of [S2]
for $\ell$-adic cohomology.
Together with this,
it gives a motivic formula,
Theorem 3 in Section 5.
It particularly implies that the
category of motives of rank 1
associated to an algebraic Hecke character
is closed under taking the determinant of cohomology.
Hence it gives a support to a conjecture of Deligne,
Conjecture 8.1 iii [D4];
a motive of rank 1 is associated to an algebraic
Hecke character.
A typical example of our formula
is that the period of Fermat hypersurfaces are
product of special values of $\Gamma$-function
(cf. proof of Lemma 5.5).
The main theorem for $X=\PP$
is a reformulation of a theorem of
the second named author [T] Theorem 1.2.
The theorem is proved by reducing to this case
by induction on dimension using a Lefschetz pencil.

Content of each section is as follows.
In Section 1,
we review generalities on integrable connections
and local systems as in [D1]
and define the determinant of period integrals.
In Section 2, first
we define relative Chow group and
relative top chern class and give
its basic properties as in [S1] Section 1.
The rest of the section is
devoted to a proof of an adelic description, Proposition 1,
of the relative Chow group.
In Section 3, the pairing with a relative 0-cycles are defined.
To define it, we develop a theory of tame symbols.
Its philosophical back ground is explained in
Remark in the text.
In Section 4, we prove the main theorem, Theorem 1.
For curves,
a stronger version, Theorem 2,
involving Deligne-Riemann-Roch [D3]
is also proved.
{}From the point of view of Riemann-Roch,
the reason why the canonical class appears in the formula
is the equality
$\sum(-1)^qch(\Omega^q_X(\log D))=
(-1)^nc_n(\Omega^1_X(log D))$ for $n=\dim X$.
In the final Section 5,
first we reformulate a formula, Theorem S (cf. Theorem 1 [S2]),
for $\ell$-adic cohomology using relative Chow group.
Finally a motivic version, Theorem 3, is proved.
There Jacobi sum \`a la Anderson [A] appears.

This paper is a full version of [ST]

\demo{Remark}
We have changed the definition of
the determinant of the period in [ST]
to its inverse,
in order to make it compatible
with [D4] 8.7.
Because of this change,
the sign in Theorem loc.cit and
Corollary of Theorem 1 in Section 4
looks apparantly inconsistent.
We apologize for inconvenience.
\enddemo

\heading
1. Connections and local systems.
\endheading

Let $k$ be a field of characteristic 0 and
$U$ be a smooth separated scheme over $k$.
An integrable connection $\nabla$ on a
locally free $\Cal O_U $-module $\EE$ of finite rank
is an additive morphism $\nabla:\EE\to\EE\otimes\OmU$
satisfying $\nabla(fe)=f\nabla e+e\otimes df$
for local sections of $f$ of $\Cal O_U $ and $e$ of $\EE$
and $\nabla^2=0:\EE\to\EE\otimes\Omega^2_U$.
Let $X$ be a smooth separated
scheme over $k$ of dimension $n$ including $U$
as the complement of a divisor
$D$ with simple normal crossings.
A divisor $D$ is called to have simple normal crossings
if its irreducible components
$D_{i,i\in I}$ are smooth and
their intersections $D_J=\bigcap_{i\in J}
D_i$ for $J\subset I$ are transversal.
Let $\Omega^q_X(\log D)$
be the sheaf of diffential
$q$-forms with logarithmic poles
along $D$.
For a coherent $\Cal O_X $-module
$\EX$, we call a logarithmic integrable connection
an additive morphism $\nabla:\EX\to\EX\otimes\Omega^1_X(log D)$
satisfying $\nabla(fe)=f\nabla e+e\otimes df$
for local sections of $f$ of $\Cal O_X $ and $e$ of $\EX$
and $\nabla^2=0$.
We call a logarithmic integrable connection $(\EX,\nabla)$
an extention of an integrable connection $(\EE,\nabla)$
on $U$ to $X$ when an isomorphism
$(\EX,\nabla)|_U\isom(\EE,\nabla)$
is given.
\proclaim{Lemma 1.1}
For an integrable connection $(\EE,\nabla)$
on $U$,
the following conditions (1)-(3) are equivalent.
\roster
\item There exists an extention $(\EX,\nabla)$ to $X$.
\item There exists an extention $(\EX,\nabla)$ to $X$
such that $\EX$ is reflexive.
\item For the generic point $\eta_i$ of each
irreducible component $D_i$ of $D$,
there exists a free
$\Cal O_{X,\eta_i}$-module
$\Cal E_i$ of finite rank
and a logarithmic integrable connection
$\nabla:\Cal E_i\to
\Cal E_i\otimes\Omega^1_X(log D)_{\eta_i}$
extending the generic fiber
$\nabla:\EE_\xi\to \EE_\xi\otimes\Omega^1_{X,\xi}$
at the generic point $\xi$ of $X$.
\endroster
\endproclaim

\demo{Proof}
It is clear that (2)
implies (1) and (3).
On the other way,
each of (1) and (3)
implies that
there exists an open subscheme
$U'$ of $X$ including $U$
such that the codimension of $X-U'\ge 2$
and a locally free extension
$\EE_{U'}$
of $\EE$ to $U'$.
If $j:U'\to X$ denotes the open immersion,
$\EX=j_*\EE_{U'}$
is a reflexive extension of $\EE$ to $X$.
\enddemo
We say an integrable connection $(\EE,\nabla)$
on $U$ is regular along $D$ if
the equivalent conditions in Lemma 1.1
are satisfied.

\proclaim{Corollary}
Let $X$ and $X'$ be proper smooth schemes
over $k$
including $U$ as the complement of
divisors $D$ and $D'$
with simple normal crossings respectively.
For an integrable connection
$(\EE,\nabla)$ on $U$,
it is regular along $D$
if and only if so is along $D'$.
\endproclaim

\demo{Proof}
We may assume there is a proper morphism
$X'\to X$ inducing identity on $U$.
Then by (1),
regularity along $D$ implies that along $D'$.
By (3) we have the inverse implication.
\enddemo
We say an integrable connection
$(\EE,\nabla)$ on $U$
is regular along boundary,
if there is one (hence for any) proper $X$ including $U$
as the complement of a divisor $D$
with simple normal crossings such that
it is regular along $D$.
The integrable connections on $U$
regular along boundary form a neutral Tannakian category.

In the following,
we always consider reflexive extension
unless we say otherwise explicitly.
Let $(\EX,\nabla)$ be a reflexive
logarithmic integrable connection
on $(X,D)$.
The de Rham complex
$$\gather
DR(\EX)=
(\EX\otimes \Omega_X^\bullet(\log D))\\
=[\EX
\overset{\nabla}\to \to
\EX\otimes\Omega^1_X(log D)
\overset{\nabla}\to \to
\EX\otimes\Omega_X^2(\log D)\to\cdots]
\endgather$$
is defined as usual.
For each component
$D_i$ of $D$,
the residue
$\res_i\nabla\in$ \linebreak
$End_{\Cal O_{D_i}}(\EX|_{D_i})$
of $\nabla$ at $D_i$
is defined as the map induced by
$$\EX
\overset{\nabla}\to \to
\EX\otimes\Omega^1_X(log D)
\overset{id\otimes\res_i}\to \to
\EX\otimes_{\Cal O_X }\Cal O_{D_i}.$$
We define the characteristic polynomial
$\Phi_{\EE_{X,i}}(T)=
\det (T-\res_i\nabla)
\in \kappa(\eta_i)[T]$.
It is in $k_i[T]$
where $k_i$ is the constant field of
$D_i$ by [D1] Chapitre II Proposition 3.10 (ii).
We say an extension $(\EX,\nabla)$
of $(\EE,\nabla)$
is small (resp. big) at $D_i$
if any integer
$n\le0$ (resp. $n\ge1$)
is not a root of the polynomial
$\Phi_{\EX,i}(T)=0.$
We say $(\EX,\nabla)$
is small if it is small
at every irreducible component $D_i$.
The extention $(\Cal O_X (-D),d)$
of $(\Cal O_U ,d)$ is small since
$\Phi_i(T)=T-1$.
An extention $\EX$ is small if and only if
its dual $\EX^*=\Hom(\EX,\Cal O_X (-D))$ is big.

\proclaim{Lemma 1.2}
Let $(\EX,\nabla)$
be a reflexive extension.
If $\res_i\nabla\in
End_{\Cal O_{D_i}}(\EX|_{D_i})$
is an isomorphism
at $\eta_i$,
it is an isomorphism on $D_i$.
\endproclaim

\demo{Proof}
We use the following elementary fact.

\proclaim{Lemma 1.3}
Let $\EE$ be a reflexive $\Cal O_X $-module on a
regular scheme $X$.
Then the complement of the open subscheme
$W$ of $X$
where $\EE$ is locally free is of
codimension $\ge 3$.
\endproclaim
\demo{Proof}
It is sufficient to show that $\EE$ is free
at a point $x$ of codimension 2.
Taking a regular divisor $D\ni x$,
it is enough to show $\EE|_D$
is torsion free.
Let $j:W\cap D\to D$ be the open immersion.
Since $\EE$ is reflexive,
we see that $\EE|_D$
is a submodule of
$j_*j^*(\EX|_D)$ and is torsion free.
\enddemo
Let $j:W\cap D_i\to D_i$
be the open immersion.
Since $\EX$ is reflexive,
we see that $\EX|_{D_i}$
is a submodule of the reflexive sheaf
$j_*j^*(\EX|_{D_i})$.
By the assumption,
$\det\res_i\nabla$ is a non-zero constant
on $W\cap D_i$. Hence
$\res_i\nabla$ induces
an automorphism of
$j_*j^*(\EX|_{D_i})$.
Now the assertion is a consequence
of the following elementary fact.

\proclaim{Lemma 1.4}
Let $f$ be an automorphism of a coherent
$\Cal O_X $-module $\MM$
on a \linebreak
noethrian scheme $X$ and
$\Cal N\subset \MM$ be
a sub-coherent $\Cal O_X $-module
stable by $f$.
Then the restriction
$f|_{\Cal N}$ is an automorphism.
\endproclaim
\demo{Proof}
We have an increasing sequence of submodules
$\Cal N
\subset f^{-1}(\Cal N)
\subset\cdots$.
Hence we have
$f^{-n}(\Cal N)= f^{-n-1}(\Cal N)$
for some $n$ and $f(\Cal N)=\Cal N$.
\enddemo
\enddemo
For an integrable connetion $(\EE,\nabla)$
on $U$ regular along the boundary,
the category of small extensions
(resp. big extensions) of $(\EE,\nabla)$
to $X$ are non empty and cofinal since
$\Phi_{\EX(\sum n_jD_j),i}(T)=
\Phi_{\EX,i}(T+n_i).$

\proclaim{Lemma 1.5}
For small (resp. big) extensions
$\EX\subset\EE'_X$ of $(\EE,\nabla)$,
the inclusion of the de Rham complexes
$DR(\EE'_X)\subset DR(\EX)$ is a quasi-isomorphism.
\endproclaim
\demo{Proof}
We give a proof for the small case
and the big case is done in the same way.
There is a sequence of
extensions
$\EX=\EE_0\subset\EE_1\subset
\cdots\subset \EE_\ell\subset\cdots\subset
\EE_m=\EX'$
of extensions of $(\EE,\nabla)$
such that
$\EE_\ell/\EE_{\ell-1}$
is an $\Cal O _{D_{i(\ell)}}$-module
for $1\le\ell\le m$.
In fact, a reflexive extension
$\EX\subset\EE"_X\subset
\EX'$
is determined by
a family of
sub $\Cal O _{X,\eta_i}$-modules
$\EE_{X,i}\subset\EE"_{X,i}\subset
\EE'_{X,i}$
such that
$\nabla(\EE"_{X,i})\subset
\EE"_{X,i}\otimes\Omega^1_X(log D)$
for each irreducible component $D_i$
(cf. proof of Lemma 1.1).

By induction, we may assume that
there exists an irreducible component $D_i$
such that the quotient
$\EX/\EE'_X$ is an $\Cal O_{D_i}$-module.
We show that the quotient
$DR(\EX)/DR(\EE'_X)
=(\EX/\EE'_X\otimes_{\Cal O_{D_i}}
\Omega^\bullet_X(\log D)|_{D_i})$
is acyclic.
We define a decreasing filtration $F$ on
$DR(\EX)/DR(\EE'_X)$ by
$$
F^p(DR(\EX)/DR(\EE'_X))^q=
\cases
0& q<p\\
\EX/\EE'_X\otimes_{\Cal O_{D_i}}\Omega^q_{D_i}(\log D|_{D_i})& q=p\\
\EX/\EE'_X\otimes_{\Cal O_{D_i}}\Omega_X^q(\log D)|_{D_i}& q>p.
\endcases
$$
Here $D|_{D_i}=
\bigcup_{j\ne i}D_j\cap D_i$.
The $gr$-complex
$Gr_F^p$
is given by
$$\res_i\nabla\otimes id:
\EX/\EE'_X\otimes_{\Cal O_{D_i}}\Omega^p_{D_i}(\log D|_{D_i})\to
\EX/\EE'_X\otimes_{\Cal O_{D_i}}\Omega^p_{D_i}(\log D|_{D_i}).$$
By the assumption that
$\EX$ is small and by Lemma 1.4,
$\res_i\nabla$ is an isomorphism
on $\EX/\EX'$. Thus Lemma 1.5 is proved.
\enddemo

For an integrable connection
$(\EE,\nabla)$
on $U$
regular along boundary,
we define the de Rham cohomology with
compact support
$H^q_c(U,DR(\EE))$
for $q\in \ZZ$.
Let $X$ be a proper smooth
$k$-scheme including $U$
as the complement of
a divisor $D$ with simple normal croosings.
Then we define
$$H^q_c(U,DR(\EE))=
\projlim_{(\EX,\nabla): \text{ small}}
H^q(X,DR(\EX)).$$
Here $(\EX,\nabla)$
are the small extension of $(\EE,\nabla)$
to $X$.
By Lemma 1.5, all the transition morphisms are isomorphisms.
It is indepenent of choice of $X$.
This fact is proved for example
by the Serre duality below.
By Lemma 1.5 and by the isomorphism
$Rj_*DR(\EE)\simeq
j_*DR(\EE)$,
the usual de Rham cohomology
$H^q(U,DR(\EE))$
is isomorphic to
$H^q(X,DR(\EX))$
for big $\EX$.
For a logarithmic connection
$\EX$ and its dual $\EX^*=\Hom(\EX,\Cal O_X (-D))$,
the pairing
$$H^p(X,DR(\EX))\otimes
H^q(X,DR(\EX^*))
\overset\cup\to \to
H^n(X,\Omega^n_X)
\overset{\Tr_{X/k}}\to\to k$$
for $p+q=2n$
is perfect by Serre duality.
Therefore the pairing
$$H^p_c(U,DR(\EE))\otimes
H^q(U,DR(\EE^*)) \to k$$
for $p+q=2n$
is perfect.
By Riemann-Roch,
the Euler number
$\chi_c(U,DR(\EE))=
\sum_p(-1)^p
\dim_kH^p_c(U,DR(\EE))$
is rank $\EE\times \chi_c(U)$,
where $\chi_c(U)$ is the Euler number
$\chi_c(U,DR(\Cal O_U ))$.

Let $k_0$ be a subfield of the complex number field $\CC$.
We assume $k$ is a finite extension of $k_0$.
Let $X, U$ and $D$ over $k$ be as above.
Let $X^{an}$ be the complex manifold
with underlying set $X(\CC)=Hom_{k_0-alg}(\CC,X)$.
An integrable connection
$(\EE,\nabla)$ on $U$
defines an analytic integrable connection
$(\EE^{an},\nabla^{an}:\EE^{an}\to
\EE^{an}\otimes\Omega^{1\ an}_U)$
on $U^{an}$.
The sheaf $\Ker \nabla^{an}$
of flat sections is a local system of
$\CC$-vector spaces on $U^{an}$.
The canonical map
$\Ker \nabla^{an}\to
DR(\EE^{an})
=(\EE\an\otimes
\Omega^{\bullet an}_U)$
is a quasi-isomorphism by
Poincar\'e's lemma.

\proclaim{Lemma 1.6}
A reflexive extension $\EX$
is small if and only if
the analytic de Rham complex
$DR(\EX)\an$
is quasi-isomorphic to
$j_!\an DR(\EE)\an$.
\endproclaim
\demo{Proof}
If $\EX$ is small,
we have
$DR(\EX)\an|_{D\an}=0$
by Lemma 1.5.
We show that, if $\EX$
is not small, the cohomolgy sheaf
$\Cal H^0(DR(\EX)\an)|_{D\an}\ne0$.
Take a maximal small extension
$\EX'\subset \EX$.
By replacing $\EX$
by $\EX\cap\EX'(D_i)$
and shrinking $X$ if necessary,
we may assume that
$\EX$ and $\EX'$ are free,
$D=D_i$
is irreducible and that
$\FF=\EX/\EX'$ is a
free $\Cal O _D$-module.
By the maximality of $\EX'$,
the characteristic polynomial
$\det(T-\res\nabla,\FF)
\in k_i[T]$
has 0 as a solution.
We show that the coherent sheaf
$\EX"=\Ker(\res\nabla:\EX\to\FF)$
is an extension of $(\EE,\nabla)$.
The composition
$\EX\otimes\Omega^1_X(log D)
\overset\nabla\to\to
\EX\otimes\Omega_X^2(\log D)
\overset{\res}\to\to
\EX/\EX'\otimes\Omega^1_D$
is given by $\res\nabla\otimes id$.
Hence on $\EX"$,
we have
$(\res\nabla\otimes id)\circ\nabla=
\res\circ\nabla^2=0$.
This means $\EX"$ is an extension.
By replacing $\EX$
by $\EX"$,
we may further assume that
$\res\nabla$ is 0
on $\Cal F=\EX/\EX'$.
The logarithmic connection $\nabla$
on $\EX$
induces an integrable connection
$\nabla:\Cal F\to \Cal F\otimes\Omega^1_D$.
Then we have exact sequences
of complexes
$$\gather
0\to DR_X(\EX')
\to DR_X(\EX)\to
\Cal F\otimes \Omega^\bullet_X(\log D)\to0\\
0\to DR_D(\Cal F)
\to\Cal F\otimes \Omega^\bullet_X(\log D)\to
DR_D(\Cal F)[-1]\to 0.
\endgather$$
By the assumption $\EX'$ is small we obtain
$$\Cal H^q(DR_X(\EX)\an)\simeq
\cases \Ker\nabla\an\text{ on }\Cal F\an & q=0,1\\
0&\text{otherwise}.
\endcases$$
Thus Lemma 1.6 is proved.
\enddemo

The local system $\Ker \nabla^{an}$
defines a monodromy representation of the fundamental group
$\pi_1(U^{an},\bar x)\to
Aut(\EE(\bar x))$ for $\bar x\in U^{an}$.
For an open subscheme $D_i'$ of an
irreducible component $D_i$ of $D$
where $\EX|_{D_i^{\prime }}$ is locally free,
let $T_i\in Aut(\EX^{an}|_{D_i^{\prime an}})$
be the limit of the positive generator
of the monodromy along $D_i^{an}$.
Then we have
$$T_i=\exp(-2\pi\sqrt{-1}\ \res_i\nabla^{an})$$
by [D1] Chapitre II Proposition 3.11.

Let $F$ be another subfield of $\CC$.
For a smooth separated scheme $U$ over $k_0$,
we define
$M_{k_0,F}(U)$ to be the category
consisting of triples
$\Cal M=((\Cal E,\nabla),V,\rho)$
\roster
\item A locally free $\Cal O_U $-module $\Cal E$
of finite rank with an integrable connection
$\nabla:\Cal E\to \Cal E\otimes\OmU$ which is
regular along the boundary.
\item A local system $V$ of $F$-vector spaces
on the complex manifold
$U^{an}$.
\item A morphism $\rho:V\to \EE\an$
on $U\an$ inducing an isomorphism
$V\otimes_F\Bbb C\isom
Ker \nabla^{an}$ of local systems of $\Bbb C$-vector spaces on $U^{an}$.
\endroster
We call $(\Cal E,\nabla)$
the de Rham component,
$V$ the Betti component and
$\rho$ the comparison isomorphism of
a triple $\MM$.
The category $M_{k_0,F}(U)$
has a natural structure of
neutral tannakian category.
The unit object is
$\bold 1=((\Cal O_U ,d),F, \text{ canonical map})$.
If $U$ is defined over a finite extension
$k$ of $k_0$,
for an element
$p\in (k\otimes_{k_0}\CC)^\times$
an object
$[p]=((\Cal O_U ,d),F, p^{-1}\times\ \text{canonical map})$
is defined.
For an object
$\MM=((\EE,\nabla),V,\rho)$
of $M_{k_0,F}(U)$,
we call $\rank_{\Cal O_X }\EX=\rank_FV$
the rank of $\MM$ and define the determinant
$\det \MM$ to be the triple
$((\det_{\Cal O_X }\EE,\text{trace }\nabla),\det V,\det\rho)$.
The objects of rank 1
of
$M_{k_0,F}(U)$
and the isomorphisms
form a commutative Picard category
$P_{k_0,F}(U)$
with respect to the tensor product.
Following the general terminology [D3] (1.4.1), [D7] (4.1),
a category $P$ is called a commutative Picard
category, if the following condition is satisfied.
The category $P$
is non empty,
every morphism is an isomorphism,
it has a functor
$\otimes:P\times P\to P$,
an associativity constraint
$\otimes\circ(\otimes\times id)\simeq
\otimes\circ(id\times\otimes):
P\times P\times P\to P$
and a commutativity constraint
$\otimes\simeq \otimes \circ c:
P\times P\to P$ compatible with
the associativity constraint
[Sa] I 1.2
and for each object $L$ of $P$,
the functors
$L\otimes\ $
and $\ \otimes L$ are self-equivalence of $P$.
The isomorphism classes of objects of
$P$ form a commutative group
called the class group $Cl(P)$
of $P$.
We write the class group of
$P_{k_0,F}(U)$
by
$MPic_{k_0,F}(U)$.
For a finite extension $k$
of $k_0$,
the class group
$MPic_{k_0,F}(\Spec k)$
is canonically isomorphic to
$k^\times\backslash
(k\otimes_{k_0}\CC)^\times
/(F^\times)^{Hom_{k_0}(k,\CC)}$
by $p\mapsto[p]$.
For an object $\MM=(\EE,V,\rho)$
of
$P_{k_0,F}(k)$
with bases $e$ of $\EE$ and
$v$ of $V$,
its class
$[\MM]$ is the class of
$e/\rho(v)$.
Here we identify
$k\otimes\CC\simeq
\CC^{Hom_{k_0}(k,\CC)}$
by
$x\otimes y\mapsto
(\sigma(x)y)_\sigma.$

For a smooth separated scheme $U$
over a finite extension $k$ of $k_0$ and
$\Cal M=
((\EE,\nabla),V,\rho)
\in M_{k_0,F}(U)$,
we define the determinant of
periods
$$per (\Cal M),
\ per_c (\Cal M)
\in
k^\times\backslash
(k\otimes_{k_0}\CC)^\times
/(F^\times)^{Hom_{k_0}(k,\CC)}$$
as the class
$[\det R\Gamma(U,\Cal M)]$ and
$[\det R\Gamma_c(U,\Cal M)]\in
MPic_{k_0,F}(\Spec k)$
defined below.
Let $X$ be a proper smooth scheme
over $k$ including $U$ as the complement
of a divisor $D$ with simple normal crossings.
For an extension $(\EX,\nabla)$ of $(\EE,\nabla)$
to $X$, we have an isomorphism
$H^q(X,DR(\EX))\otimes_{k_0}\Bbb C\simeq
H^q(X^{an},DR(\EX)^{an})$ by GAGA.
The isomorphism $\rho$ induces
an isomorphism
$H^q(X^{an},DR(\EX)^{an})\simeq
H^q_c(U^{an},V)\otimes_F\Bbb C$
(resp. $H^q(U^{an},V)\otimes_F\Bbb C$)
if $\EX$ is small (resp. big).
Hence we obtain canonical isomorphisms
$H^q_c(\rho):H^q_c(U,DR(\Cal E))\otimes_k\Bbb C\simeq
H^q_c(U^{an},V)\otimes_F\Bbb C$ and
$H^q(\rho):H^q(U,DR(\Cal E))\otimes_k\Bbb C\simeq
H^q(U^{an},V)\otimes_F\Bbb C$.
The triples
$$\align
H^q_c(U,\Cal M)
=&(H^q_c(U,DR(\Cal E)),H^q_c(U^{an},V),H^q_c(\rho))\\
H^q(U,\Cal M)
=&(H^q(U,DR(\Cal E)),H^q(U^{an},V),H^q(\rho))
\endalign$$
are objects of $M_{k_0,F}(\Spec k)$.
Taking the alternating tensor product of the determinant,
$$\align\det R\Gamma_c(U,\Cal M)=&
\otimes_q(\det H^q_c(U,\MM))^{\otimes(-1)^q},\\
\det R\Gamma(U,\Cal M)=&
\otimes_q(\det H^q(U,\MM))^{\otimes(-1)^q}
\endalign$$
$\in M_{k_0,F}(k)$ are defined.
Thus the periods
$$
per_c (\Cal M)=
[\det R\Gamma_c(U,\Cal M)],\quad
per (\Cal M)=
[\det R\Gamma(U,\Cal M)]
$$
$\in k^\times\backslash(k\otimes_{k_0}\Bbb C)^\times/
(F^\times)^{Hom_{k_0}(k,\CC)}$
are defined.
More concretely,
if $k=k_0$,
and, for each $q$, if $e^q$
is a basis of $k$-vector space
$H^q_c(U,DR(\Cal E))$,
$v^q$ is a basis of $F$-vector space
$H^q(U^{an},V)$
and if $P^q$
is the matrix representing $H^q_c(\rho)$
with respect to the bases
$e^q$ and $v^q$,
then
$per_c(\MM)=
\prod_q(\det P^q)^{(-1)^{q+1}}
\in k^\times\backslash\Bbb C^\times/
F^\times$.
Our main theorem
(Theorem 1 in Section 4)
is a formula for
the period
$per_c(\MM)\cdot
per_c(1)^{-\rank \MM}$.
By Serre duality and Poincar\'e duality,
the canonical pairing
$$H^q_c(U,\Cal M)\times
H^q(U,\Cal M^*)\to\bold1(-\dim X)$$
is perfect for $p+q=2\dim X$.
Hence we obtain
\proclaim{Lemma 1.7}
The canonical pairing
$$\det R\Gamma_c(U,\Cal M)\times
\det R\Gamma(U,\Cal M^*)\to\bold 1(-\dim X\cdot\chi_c(U,\MM))$$
is perfect. The equality
$$per_c(\MM)\cdot
per_c(1)^{-\rank \MM}=
(per(\MM^*)\cdot
per(1)^{-\rank \MM})^{-1}$$
holds.
\endproclaim

We define a variant
$\det^\Gamma R\Gamma_c(U,\MM)$
of
$\det R\Gamma_c(U,\MM)$
as follows.
For a field $k$,
let $B_k$
be the image of the homomorphism
$\partial:k(T)^\times
\to k(T)^\times;
f(T)\mapsto f(T)/f(T-1)$.
The map $\partial$
induces an isomorphism
$k(T)^\times/k^\times
\to B_k$.
Let $\Delta:k(T)^\times\to k^\times$
be the homomorphism characterized by
$\Delta f=(-1)^na_na_0^{-1}$
for
$f(T)=\sum_{i=0}^n a_iT^{n-i},
a_na_0\ne0$
and $\Delta T=1$.
It is well-defined since
it is equal to the product of the tame symbols
$(\ ,T)_0\cdot (\ ,T)_\infty:k(T)^\times\to k^\times$.
Since $\Delta|_{k^\times}=1$,
there is a homomorphism
$\Gamma:B_k\to k^\times$
characterized by the property
$\Gamma\circ\partial=1/\Delta$.
Let $\Gamma:\CC(T)^\times\to\CC^\times$
be the homomorphism characterized by
$T-\alpha\mapsto\Gamma(\alpha)$
for $\alpha\notin\Bbb N$,
$T+n\mapsto(-1)^n(n!)^{-1}$ for
$n\in \Bbb N$
and $\Gamma|_{\CC^\times}=1$.
By $\Gamma(\alpha+1)=\alpha\Gamma(\alpha)$,
it is an extension of $\Gamma:B_{\CC}\to\CC^\times$
above.
For $k_0\subset \CC$ and a finite extension $k$ of $k_0$
as above,
we define
$\Gamma:k(T)^\times\to
(k\otimes_{k_0}\CC)^\times
=(\CC^\times)^{Hom_{k_0}(k,\CC)}$
as the composite of the natural map
$k(T)^\times\to(\CC(T)^\times)^{Hom_{k_0}(k,\CC)}$
and $\Gamma$ for $\CC(T)^\times$.
we have a commutative diagram
$$\CD
B_k@>\Gamma>> k^{\times}\\
@VVV@VVV\\
k(T)^\times @>>\Gamma>(k\otimes_{k_0} \Bbb C)^\times.
\endCD$$
For a matrix $M \in M(k,n)$, we define
$\Gamma (M) = \Gamma (\det (T-M)) \in (k \otimes_{k_0}\CC )^{\times}$.
For an extension
$(\EX,\nabla)$
of an integrable connection on $U$,
we put
$$
\Phi_{\EX}(T)=
\prod_{i\in I}
N_{k_i/k}
\Phi_{\EX,i}(T)^{c_i}
$$
$\in k(T)^\times$.
Here $c_i$ is the Euler number
$\chi_c(D_i^*)=
\dim_{k_i}R\Gamma_c(D_i^*,DR(\Cal O _{D_i^*}))$
of $D_i^*=D_i-\bigcup_{j\ne i}D_j$.
When $k$ is a finite extension of
a subfield $k_0$ of $\CC$,
we put
$\Gamma(\nabla:\EX)=
\Gamma (\Phi_{\EX}(T))
\in (k\otimes_{k_0}\CC)^\times$ and
$\Gamma_{D_i} (\nabla:\EX)=
\Gamma (\Phi_{\EX ,i}(T))
\in (k_i\otimes_{k_0}\CC)^\times$.

\proclaim{Lemma 1.8}
For extensions $\EX'\subset \EX$ of
an integrable connection $\EE$,
we have
$$
\Phi_{\EX}(T)/
\Phi_{\EX'}(T)
\in B_k.$$
If $\FF=\EX/\EX'$
is an $\Cal O_{D_i}$-module
for an irreducible component $D_i$ of $D$,
the residue $\res_i\nabla$ induces
an $\Cal O_{D_i}$-linear endomorphism of $\FF$ and
$$\partial\det(T-\res_i\nabla:\FF)=
\Phi_{\EX,i}(T)/
\Phi_{\EX',i}(T)$$
in $k_i(T)$.
\endproclaim
\demo{Proof}
The second assertion
follows from the exact sequence
$$
0\to \FF\otimes\Cal O (-D_i)
(=\Cal Tor_1^{\Cal O_X }(\FF,\Cal O_{D_i}))
\to
\EX'|_{D_i}\to\EX|_{D_i}
\to\FF\to 0.$$
The first assertion is deduced from
the second one by the argument as in
the proof of Lemma 1.5.
\enddemo

We put
$\det (\nabla:\EX/\EX')=
\Gamma(
\Phi_{\EX(T)}/
\Phi_{\EX'(T)})
\in k^\times$.
Let $\MM=((\EE,\nabla),V,\rho)$ be an object of
$M_{k_0,F}(U)$.
Then we define
$\det^\Gamma
R\Gamma_c(U,\MM)$
to be
$$\projlim_{(\EX,\nabla):\text{small}}
(\det R\Gamma(X,DR(\EX)),
\det R\Gamma_c(U\an,V),
\Gamma(\nabla:\EX)^{-1}\times
\det R\Gamma_c(U,\rho)).$$
Here the transition map
$\det R\Gamma(X,DR(\EX'))\to
\det R\Gamma(X,DR(\EX))$
for small extensions $\EX'\subset\EX$ of $\EE$
is
$\det (\nabla,\EX/\EX')\in k^\times$
times the natural map induced by
the quasi-isomorphism
$DR(\EX')\to DR(\EX)$.
We also put
$\det^\Gamma
R\Gamma_c(U,DR(\EE))=
\projlim_{(\EX,\nabla):\text{small}}
\det R\Gamma(X,DR(\EX)).$
We define $per^\Gamma_c(\MM)=
[\det^\Gamma R\Gamma_c(U,\MM)]
\in MPic_{k_0,F}(k)$.
The class
$\Gamma(\nabla:\MM)
\in MPic_{k_0,F}(k)$
of
$\Gamma(\nabla:\EX)
\in (k\otimes_{k_0}\Bbb C)^\times$
is independent of the choice of
a small extension $\EX$.
It is clear that
$$
per^\Gamma_c(\MM)=
per_c(\MM)\cdot\Gamma(\nabla:\MM).
$$

\heading
2. relative Chow group
\endheading

In this section, we define and study the relative Chow
group $CH^n(X\bmod D)$ of dimension 0
and the relative canonical cycle $c_{X\bmod D}
\in CH^n(X\bmod D)$.
They are slight modifications of those in [S1] Section 1.
Let $X$ be a smooth scheme over a field $k$.
We call a finite family
$D=(D_i)_{i\in I}$
of regular subschemes of $X$
a regular family in $X$.
A typical example is the family
of irreducible component of
a divisor $D$ with simple normal crossings.
By abuse of notation, we use the same letter $D$
for the divisor itself and for the family of its components.
Let $\Cal K_n(X)$ denote the sheaf of Quillen's K-group on $X_{Zar}$.
Namely the Zariski sheafification of the presheaf $U\mapsto K_n(U)$ [Q].
For a regular family $D$ in $X$, let
$\Cal K_n(X \bmod D)$
be the complex
$[\Cal K_n(X)\to \oplus_i\Cal K_n(D_i)]$.
Here $\KK_n(X)$ is put on
degree 0 and $\KK_n(D_i)$ denotes their direct image on $X$.
It is the truncation at degree 1
of the complex $\KK_{n,X,D}$ studied in [S1]
and there is a natural map $\KK_{n,X,D}\to\KK_n(X\bmod D)$.
For $n=\dim X$,
we call the hypercohomology
$H^n(X,\KK_n(X\bmod D))$
the relative Chow group of dimension 0 and write
$$CH^n(X \bmod D)=H^n(X,\KK_n(X\bmod D)).$$

We recall the definition of the relative top chern class.
Let $\EE$ be a locally free $\Cal O_X $-module of rank $n$.
We call a family of surjective morphisms
$r_i:\EE|_{D_i}\to\Cal O_{D_i}$ for $i\in I$
a partial trivialization of $\EE$ on $D$.
Let $V=\Bbb V(\EE)=\text{Spec }(S(\EE^*))$
be the covariant vector bundle
associated to $\EE$.
For $i\in I$,
let $\Delta_i=r_i^{-1}(1)$, where
$r_i:V|_{D_i} \to \Bbb A^1_{D_i}$ is
the induced map by
$r_i:\EE|_{D_i}\to\Cal O_{D_i}$
and $1\subset\Bbb A^1$ is the 1-section.
The family $\Delta$ is a regular family in $V$.
Let $\{ 0\} \subset V$
denote the zero section.
Then the horizontal arrows in the diagram below are isomorphisms
$$\CD
H^n_{\{ 0\} }(V, \KK_n(V\bmod\Delta ))
@<<< H^n_{\{ 0\} }(V, \KK_n(V))
@>>> H^0(X ,\Bbb Z)\\
@VVV @. @.\\
H^n(V, \KK_n(V \bmod \Delta ))
@<<< H^n(X, \KK_n(X \bmod D))
@=
CH^n(X\bmod D)
\endCD$$
by the purity and homotopy property of K-cohomology.
The relative top chern class
$c_n(\EE,r)\in CH^n(X \bmod  D)$
is defined as the image of $1\in H^0(X,\Bbb Z)$.
Let $D$ be a divisor with simple normal crossings
and $n=\dim X$.
Then the relative canonical class
$c_{X \bmod D}$
is defined as a relative top chern class by
$$c_{X \bmod D}=
(-1)^nc_n(\Omega^1_X(log D),\res)\in CH^n(X \bmod D).$$
Here the Poincare residue
$\res_i:\Omega^1_X(log D)|_{D_i}\to \Cal O_{D_i}$
defines a partial trivialization
$\res$ of $\Omega^1_X(log D)$ on $D$.

We give some basic properties (cf. [S1] Section 1)
of relative chern classes
used in the proof of main results in Section 4.
For $i\in I$,
the natural map
$\Kn(D_i)[-1]\to
\Kn(X\bmod D)$
induces a map
$i_{i*}:H^{n-1}(D_i,\Kn(D_i))\to
H^n(X,\Kn(X\bmod D))$.

\proclaim{Lemma 2.1}
Let $(\EE,\rho)$
be a locally
free $\Cal O_X $-module of rank $n$
partially trivialized on
a regular family $D$.
For $i\in I$,
let $f_i\in \Gamma(D_i,\Gm)$
and $\sigma=(\sigma_i)_{i\in I}$
be another partial trivialization
of $\EE$ defined by
$\sigma_i=f_i^{-1}\rho_i$.
Then we have
$$c_n(\EE,\sigma)=
c_n(\EE,\rho)-
\sum_{i\in I}
i_{i*}(f_i\cup
c_{n-1}(\EE_i)).$$
Here
$c_{n-1}(\EE_i)\in
H^{n-1}(D_i,\KK_{n-1})$
is the chern class of
the locally free $\Cal O_{D_i}$-module
$\EE_i=\Ker(\rho_i:\EE|_{D_i}\to \Cal O_{D_i})$
and $\cup$ is the cup-product.
\endproclaim

In Lemma 2.2, we consider a partition
$I=I_1\amalg I_2$
of the index
set $I$ of
a regular family $D$.
We put $D_1=(D_i)_{i\in I_1}$
and
$D_2=(D_i)_{i\in I_2}$.
There is a cup-product
$$
H^m(X,\KK_m(X\bmod D_1))\times
H^r(X,\KK_r(X\bmod D_2))\to
H^n(X,\Kn(X\bmod D))$$
for $n=m+r$.
\proclaim{Lemma 2.2}
Let $0\to \Cal F\to\EE\to\Cal G\to 0$
be an exact sequence of locally free
$\Cal O_X $-modules of rank $m,n$ and $r$
respectively, $\rho$
be a partial trivialization of $\EE$
on a regular family $D$.
Assume that on $D_1$,
the partial trivialization $\rho$
induces a partial trivialization
$\sigma=\rho|_{\Cal F}$ of $\Cal F$
and that on $D_2$
it induces $\tau$ of $\Cal G$.
In other word,
for $i\in I_1$,
the restriction $\rho_i:\Cal F\to \Cal O_{D_i}$
is surjective and, for $i\in I_2$
the restriction $\rho_i:\Cal F\to \Cal O_{D_i}$
is 0.
Then we have
$$
c_n(\EE,\rho)=
c_m(\Cal F,\sigma)\cup
c_r(\Cal G,\tau).$$
\endproclaim

\proclaim{Lemma 2.3}
Let $E$ be a regular divisor of $X$ and
let $0\to \Cal E\overset\varphi\to\to\Cal F\to\Cal G\to 0$
be an exact sequence of locally free
$\Cal O_X $-modules $\EE$ and $\Cal F$ of the same rank $n$
and a locally free $\Cal O _E$-module $\Cal G$ of rank $m$.
Let $\varphi:V=\Bbb V(\EE)\to W=\Bbb V(\FF)$
be the associated map,
$j:K\to V$ be the inclusion of
the sub-vector bundle
$K=\Bbb V(\Ker \varphi|_E)\subset V_E$
of rank $m$
and $\Cal N^*_{E/X}$ be the
conormal invertible $\Cal O _E$-module.
Then we have
$$
\varphi^*([0_W])=
[0_V]+j_*(\sum_{k=0}^{m-1}
c_k(\Cal G)c_1(\Cal N^*_{E/X})^{m-1-k})$$
in $H^n(V,\Kn(V\bmod \Delta))$
for a regular family $\Delta\subset V$
such that $\Delta\cap(\{0_V\}\cup K)=\emptyset$.
Here $[0_V]$ is the image of 1 by
$$H^0(\{0_V\},\ZZ)\simeq
H_{\{0_V\}}^n(V,\Kn) \to
H^n(V,\Kn(V\bmod \Delta))$$
and $j_*$ denotes the Gysin map
$$H^{m-1}(K,\KK_{m-1})(\simeq
H^{m-1}(E,\KK_{m-1}))\to
H_K^n(V,\Kn)\to
H^n(V,\Kn(V\bmod \Delta)).$$
\endproclaim

\proclaim{Corollary 1}
Let $\rho$ and $\sigma$
be partial trivializations
on $\EE$ and $\FF$ respectively
on a regular family $D$.
Assume the divisor $E$ does not meet $D$
and $\sigma|_{\EE}=\rho$.
Let $\Cal I=\text{Im }\varphi|_E$
and $i_{E*}:H^{n-1}(E,\KK_{n-1}(E))\to
H^n_E(X,\Kn(X))\to
H^n(X,\Kn(X\bmod D))$
 be the Gysin map.
Then we have
$$c_n(\FF,\sigma)
=c_n(\EE,\rho)
+i_{E*}(c_{n-m}(\Cal I)\sum_{k=0}^{m-1}
c_k(\Cal G)c_1(\Cal N^*_{E/X})^{m-1-k})$$
in $H^n(X,\Kn(X\bmod D))$.
\endproclaim
\demo{Proof}
Clear from the commutative diagram
$$\CD
H^{m-1}(K,\KK_{m-1})
@>>>
H^{n-1}(V_E,\KK_{n-1})
@>>>
H^n(V,\Kn(V\bmod\Delta))\\
@AAA@AAA@AAA\\
H^{m-1}(E,\KK_{m-1})
@>>{\cup c_{n-m}(\Cal I)}>
H^{n-1}(E,\KK_{n-1})
@>>{i_{E*}}>
H^n(X,\Kn(V\bmod D))
\endCD$$
where the vertical arrows are the pull-back and
the upper horizontal arrows are the Gysin maps.
\enddemo

In Corollary 2, we assume
the regular divisor $E$ of $X$
meets transversely with $D_i$
for $i\in I$.
Let $D_E$ be the regular family
$(D_i\cap E)_{i\in I}$ and
let $i_{E*}$ denote
the canonical map
$$\align
i_{E*}:H^{n-1}(E,\KK_{n-1}(E\bmod D_E))\simeq
&\ H^n_E(X,\Kn(X\bmod D))\\
&\to
H^n(X,\Kn(X\bmod D)).
\endalign$$

\proclaim{Corollary 2}
Let $0\to \Cal E\overset\varphi\to\to\Cal F\to\Cal G\to 0$
be an exact sequence of locally free
$\Cal O_X $-modules $\EE$ and $\Cal F$ of the same rank $n$
and an invertible $\Cal O _E$-module $\Cal G$.
Let $\rho$ and $\sigma$
be partial trivialization of $\EE$ and $\Cal F$
on $D$ such that $\rho=\sigma|_{\EE}$.
Then $\sigma$ induces a partial trivialization
$\sigma|_{\Cal I}$ on the
locally free $\Cal O _E$-module
$\Cal I=\text{Image }\varphi|_E$
of rank $n-1$
and we have
$$
c_n(\Cal F,\sigma)
=
c_n(\Cal E,\rho)
+i_{E*}c_{n-1}(\Cal I,\sigma|_{\Cal I}).$$
\endproclaim

In Lemma 2.4, we fix an element $0\in I$
and assume that $E=D_0$
is a regular divisor of $X$
and that $E$ meets transversely with $D_i$
for $i\in I,\ne0$.

\proclaim{Lemma 2.4}
Let $0\to \Cal E\overset\varphi\to\to\Cal F\to\Cal G\to 0$
be an exact sequence of locally free
$\Cal O_X $-modules $\EE$ and $\Cal F$ of the same rank $n$
and locally free $\Cal O _E$-module $\Cal G$ of rank $m$.
Let $\rho$ and $\sigma$
be partial trivialization of $\EE$ and $\Cal F$
on $D$ such that $\rho=\sigma|_{\EE}$.
For $i\in I,\ne0$,
we assume that on $D_i\cup E$
either $\rho_0\oplus\rho_i:
\EE\to\Cal O _{D_i\cap E}^{\oplus 2}$
is surjective or
$\rho_0 |_{D_i\cap E}=\rho_i |_{D_i\cap E}$.
Then we have
$$
c_n(\Cal F,\sigma)
=
c_n(\Cal E,\rho)
.$$
\endproclaim

Proofs of Lemmas 2.1,2.2 and 2.3 and of Corollary 2
are done in the same way
as Propositions 1,2 and 3 and Corollary
loc.cit respectively and we omit them.
Proof of Lemma 2.4 can be simplified considerably and
we give it here.
\demo{Proof of Lemma 2.4}
Let $V=\Bbb V(\EE)$ and
$W=\Bbb V(\Cal F)$
and let $\Delta_V$
and $\Delta_W$
be the regular family of
closed subschemes of $V$
and $W$ defined by partial trivializations
$\rho$ and $\sigma$ respectively.
Since $K \cap\rho^{-1}(1)=\emptyset$, we have a commutative diagram
$$\CD
H^n(W,\Kn(W\bmod \Delta_W))
@>{\varphi^*}>>
H^n(V,\Kn(V\bmod \Delta_V))\\
@AAA @AAA\\
H^n(X,\Kn(X\bmod D))
@=
H^n(X,\Kn(X\bmod D))
\endCD$$
of isomorphisms.
By Lemma 2.3,
it is enough to show that
$j_*:H^{m-1}(K,\KK_{m-1})
\to
H^n(V,\Kn(V\bmod \Delta_V))$
is the 0-map.
We consider a complex
$$
\Kn(V\bmod \Delta)'=
[\Kn(V)\to
\bigoplus_{i\in I}
\Kn(\Delta_i)
\to
\bigoplus_{i\in I,\ne0}
\Kn(\Delta_i\cap \Delta_0)].$$
Since $j_*$ factors
$H^n(V,\Kn(V\bmod \Delta)')$,
it is enough to show
$j_*:H^{m-1}(K,\KK_{m-1})
\to
H^n(V,\Kn(V\bmod \Delta_V)')$
is the 0-map.
It is proved by the same argument as in the last 6 lines of
the first paragraph of the proof of Proposition 4
loc. cit.
\enddemo

In the rest of this section,
we give an adelic presentation
of the form
$$CH^n(X \bmod D)\simeq
\text{Coker }(\partial:\bigoplus_{y\in X_1}B_y
\to\bigoplus_{x\in X_0}A_x).$$
Here $X_i$ denotes the set of the points of $X$ of dimension $i$.

To define the right hand side, we
first introduce some terminology on the tame symbol.
For a vector space $L$
of dimension 1
over a field $K$,
we put $A(L)=
(\bigoplus_{m\in \ZZ}
L^{\otimes m})^\times.$
It is an extension of $\ZZ$
by $K^\times$
since the algebra
$\bigoplus_{m\in \ZZ}
L^{\otimes m}$
is isomorphic to the Laurent
polynomial ring
$K[T,T^{-1}]$.
Its underlying set
is the disjoint union
$\amalg_{m\in \ZZ}
\{K\text{-basis of }L^{\otimes m}\}$.
Let $R$ be a noetherian local integral ring of dimension 1
such that the normalization
$\tilde R$ is finite over $R$.
Let $K$ be the fraction field and
$\kappa$
be the residue field of $R$.
We review the definition of
the order
$\ord:K^\times\to \ZZ$
and the tame symbol
$\partial:
K_2(K)\to \kk^\times.$
For a maximal ideal
$\tilde x$
of $\tilde R$, let
$\ord_{\tilde x}:
K^\times\to \ZZ$ and
$\partial_{\tilde x}:
K_2(K)\to \kk(\tilde x)^\times$
be the usual order
and the tame symbol
$$\partial_{\tilde x}
(\{f,g\})=
(-1)^{\ord_{\tilde x}f
\cdot\ord_{\tilde x}g}
f^{\ord_{\tilde x}g}
g^{-\ord_{\tilde x}f}
(\tilde x)$$
for the discrete valuation ring
$\tilde R_{\tilde x}$.
Then $\ord
=\sum_{\tilde x\mapsto x}
[\kk(\tilde x):\kk]\ord_{\tilde x}$
and
$\partial=
\prod_{\tilde x\mapsto x}
N_{\tilde x/x}\circ \partial_{\tilde x}$.
For an invertible $R$-module
$L$,
we define a pairing
$$(\ ,\ ):
A(L_K)\times K^\times\to
A(L_{\kk})$$
generalizing the tame symbol
as follows.
We put
$A(L)
=\amalg_{m\in \ZZ}
\{R\text{-basis of }L^{\otimes m}\}$
with the structure of abelian group
defined by the tensor product.
The abelian group
$A(L_K)$
is the amalgameted sum
$A(L)\oplus_{R^\times}
K^\times$.
We define the pairing by requiring that,
on $A(L)\times K^\times$,
it is the composition of
(the reduction)$\times\ord:
A(L)\times K^\times\to
A(L_\kk)\times \ZZ$
and the obvious pairing
and that, on
$K^\times\times K^\times$,
it is the tame symbol
$\partial:K^\times\times K^\times
\to \kk^\times$.
Since they coincide on
$R^\times\times K^\times,$
the pairing is well-defined.

The pairing
$(\ ,\ ):
A(L_K)\times K^\times\to
A(L_{\kk})$
is computed on the normalization $\tilde R$ as follows.
Let $\tilde x$ be a maximal ideal of $\tilde R$.
The tame symbol
$(\ ,\ )_{\tilde x}:
A(L_K)\times K^\times\to
A(L_{\kappa(\tilde x)})$
is defined in the same way as above.
We also have the norm map
$A(L_{\kappa(\tilde x)})
\to A(L_\kk)$.
A basis
$\ell$ of
$L_{\kappa(\tilde x)}^{\otimes e}$
is of the form
$a\cdot \ell_0^{\otimes e}$
for $a\in \kappa(\tilde x)^\times$
and a basis $\ell_0$ of $L_\kk$.
The basis $N(\ell)=N_{\kappa(\tilde x)/\kk}(a)
\cdot\ell_0^{\otimes e\cdot[\kappa(\tilde x):\kk]}$
of
$L_\kk^{\otimes e\cdot[\kappa(\tilde x):\kk]}$
is independent of the choice of $\ell_0$
and well-defined.
The norm map $\ell\mapsto N(\ell)$
defines a homomorphism
$N_{\kxt/x}:A(L_{\kappa(\tilde x)})
\to A(L_\kk)$ and the diagram
$$\CD
0@>>>\kxt^\times
@>>> A(L_{\kxt})
@>>>\ZZ
@>>>0\\
@.
@V{N_{\kxt/x}}VV
@V{N_{\kxt/x}}VV
@VV{[\kxt:\kk]}V
@.\\
0@>>>\kk^\times
@>>> A(L_{\kk})
@>>>\ZZ
@>>>0
\endCD$$
is commutative.
We have
$$(\ ,\ )_\kk=
\prod_{\tilde x}N_{\kxt/\kk}\circ(\ ,\, )_{\tilde x}.$$

We define $A_x$ for $x\in X$.
For $x\in X$, we put $I_x=\{i;x\in D_i\}$.
For $i\in I_x$, let $N_i(x)$ be the one-dimensional
$\kx$-vector space $\Cal O_X (-D_i)\otimes\kx$.
We put $A_{x,i}=A(N_i(x))$
and write the canonical map
$\ord_i:A_{x,i}\to \ZZ$.
The abelian group
$A_x$ is the fiber product
$\left(\prod_{i\in I_x}\right)_{\ZZ} \Axi$
with respect to $\ord_i$.
It is an extention of $\ZZ$ by
$\bigoplus_{i\in I_x}\kxx$.
We define $B_y$ for $y\in X$.
For $i\in I_y$,
let $\Byi$
be the amalgamated sum
$(A_{y,i}\otimes\kyx)\oplus_
{\kyx\otimes\kyx}
K_2(y)$.
It is an extension of
$\kyx$ by $K_2(y)$.
The group $B_y$ is defined as the fiber product
$\left(\prod_{i\in I_y}\right)_{\kyx}
\Byi$ with respect to $\ord \otimes id :\Byi\to\kyx$.

The homomorphism $\partial$ is the direct sum
of the $(x,y)$-components
$\partial_{x,y}:B_y\to A_x$
for $x\in X_0$ and $y \in X_1$.
This fits in the commutative diagram
$$\CD
0@>>>
\bigoplus_{i\in I_y}K_2(\ky)@>>>
 B_y@>>>
 \kyx
@>>> 0\\
@.@V{\bigoplus_{i\in I_y}(\ ,\ )_x}VV
@VV{\partial_{x,y}}V
@VV{\ord_x}V\\
0@>>>
\bigoplus_{i\in I_x} \kxx @>>>
A_x@>>>
\Bbb Z
@>>> 0
\endCD$$
and is 0
unless $x$ is not in the closure $Y$ of $\{y\}$.
Here $\ord_x:\kyx\to \Bbb Z$ is the usual order
and $(\ ,\ )_x:K_2(\ky)\to \kxx$ is the tame symbol.
The local component
$\partial_{x,y}:B_y\to A_x$
is determined by its $i$-th components
$\partial_{x,y,i}:B_{y,i}\to A_{x,i}$
for $i\in I_y$ and
$\partial_{x,y,i}:\kyx\to A_{x,i}$
for $i\in I_x- I_y$.
We use the tame symbol $(\ ,\ )$ defined above for the
local ring $\Cal O _{Y,x}$
and the invertible module
$\Cal O (-D_i)$.
For $i\in I_y,$
the $i$-th component
$\partial_{x,y,i}$
is induced by the pairing
$(\ ,\ ):
A(N_i(y))\times\kyx
\to
A(N_i(x))=A_{x,i}$.
For
$i\in I_x- I_y$,
it is defined as
$(-1)^{\ord_x}\times(1,\ ):
\kyx\to
A(N_i(x))=A_{x,i}$.
Here $(1 ,\ )$ is the pairing with
1 regarded as a
$\ky$-basis of
$\Cal O (-D_i)\otimes \ky$
which is an element of
$A(\Cal O (-D_i)\otimes \ky)$.
It is clear that
we have a commutative diagram
$$\CD
0@>>>
 K_2(\ky)@>>>
 B_{y,i}@>>>
 \kyx
@>>> 0\\
@.@V{\partial_x}VV
@VV{\partial_{x,y,i}}V
@VV{\ord_x}V\\
0@>>>
 \kxx @>>>
A_{x,i}@>>>
\Bbb Z
@>>> 0
\endCD$$
for $i\in I_y$
and that the composition
$\kyx\overset{\partial_{x,y,i}}\to\to
\Axi\overset{\ord_i}\to\to\ZZ$
is $\ord_x$ for
$i\in I_x- I_y$.
Hence $\partial_{x,y}:
B_y\to A_x$
is defined
and satisfies the commutative diagram above.

We give an explicit computation of
$\partial_{x,y}$.
For an inverse image
$\tilde x$ of $x$ in the normalization
$\tilde Y$
of the closure $Y$ of $\{y\}$,
let
$A_{\tilde x,i}
=A(N_i(x)\otimes_{\kx}\kxt)$
and define the fiber product
$A_{\tilde x}
={\prod_{i\in I_x}}_{\ZZ}
A_{\tilde x,i}$
and the boundary map
$\partial_{\tilde x,y}:
B_y\to A_{\tilde x}$
in the same way as above.
The norm maps defined above induces
$N_{\kxt/\kx}
:A_{\tilde x}\to A_x$
and we have
$$\partial_{x,y}=
\prod_{\tilde x}
N_{\kxt/\kx}
\partial_{\tilde x,y}.$$
We compute
$\partial_{\tilde x,y}$
explicitely
by taking a basis
$\pi_i$ of
$\Cal O (-D_i)|_Y$ at $x$ for $i\in I_x$.
The basis define splittings
$A_{\tilde x}=
\ZZ\oplus
\bigoplus_{i\in I_x}\kxt^\times$ and
$B_y=\kyx\oplus
\bigoplus_{i\in I_y}K_2(y)$ by
$1\mapsto \pi_i$.
Since the other components of $\partial_{\tilde x,y}$
are evident,
we compute the component
$\kyx\to\kxt^{\times}$.
If $i\in I_y$, it is 0
and if $i\in I_x-I_y$,
it is given by
$$f\mapsto
(-1)^{\ord_{\tilde x}f}
(\pi_i,f)_{\tilde x}^{-1}
=(-\pi_i,f)_{\tilde x}^{-1}.$$
In fact in the latter case,
the rational section 1
of $\Cal O (-D_i)|_Y$ is
$\pi_i^{-1}$-times the basis $\pi_i$.

\proclaim{Proposition 1}
Let $X$ be a smooth
scheme of dimension $n$
over a field $k$
and $D$ be a divisor with simple
normal crossings.
Under the notation above,
there exists a canonical isomorphism
$$CH^n(X \bmod D)
\simeq
\text{Coker }(\partial:\bigoplus_{y\in X_1}B_y
\to
\bigoplus_{x\in X_0}A_x).$$
\endproclaim

\demo{Proof}
We consider the spectral sequence
$$E_1^{p,q}=\bigoplus_{x\in X^{p}}
H^{p+q}_x(X,\KK_n(X\bmod D))\Rightarrow
H^{p+q}(X,\KK_n(X\bmod D)).$$
By the long exact sequence
$$
\to
\bigoplus_i
H_x^{r-1}(X,\Kn(D_i))
\to
H_x^r(X,\Kn(X\bmod D))
\to
H_x^r(X,\Kn(X))
\to
$$
and by the purity
$$\align
H_x^r(X,\Kn(X))=&
\cases
K_{n-r}(x)&x\in X^r\\
0&otherwise
\endcases\\
H_x^{r-1}(X,\Kn(D_i))=&
\cases
K_{n-r+1}(x)&x\in X^r\cap D_i\\
0&otherwise,
\endcases
\endalign$$
we have
$E_1^{p,q}=0$
for $q\ne0$
and
an exact sequence
$$0\to
\bigoplus_{x\in X^p}
(\bigoplus_{i\in I_x}
K_{n-p+1}(x))
\to
E_1^{p,0}
\to
\bigoplus_{x\in X^p}
K_{n-p}(x)
\to 0.$$
Therefore the spectral sequence
degenerates at $E_2$-terms and
we have an exact sequence
$$\align
\bigoplus_{y \in X_1}
H^{n-1}_y(X,\Kn(X\bmod D))
&\to
\bigoplus_{x\in X_0}
 H^n_x(X,\KK_n(X\bmod D))\\
\to &CH^n(X\bmod D)
\to 0.
\endalign$$
Hence it is sufficient to define isomorphisms
$\alpha_x:\Ax\simeq
H^n_x=H^n_x(X,\KK_n(X\bmod D))$
for $x\in X_0$
and
$\beta_y:B_y\simeq
H^{n-1}_y=H^{n-1}_y(X,\Kn(X\bmod D))$
for $y\in X_1$
such that the diagram
$$\CD
B_y       @>{\partial_{x,y}}>>       A_x\\
@V{\beta_y}VV                 @VV{\alpha_x}V\\
H^{n-1}_y @>>{d_{x,y}}>              H^n_x
\endCD$$
is commutative.

First we prove
\proclaim{Lemma 2.5}
Let $X$ be a smooth scheme
over a field $k$
and $i:D\to X$ be the closed immersion
of a smooth divisor.
Then there is a canonical pairing
$$
R^1i^!\KK_p(X\bmod D)
\times
\KK_q(D)
\to
R^1i^!\KK_{p+q}(X\bmod D).$$
It is compatible with the
natural pairing
$$
R^1i^!\KK_p(X\bmod D)
\times
i^*\KK_q(X)
\to
R^1i^!\KK_{p+q}(X\bmod D)$$
induced by the cup-product.
\endproclaim
\demo{Proof}
It is sufficient to define a functorial pairing
$$
H^1_{D\cap U}(U,\KK_p(X\bmod D))
\times
H^0(D\cap U,\KK_q(D))
\to
H^1_{D\cap U}(U,\KK_{p+q}(X\bmod D))
$$
for open subsets $U\subset X$
such that the diagram
$$\CD
H^1_{D\cap U}(U,\KK_p(X\bmod D))
&\times&
H^0(D\cap U,\KK_q(D))
@>>>
H^1_{D\cap U}(U,\KK_{p+q}(X\bmod D))\\
@|  @AAA @|\\
H^1_{D\cap U}(U,\KK_p(X\bmod D))
&\times&
H^0(U,\KK_q(X))
@>>>
H^1_{D\cap U}(U,\KK_{p+q}(X\bmod D))
\endCD$$
is commutative.
Since it is easy to check the functoriality,
we write $U=X$ for simplicity.

We consider the deformation to the normal
bundle of $D$ in $X$.
Let $\tilde X$ be the complement of the proper
transform of
$X\times\{0\}$
in the blowing-up
of
$\Bbb A^1_X$
at the center
$D\times\{0\}$.
Let $V=\Bbb V(N_{D/X})\subset \tilde X$
be the exceptional divisor
and $\tilde D
\simeq \Bbb A^1_D$
be the proper transform of
$\Bbb A^1_D\subset \Bbb A^1_X$.
We obtain a commutative diagram
$$\matrix
H^1_D(V,\KK_p(V\bmod D))
&\times&
H^0(V,\KK_q(V))
&@>>>&
H^1_D(V,\KK_{p+q}(V\bmod D))\\
@AAA @AAA  @AAA\\
H^1_{\tilde D}(\tilde X,\KK_p(\tilde X\bmod\tilde D))
&\times&
H^0(\tilde X,\KK_q(\tilde X))
&@>>>&
H^1_{\tilde D}(\tilde X,\KK_{p+q}(\tilde X\bmod\tilde D))\\
@VVV   @VVV       @VVV\\
H^1_D(X,\KK_p(X\bmod D))
&\times&
H^0(X,\KK_q(X))
&@>>>&
H^1_D(X,\KK_{p+q}(X\bmod D)).
\endmatrix$$
Here the upper vertical arrows are
induced by the inclusion $V\to \tilde X$ at 0
and the lower ones are induced by
$X\to \tilde X$ at 1.

We show that the vertical arrows in
the left and right columns are
isomorphisms.
By the purity of $K$-groups
$R^qi^!\KK_p(X)=\KK_{p-1}(D)$
for $q=1$ and $=0$ for otherwise,
we have an exact sequence
$$
0\to
\KK_p(D)
\to
R^1i^!\KK_p(X\bmod D)
\to
\KK_{p-1}(D)
\to
0$$
and
$R^qi^!\KK_p(X\bmod D)=0$
for $q\ne1$.
By the spectral sequence
$$E_2^{p,q}=
H^p(D,R^qi^!\KK_p(X\bmod D))
\Rightarrow
H^{p+q}_D(X,\KK_p(X\bmod D)),$$
we have an exact sequence
$$
\align
0\to
H^0(D,\KK_p(D))
\to
H^1_D(X,\KK_p(X\bmod D))
& \to
H^0(D,\KK_{p-1}(D)) \\
& \to
H^1(D,\KK_{p}(D)). \\
\endalign $$
We have similar exact sequences for
$\tilde D
\subset
\tilde X$
and $D\subset V$.
Now by the homotopy property of
$K$-cohomology,
the vertical arrows in
the left and right columns are
isomorphisms.
Further by
the homotopy property,
the upper middle
$H^0(V,\KK_q(V))$
is canonically isomorphic to
$H^0(D,\KK_q(D))$.
By these isomorphisms,
the upper row gives the pairing
$$
H^1_D(X,\KK_p(X\bmod D))
\times
H^0(D,\KK_q(D))
\to
H^1_D(X,\KK_{p+q}(X\bmod D)).
$$

We check the compatibility.
The lower arrow
$H^0(\tilde X,\KK_q(\tilde X))
\to
H^0(X,\KK_q(X))$
in the middle column
is a surjection
since it is induced by a section of the
canonical map
$\tilde X\to X$.
 Therefore it is sufficient to prove the commutativity
of the diagram
$$\CD
H^0(\tilde X,\KK_q(\tilde X))
@>>>
H^0(X,\KK_q(X))\\
@VVV  @VVV\\
H^0(V,\KK_q(V))
@>>>
H^0(D,\KK_q(D)).
\endCD$$
This follows from the commutative diagram
$$\matrix
\tilde X &     &  @<<<   &     & X   \\
         &     &         &     & @AAA\\
@AAA          \tilde D &@<<0<& D   \\
         &     &  @AA1A        &     \\
  V      &@<<< & D       &     & .
\endmatrix$$
Here the arrows
with tags
0 and 1 are the 0-section and 1-section
respectively
and they induces the same isomorphisms
$H^0(\tilde D,\KK_q(\tilde D))\to
H^0(D,\KK_q(D)).$
Thus Lemma 2.5 is proved.
\enddemo
\enddemo

We define the isomorphisms
$\alpha_x$ and $\beta_y$
using Lemma 2.5.
First we consider
$\alpha_x:
A_x\to
H^n_x(X,\Kn(X\bmod D))$
for
$x\in X_0$.
It is sufficient to define
an isomorphism
$$\alpha_{x,i}:\Axi\to
H^n_x(X,\Kn(X\bmod D_i))=
H^n_{x,i}$$
such that the diagram
$$\CD
\Axi              @>>>     \ZZ  \\
@V{\alpha_{x,i}}VV         @|      \\
H^n_{x,i}         @>>> H^n_x(X,\Kn(X))
\endCD$$
is commutative for each $i\in I_x$.
Hence we may assume $D=D_i$ is irreducible and
$x\in D$.
By Lemma 2.5, we obtain a pairing
$$\align
H^0(D_x,
R^1i^!\KK_1(X\bmod D))
& \times
H^{n-1}_x(D,\KK_{n-1}(D))  \\
& \to
H^{n-1}_x(D_x ,
R^1i^!\Kn(X\bmod D)).\\ \endalign $$
We show that
$$\gather
\Gamma(X_x-D_x,\Gm)/(1+\Cal I_x)\simeq
H^0(D_x,
R^1i^!\KK_1(X\bmod D)),\\
\ZZ\simeq
H^{n-1}_x(D,\KK_{n-1}(D)),\\
H^n_x(X_x,
\Kn(X\bmod D))\simeq
H^{n-1}_x(D_x,
R^1i^!\Kn(X\bmod D))
\endgather$$
for $\Cal I=\Cal O (-D)$.
Under the natural identification
$A_x\simeq
\Gamma(X_x-D_x,\Gm)/1+m_x$,
we also show that the pairing induces
the required map
$\alpha_x:
A_x\to
H^n_x(\Kn(X\bmod D))$.

By purity, we have
$H^{n-1}_x(D,\KK_{n-1}(D))=\ZZ$.
Since
$R^qi^!\Kn(X\bmod D)=0$
for $q\ne1$,
the spectral sequences
$$E_2^{p,q}=
H^p(D_x,
R^qi^!\KK_1(X\bmod D))
\Rightarrow
H^{p+q}_{D_x}
(X_x,\KK_1(X\bmod D))$$
$$E_2^{p,q}=
H^p_x(D_x,
R^qi^!\Kn(X\bmod D))
\Rightarrow
H^{p+q}_x
(X,\Kn(X\bmod D))$$
degenerate at $E_2$
and we have isomorphisms
$$\align
H^0(D_x,
R^1i^!\KK_1(X\bmod D))
\simeq
&H^1_{D_x}
(X_x,\KK_1(X\bmod D))\\
H^{n-1}_x(D_x,
R^1i^!\Kn(X\bmod D))
\simeq
&H^n_x
(X,\Kn(X\bmod D)).
\endalign$$
Further
by the long exact sequence
$$
\to
H^p(X_x,\KK_1(X\bmod D))
\to
H^p(X_x-D_x,\KK_1(X))
\to
H^{p+1}_{D_x}(X_x,\KK_1(X\bmod D))
\to
$$
and by the equalities
$H^0(X_x,\KK_1(X\bmod D))
=1+\Cal I_x$
for $\Cal I=\Cal O (-D)$,
$H^0(X_x-D_x,\KK_1(X))
=\Gamma(X_x-D_x,\Bbb G_m)$
and by
$H^1(X_x,\KK_1(X\bmod D))
=0$,
we have
$H^1_{D_x}(X_x,\allowbreak\KK_1(X\bmod D))
=\Gamma(X_x-D_x,\Bbb G_m)
/1+\Cal I_x.$
Hence the pairing above is a
homomorphism
$\tilde \alpha_x:
\Gamma(X_x-D_x,\Bbb G_m)
/1+\Cal I_x\to
H^n_x
(X,\Kn(X\bmod D)).$

By the exact sequence
$$0\to
\KK_q(D)\to
R^1i^!\KK_q(X\bmod D)
\to
\KK_{q-1}(D)
\to0,
$$
the map
$\tilde \alpha_x$
sits in the commutative diagram
$$\CD
0
@>>>
\Gamma(D_x,\Bbb G_m)
@>>>
\Gamma(X_x-D_x,\Bbb G_m)
/1+\Cal I_x
@>>>
\ZZ
@>>>0\\
@. @VVV @V{\tilde\alpha_x}VV  @| @.\\
0
@>>>
\kxx
@>>>
H^n_x(X,\Kn(X\bmod D))
@>>>
\ZZ
@>>>0.
\endCD$$
Therefore
$\tilde\alpha_x$
annihilates $1+m_x$.
Since
$A_x\simeq
\Gamma(X_x-D_x,\Bbb G_m)
/1+ m_x$,
the map $\tilde\alpha_x$
induces the required isomorphism
$\alpha_x:A_x\to H^n_x.$

Next we consider
$\beta_y$
for $y\in X_1$.
Similarly as above,
it is sufficient to define
an isomorphism
$\beta_y:
B_y\to
H^{n-1}_y
(X,\Kn(X\bmod D))$
assuming $D$ is irreducible and $y\in D$
such that the diagram
$$\CD
B_y
@>>>
\kyx\\
@V{\beta_y}VV @VV{\wr}V\\
H^{n-1}_y
(X,\Kn(X\bmod D))
@>>>
H^{n-1}_y
(X,\Kn(X))
\endCD$$
is commutative.
By Lemma 2.5,
we obtain a pairing
$$ \align
H^0(D_y,
R^1i^!\KK_1(X\bmod D))
& \times
H^{n-2}_y(D,\KK_{n-1}(D)) \\
& \to
H^{n-2}_y(D_y,
R^1i^!\Kn(X\bmod D)).\\ \endalign $$
Similarly as above,
we have isomorphisms
$H^0(D_y,
R^1i^!\KK_1(X\bmod D))
=\Gamma(X_y-D_y,\Bbb G_m)
/1+\Cal I_y$,
$H^{n-2}_y(D,\KK_{n-1}(D))
=\kyx$
and
$H^{n-2}_y(D_y,
R^1i^!\Kn(X\bmod D))
=H^{n-1}_y
(X,\Kn(X\bmod D))$.
The same argument as above
shows that the pairing
induces a morphism
$\tilde\beta_y:
A_y\otimes\kyx
\to
H^{n-1}_y
(X,\Kn(X\bmod D))$
and further an isomorphism
$\beta_y:B_y\to
H^{n-1}_y
(X,\Kn(X\bmod D))$.

Finally we show that the boundary map
$d=d_1^{n-1,0}$
of the spectral sequence
coincides with
$\partial:
\bigoplus B_y\to
\bigoplus A_x$
under the identification by $\alpha$ and $\beta$.
It is sufficient to show that
the $(x,y)$-component
$d_{x,y}$ and $\partial_{x,y}$
are the same for
each $x\in X_0$ and $y\in X_1$.
If $x$
is not in the closure $Y$ of
$\{y\}$,
we have
$d_{x,y}=\partial_{x,y}=0$.
Hence we assume
$x\in Y$.
If $x\notin D$,
then $y\notin D$ and
$d_{x,y}=\partial_{x,y}=
\ord_x:\kyx\to\ZZ$.
Hence we assume $x\in D$.
By considering each
$i$-th component for $i\in I_x$,
we may assume $D$ is irreducible.
First we check the case where
$Y$ meets transversally $D$ at $x$.
Then by Lemma 2.5,
we obtain a commutative diagram
$$\CD
H^0(X_x-D_x,\Bbb G_m)
\times
H^{n-1}_Y(X,\KK_{n-1}(X))
@>>>
H^{n-1}_y(X,\Kn(X))\\
@VVV @VV{d_{x,y}}V\\
H^1_{D_x}(X_x,\KK_1(X\bmod D))
\times
H^{n-1}_x(D,\KK_{n-1}(D))
@>>>
H^n_x(X,\Kn(X\bmod D))
\endCD$$
since
$$\align
H^1_{D_x}(X_x,\KK_1(X\bmod D))
=&
H^0(D_x,R^1i^!\KK_1(X\bmod D))\\
H^n_x(X,\Kn(X\bmod D))
=&
H^{n-1}_x(D_x,R^1i^!\KK_n(X\bmod D)).
\endalign$$
Hence the diagram
$$\CD
\Gamma(X_x-D_x,\Bbb G_m)
@>>>
\kyx\\
@VVV @VV{d_{x,y}}V\\
A_x @>{\sim}>{\alpha_x}>
H^n_x
\endCD$$
is commutative.
Since the diagram
$$\CD
\Gamma(X_x-D_x,\Bbb G_m)
&@>>>&
\kyx\\
@VVV& \swarrow\partial_{x,y}&@.\\
A_x& & @.
\endCD$$
is commutative by definition
and since
$\Gamma(X_x-D_x,\Bbb G_m)
\to\kyx$
is surjective,
the differentials
$d_{x,y}$
and $\partial_{x,y}$
are the same.

We consider the general case.
First we show that
the problem is etale local.
Let
$\varphi:
X'\to X$
be an etale morphism
$\varphi(x')=x$
and
$\eta_{x'}=
\{y'\in X';\varphi(y')=y\text{ and }
\overline {\{y'\}}\ni x' \}$.
Then we have commutative diagrams
$$\CD
H^{n-1}_y
@>{d_{x,y}}>>
H^n_x
@.\quad
B_y
@>{\partial_{x,y}}>>
A_x\\
@VVV @VVV @VVV @VVV\\
\bigoplus_{y'\in\eta_{x'}}
H^{n-1}_{y'}
@>{\bigoplus d_{x',y'}}>>
H^n_{x'}
@.\quad
\bigoplus_{y'\in\eta_{x'}}
B_{y'}
@>{\bigoplus \partial_{x',y'}}>>
A_{x'}.
\endCD$$
Since
$A_x\to A_{x'}$ is injective,
to prove
$d_{x,y}=\partial_{x,y}$,
it is enough to show
$d_{x',y'}=\partial_{x',y'}$
for $y'\in\eta_{x'}$.

\proclaim{Claim 1}
By taking an etale neighborhood,
we find a closed integral subschemes
$W\subset Z\subset X$ satisfying the following conditions
\roster
\item $x\in W$ and $y\in Z$.
\item $W$ is regular of dimension 1
and meets transversally with $D$ at $x$.
\item $Z$ is of dimension 2.
$Z-W$ is regular and meets transversally with $D$.
\item The normalization $\tilde Z$
is regular and the
inverse image
$\tilde T$ of $D$ in $\tilde Z$
is regular.
\item The reduced inverse image $\tilde W$ of $W$
in $\tilde Z$ is regular and
$\tilde W$ meets transversally
with $\tilde T$
\item $\tilde T\cap\tilde W$
=(inverse image of $x$ in $\tilde Z$)
contains only one point $\tilde x$.
\endroster
\endproclaim
\demo{Proof of Claim 1}
By localizing if necessary,
we may take a smooth projection
$\pi:X\to D$
such that the immersion $D\to X$
is a section.
Since it is trivial if $\pi(y)=x$,
we assume $\pi(Y)=T$ is a curve.
By further localizing if necessary,
we may assume that
$T-\{x\}$ is normal and that
the inverse image of $x$
in the normalization $\tilde T$
of $T$ consists of only one point
$\tilde x$.
Then $W=\pi^*(x)$ and $Z=\pi^*(T)$
satisfy the conditions.
\enddemo

Changing the notation,
we use the same characters $X,Z,W$ etc. in Claim 1
for their inverse images in the localization
$X_x$ at $x$.
We consider the canonical morphism
$$d:
H^{n-1}_{Z-x}(X,\Kn(X\bmod D))
\to
H^n_x(X,\Kn(X\bmod D)).$$
Let $\zeta$ be the
generic point of $Z$ and $Z^*=Z-x$.
By the spectral sequence
$E_1^{p,q}=
\bigoplus_{z\in X^p\cap Z^*}
H^{p+q}_z(X,\Kn(X\bmod D))
\Rightarrow
H^{p+q}_{Z^*}(X,\Kn(X\bmod D))$
degenerating at $E_2$-terms,
we have an exact sequence
$$
K_2(\kk(\zeta))\to
\bigoplus_{z\in Z_1}
H^{n-1}_z(X,\Kn(X\bmod D))
\to
H^{n-1}_{Z^*}(X,\Kn(X\bmod D))
\to 0.$$
To complete the proof of Proposition 1,
it is sufficient to prove the following

\proclaim{Claim 2}
(1). By identifying by $\beta:
\bigoplus_{z\in Z_1}
B_z
\to
\bigoplus_{z\in Z_1}
H^{n-1}_z(X,\Kn(X\bmod D))$,
the map
$\bigoplus\partial_{x,z}:
\bigoplus_{z\in Z_1}
B_z
\to A_x$
annihilates the image of
$K_2(\kk(\zeta))\to
\bigoplus_{z\in Z_1}
\allowbreak H^{n-1}_z(X,\Kn(X\bmod D)).$

(2). For the generic point $w\in Z_1$ of $W$,
the canonical map
$H^{n-1}_w(X,\Kn(X\bmod D))
\to
H^{n-1}_{Z^*}(X,\Kn(X\bmod D))$
is surjective.
\endproclaim

We show that Claim 2 proves Proposition 1.
In fact Claim 2 (1) implies
$\bigoplus\partial_{x,z}$
induces a map
$\partial:
H^{n-1}_{Z^*}(X,\Kn(X\bmod D))
\to A_x$.
By identifying
by $\alpha_x:A_x\to
H^n_x(X,\Kn(X\bmod D))$,
we see $d=\partial$
on the image of
$H^{n-1}_w(X,\Kn(X\bmod D))$
since $W$ meets transversally $D$ and we have
$d_{x,w}=\partial_{x,w}$.
By Claim 2 (2),
we have
$d=\partial$ on the whole
$H^{n-1}_{Z^*}(X,\Kn(X\bmod D))$.
This proves $d_{x,z}=\partial_{x,z}$
for all $z\in Z_1$ in particular for $y$.

\demo{Proof of Claim 2}
(1). We take an arbitrary element
$a\in K_2(\zeta)$
and compute its image by
$$K_2(\zeta)\to
\bigoplus_{z\in Z_1}B_z
\to A_x.$$
By the condition (6) in Claim 1,
the intersection
$T=D\cap Z$
is integral.
Let $t$ be the generic point of $T$.
For $z\in Z_1$,
we have
$B_z=\kk(z)^\times$
for $z\ne t$.
We take a prime element
$\pi\in
\Cal O _{X,x}$ of $D$.
Then it defines splittings
$A_x=\ZZ\oplus\kxx$
and
$B_t=\kk(t)^\times \oplus K_2(t)$.
For
$z\in Z_1,\ne t$,
let $a_z=\partial_{z,\zeta}(a)
\in \kk(z)^\times$ and
for $z=t$, let
$(a^1_t,a^2_t)
=\partial_{t,\zeta}(a)\in
\kk(t)^\times\oplus K_2(t)
=B_t$.
The first component
$a^1_t\in \kk(t)^\times$ is the usual
tame symbol of $a\in K_2(\zeta)$.
Then by definition of
$\partial$,
we have
$$
\partial_{x,z}\circ\partial_{z,\zeta}(a)=
\cases
(\ord_x(a_z),
(-1)^{\ord_x(a_z)}(\pi,a_z)_x)
& z\ne t\\
(\ord_x(a^1_t),
\partial a^2_t)
& z=t.
\endcases$$
The second component
$\partial a^2_t\in \kxx$
is the usual tame symbol of
$a^2_t\in K_2(t).$
Therefore it is sufficient to prove
$$\align
\sum_{z\in Z_1,\ne t}
\ord_x(a_z)+
\ord_x(a^1_t)=&0\\
\prod_{z\in Z_1,\ne t}
(-1)^{\ord_x(a_z)}
(\pi,a_z)_x\times
\partial(a^2_t)=&1.
\endalign$$

The first equality is the usual reciprocity law for
$a\in K_2(\zeta)$. By the first one,
the second is equivalent to
$$\prod_{z\in Z_1,\ne t}
(\pi,a_z)_x\times
(-1)^{\ord_x(a^1_t)}
\partial(a^2_t)=1.
$$
To prove this,
we use the reciprocity law for
$\{\pi,a\}\in K_3(\zeta)$.
Since
$\partial_{z,\zeta}
(\{\pi,a\})=
\{\pi,a_z\}\in K_2(z)$
for $z\ne t$,
it is sufficient to show that
$$
\partial_{x,t}\circ
\partial_{t,\zeta}
(\{\pi,a\})=
(-1)^{\ord_x(a^1_t)}
\partial a^2_t.$$
in $\kxx$.
We prove this by showing
$$
\partial_{t,\zeta}
(\{\pi,a\})=
\{-1,a^1_t\}+a_t^2
$$
in $K_2(t)$.
Let $\tilde a_1
\in \Gamma(X_t,\Bbb G_m)$
be a lifting of
$a^1_t\in \kk(t)^\times.$
Then $\tilde a_2=a-\{\pi,\tilde a_1\}
\in K_2(\zeta)$
is in the kernel of the tame symbol
$K_2(\zeta)\to
\kk(t)^\times$
and hence in
$K_2(\Cal O _{Z,t})$.
By the commutativity of the diagram
$$\CD
H^1_{D_t}(X_t,\KK_1(X\bmod D))
\times
\kk(t)^\times
@>>>
H^{n-1}_t(X,\Kn(X\bmod D))\\
@AAA @AAA\\
\Gamma(X_t-D_t,\Bbb G_m)
\times
\Gamma(X_t,\Bbb G_m)
@>>>
K_2(\zeta),
\endCD$$
we have
$\partial_{t,\zeta}
(\{\pi,\tilde a_1\})=
(a^1_t,0)$.
By the commutative diagram
$$\CD
K_2(\Cal O _{Z,t})
&@>>>&
K_2(\zeta)\\
@VVV& &@VVV\\
K_2(t)
&@>>>&
H^{n-1}_t(X,\Kn(X\bmod D)),
\endCD$$
$a^2_t$ is the reduction of $\tilde a_2$.
Now we have
$$
\{\pi,a\}=
\{\pi,\pi,\tilde a_1\}+
\{\pi,\tilde a_2\}
=\{\pi,(\{-1,\tilde a_1\}+
\tilde a_2)\}.$$
Therefore we have
$$\partial_{t,\zeta}
(\{\pi,a\})=
\{-1,a_t^1\}+
a_t^2.$$
Thus the assertion (1) is proved.

(2). By the exact sequence
$$
\align
H^{n-1}_w(X,\Kn(X\bmod D)) &\to
H^{n-1}_{Z-x}(X,\Kn(X\bmod D)) \\
&\to
H^{n-1}_{Z-W}(X,\Kn(X\bmod D)), \\
\endalign
$$
and by the purity
$$
H^{n-1}_{Z-W}(X,\Kn(X\bmod D))
\simeq
H^1(Z-W,\KK_2(Z-W\bmod T-x)),
$$
it is sufficient to show that
$H^1(Z-W,\KK_2(Z-W\bmod T-x))=0$.
By the exact sequence
$$
\align
H^1(\tilde Z,\KK_2(\tilde Z\bmod\tilde T)) &\to
H^1(Z-W,\KK_2(Z-W\bmod T-x)) \\
&\to
H^2_{\tilde W}(\tilde Z,\KK_2(\tilde Z\bmod\tilde T)) \\
\endalign
$$
and by the purity
$$H^2_{\tilde W}(\tilde Z,\KK_2(\tilde Z\bmod\tilde T))
\simeq
H^1(\tilde W,\KK_1(\tilde W\bmod \tilde x)),$$
it is reduced to show that
$$\align
H^1(\tilde Z,\KK_2(\tilde Z\bmod\tilde T))=&
\text{Coker }(K_2(\Cal O _{\tilde Z})
\to(K_2(\Cal O _{\tilde T})))=0\\
H^1(\tilde W,\KK_1(\tilde W\bmod \tilde x))=&
\text{Coker }(\Cal O _{\tilde W}^\times)
\to\kk(\tilde x)^\times=0.
\endalign$$
Since
$K_2(\Cal O _{\tilde T})$
is generated by symbols,
the maps are surjective.
Thus Claim 2 and hence
Proposition 1 is proved.
\enddemo

\proclaim{Corollary}
Let $U'$ be a dense open subscheme of $U$.
The relative Chow group
$CH^n(X \bmod D)$
is generated by the classes
$[x]$ for $x\in U'_0$.
\endproclaim
\demo{Proof of Corollary}
First we assume $U'=U$.
Let $U_m=
\{x\in X;\Card I_x\le m\}$.
We have $X=U_n$ and $U=U_0$.
By induction on $m$,
it is sufficient to show that
the image of
$\bigoplus_{x\in U_m}A_x$
in $CH_0(X,D)$ is in that of
$\bigoplus_{x\in U_{m-1}}A_x$
for $m\ge 1$.
Let $x\in X_0$ such that $\Card I_x=m\ge 1$.
For each $i\in I_x$,
we put $I_i=I_x-\{i\}$ and
take a closed integral curve
$Y_i\subset D_{I_i},\ni x$
regular at $x$ and meeting $D_i$
transversally at $x$.
We also take a basis $\pi_i$
of $N_i(x)=\Cal O (-D_i)(x)$
for each $i\in I_x$.
Let $y_i$
be the generic point of $Y_i$
and $S_i\subset Y_i$
be the closed subset
$S_i=\bigcup_{j\in I-I_i}(Y_i\cap D_j)$
of $Y_i$.
Let $f$ be an arbitrary element of $A_x$.
We take a rational function
$g_i\in\kappa(y_i)^\times$
for $i\in I_x$
and a rational section
$\pi_{i,j}$
of $\Cal O (-D_j)|_{Y_i}$
for $i\ne j,\in I_x$
satisfying the following conditions.
\roster
\item"{}"
The boundary
$\partial_{x,y,i}(g_i)$
is the $i$-component
$f_i\in A_{x,i}\simeq
\kappa (y_i)^{\times}/1+m_x$
and the pull-back of $g_i$ to
the normalization of $Y_i$
is invertible
and congruent to 1 at the inverse image of
 $S_i-\{x\}$.
The section $\pi_{i,j}$
is a basis of
$\Cal O (-D_j)|_{Y_i}$
at $S_i$
and the fiber
$\pi_{i,j}(x)$
is $\pi_j$ for $j\ne i$.
\endroster
Let $h_i\in B_{y_i}$ be the element
defined by
$(\pi_{i,j})_j\otimes g_i$.
Then it is easy to check that
we have
$f=\sum_{i\in I_x}
\partial_{x,y_i}(h_i)$
and
$\partial(\bigoplus_{i\in I_x}h_i)$
is supported on
$\{x\}\cup
\bigcup_i(Y_i-S_i)$.
Since $Y_i-S_i=Y_i\cap U_{m-1}$,
Corollary for $U'=U$ is proved.

To prove Corollary for general $U'$,
it is sufficient to show that the image of
$\bigoplus_{x\in U}\ZZ$
in $CH_0(X,D)$ is in that of
$\bigoplus_{x\in U'}\ZZ$.
For a closed point $x\in U-U'$,
take a closed integral curve $Y\ni x$
regular at $x$ and meeting $U$
and $y$ be the generic point of $Y$.
Let $f$ be an arbitrary integer.
Take a rational function
$g\in\kyx=B_y$ such that $\ord_xg=f$
and that the pull-back of $g$
to the normalization of $Y$ is invertible
and congruent to 1
at the inverse image of $Y-U'$.
Then we have $\partial_{x,y}(g)=f$
and $\partial(g)$
is supported on $\{x\}\cup(Y\cap U')$.
Thus Corollary is proved.
\enddemo


\heading
3. Tame symbol.
\endheading

Let $k_0,F\subset \Bbb C$
be subfields and
$k$ be a finite extension of $k_0$.
Let $X$ be a proper smooth scheme over $k$
and $U$ be the complement of
a divisor $D$ of $X$ with simple normal
crossings as in section 1.
Recall that $MPic_{k_0,F}(U)$
denotes the class group of
$P_{k_0,F}(U)$.
The purpose of this section
is to define a canonical pairing
$$\align
(\ ,\ ):MPic_{k_0,F}(U) \otimes CH^n(X \bmod D) \to
&MPic_{k_0,F}(k)\\
&\simeq
k^\times\backslash(k\otimes_{k_0}\Bbb C)^{\times}/
(F^{\times})^{Hom_{k_0}(k,\CC)}.
\endalign$$

To state a more precise result,
we prepare some terminology on commutative Picard categories
[D3] (1.4.2)-(1.4.11).
Let $P$ and $P'$ be commutative Picard categories.
An additive functor $P\to P'$
of commutative Picard categories
is defined to be a pair $(f,\varphi)$
of a functor $f:P\to P'$
and an isomorphism
$\varphi:
\otimes\circ(f\times f)\to
f\circ\otimes$
of functors $P\times P\to P'$
compatible with the
associativity constraint $a$ and
the commutativity constraint $c$ as follows.
For objects $L,L'$ and $L^{\prime\prime}$ of $P$,
the diagrams
$$\CD
(f(L)\otimes f(L'))\otimes f(L^{\prime\prime})
@>>{a_{f(L),f(L'),f(L^{\prime\prime})}}>
f(L)\otimes (f(L')\otimes f(L^{\prime\prime}))\\
@V{\varphi_{L,L'}\otimes id_{f(L^{\prime\prime})}}VV
@VV{id_{f(L)}\otimes\varphi_{L',L^{\prime\prime}}}V\\
(f(L\otimes L'))\otimes f(L^{\prime\prime})
@.
f(L)\otimes (f(L'\otimes L^{\prime\prime}))\\
@V\varphi_{L\otimes L',L^{\prime\prime}}VV
@VV{\varphi_{L,L'\otimes L^{\prime\prime}}}V\\
f((L\otimes L')\otimes L^{\prime\prime})
@>{f(a_{L,L',L^{\prime\prime}})}>>
f(L\otimes (L'\otimes L^{\prime\prime}))
\endCD$$
and
$$\CD
f(L)\otimes f(L')
@>>{c_{f(L),f(L')}}>
f(L')\otimes f(L')\\
@V{\varphi_{L,L'}}VV
@VV{\varphi_{L,L'}}V\\
f(L\otimes L')
@>{f(c_{L,L'})}>>
f(L'\otimes L)
\endCD$$
are commutative.
The pair $1=(1,1)$
of the functor
$1:P\to P'$
sending every object to $1_{P'}$
and every morphism to $id_{1_{P'}}$
and of the isomorphism
$1:1\circ\otimes\to\otimes\circ (1\times 1)$
defined by the canonical map
$1_{P'}\to 1_{P'}\otimes 1_{P'}$
is an additive functor.
An isomorphism of additive functors
$(f,\varphi)\to(g,\psi)$
is an isomorphism of functors
$f\to g$
such that the diagram
$$\CD
f(L)\otimes f(L')@>>> f(L\otimes L')\\
@VVV @VVV\\
g(L)\otimes g(L')@>>> g(L\otimes L')
\endCD$$
is commutative.

Let $P,P'$ and $P_1$
be commutative Picard categories.
We call a biadditive functor (1.4.8) loc.cit. a
pairing. We recall the definition.
A pairing $P\times P_1\to P'$
of commutative Picard categories
is defined to be a triple
$(f,\varphi,\varphi_1)$
of a functor
$f:P\times P_1\to P'$
and isomorphisms of functors
$\varphi:
\otimes_{P'}\circ
((f\circ pr_{13})\times
(f\circ pr_{23}))\to
f\circ (\otimes_P\times 1_{P_1}):
P\times P\times P_1\to P'$
and
$\varphi_1:
\otimes_{P'}\circ
((f\circ pr_{12})\times
(f\circ pr_{13}))\to
f\circ (1_P\times \otimes_{P_1}):
P\times P_1\times P_1\to P'$
satisfying the following compatibilities.
For an object
$L_1$ of $P_1$
the pair $(f_{L_1}=f(\ ,L_1),
\varphi_{L_1}=\varphi(\ ,L_1))$
is an additive functor $P\to P'$
and
similarly for an object $L$ of $P$.
Further for objects $L,L'$ of $P$
and $L_1,L'_1$ of $P'$,
the diagram
$$\CD
(f(L,L_1)\otimes
f(L,L'_1))\otimes
(f(L',L_1)\otimes
f(L',L'_1))
&\to&
f(L,L_1\otimes L'_1)\otimes
f(L',L_1\otimes L'_1)
\\
@VVV& &\\
(f(L,L_1)\otimes
f(L',L_1))\otimes
(f(L,L'_1)\otimes
f(L',L'_1))
& &
@VV{\varphi(L,L',L_1\otimes L'_1)}V\\
@V{\varphi(L,L',L_1)\otimes
\varphi(L, L', L'_1)}VV
& &
\\
f(L\otimes L',L_1)\otimes
f(L\otimes L',L'_1)
&\to &
f(L\otimes L',L_1\otimes L'_1)
\endCD$$
is commutative.
Similarly as above,
the unit pairing
$1:P\times P_1\to P'$ is defined.
A pairing $P\times P_1\to P'$
induces a pairing
$Cl(P)\times Cl(P_1)
\to Cl(P')$
of class groups.
An isomorphism of pairings
$(f,\varphi,\varphi_1)\to(g,\psi,\psi_1)$
is an isomorphism of functors
$f\to g$
such that
$(f_{L_1},\varphi_{L_1})\to
(g_{L_1},\psi_{L_1})$
is an isomorphism of additive functor
for an object $L_1$ of $P_1$
and similaly for an object $L$ of $P$.

We consider the case where
$P_1$ is defined by a morphism
$\partial: B\to A$
of abelian groups as follows (1.4.11) loc.cit..
Let $P_\partial$
be the category with the set $A$
of objects and
$Hom_{P_\partial}(a,a')=
\{b\in B;a'=a\partial(b)\}$
for $a,a'\in A$.
Its class group $Cl(P_\partial)$
is equal to the cokernel
$C=\text{Coker }(\partial: B\to A)$.
When $B=0$, we write $P_\partial=A$ by abuse of notation.
When $\partial:B\to A$
is injective,
the natural functor
$P_{\partial}\to\text{Coker }\partial$ is
an additive equivalence.

\proclaim{Lemma 3.1}
Let $P$ and $P'$ be commutatitive Picard categories and
$\partial :B\to A$
be a homomorphism of abelian groups.

1. A pairing
$P\times A\to P'$
is the same thing as a pair
$((f_a)_{a\in A},(\varphi_{a,a'})_{a,a'\in A})$
of a family of additive functors and
and a family of isomorphisms of functors
$\varphi_{a,a'}:
f_a\otimes f_{a'}\to
f_{aa'}$
such that the diagrams below are commutative.
$$\CD
(f_a\otimes f_{a'})\otimes f_{a^{\prime\prime}}
@>>>
f_{aa'}\otimes f_{a^{\prime\prime}}
@>>>
f_{aa'a^{\prime\prime}}\\
@VVV @. @|\\
f_a\otimes( f_{a'}\otimes f_{a^{\prime\prime}})
@>>>
f_{a}\otimes f_{a'a^{\prime\prime}}
@>>>
f_{aa'a^{\prime\prime}}
\endCD$$
$$\CD
f_{a}\otimes f_{a'}
@>>>
f_{aa'}\\
@VVV  @|\\
f_{a'}\otimes f_{a}
@>>>
f_{a'a}
\endCD$$
$$\CD
(f_a(L)\otimes f_{a'}(L))
\otimes
(f_a(L')\otimes f_{a'}(L'))
@>>>
f_{aa'}(L)
\otimes
f_{aa'}(L')\\
@VVV@.\\
(f_a(L)\otimes f_{a}(L'))
\otimes
(f_{a'}(L)\otimes f_{a'}(L'))
@.\downarrow\\
@VVV@.\\
f_a(L\otimes L')
\otimes
f_{a'}(L\otimes L')
@>>>
f_{aa'}(L\otimes L').
\endCD$$
An isomorphism of pairings
$((f_a),(\varphi_{a,a'}))\to
((g_a),(\psi_{a,a'}))$
is a family of isomorphisms
$f_a\to g_a$
of additive functors such that
the diagram
$$\CD
f_a\otimes f_{a'}
@>>>
f_{aa'}\\
@VVV@VVV\\
g_a\otimes g_{a'}
@>>>
g_{aa'}
\endCD$$
is commutative for $a,a'\in A$.

2. Let $f:P\times A\to P'$ be a pairing.
Then to define a pairing
$P\times P_\partial \to P'$
inducing $f$ by the functor $A\to P_\partial$
is equivalent to to give an isomorphism of the pairings
$1\to
f\circ(1_P\times\partial)$
from the unit pairing $1$
to the composition
$f\circ(1_P\times\partial):
P\times B\to P'$.
By the assertion 1, it is further equivalent to give
a family $(e_b)_{b\in B}$ of trivializations
$e_b:1\to f_{\partial b}$
of additive functors $P\to P'$
satisfying the commutative diagram
$$\CD
1\otimes 1
@>{e_b,e_{b'}}>>
f_{\partial(b)}\otimes
f_{\partial(b')}\\
@VVV @VV{\varphi_{\partial(b),\partial(b')}}V\\
1
@>>{e_{bb'}}>
f_{\partial(bb')}.
\endCD$$
\endproclaim
\demo{Proof}
1. For a pairing $(f,\varphi,\varphi_1):
P\times A\to P'$,
the family of functors
$f_a=(f(\ ,a),\varphi(\ ,\ ,a)):P\to P'$
for $a\in A$
and the family of isomorphisms
$\varphi_{a,a'}=\varphi_1(\ ,a,a')$
satisfy the compatibility above.
It is easy to check the converse.

2. Let
$1\to
f\circ(1_P\times\partial)$
be a trivialization.
For $b\in Hom_{P_\partial}(a,a')$, $a, a'\in A$,
we define an isomorphism $f(\ ,b):f_a\to f_{a'}$
by $$f_a=f_a\otimes 1
@>{1_{f_a}\otimes e_b}>>
f_a\otimes f_{\partial(b)}
@>{\varphi(a,\partial(b))}>>
f_{a\partial(b)}=f_{a'}.$$
It is easy to check that together with
$(f,\varphi,\varphi_1):P\times A\to P'$,
it defines a pairing $P\times P_\partial\to P'$.
The rest is easy and omitted.
\enddemo

\proclaim{Corollary}
For an exact sequence of abelian group
$0\to B@>{\partial}>>A\to C\to 0$,
to define a pairing
$P\times C\to P'$
is equivalent to to give
a pairing $P\times A\to P'$
and a trivialization of the induced pairing
$P\times B\to P'$.
\endproclaim
\demo{Proof}
By Lemma 3.1.2 above,
the latter is equivalent to give
a pairing $P\times P_\partial\to P'$.
Hence by the equivalence $P_\partial \to C$,
it is proved.
\enddemo

Let $k_0,k,F,X$ and $D$ as in the beginning of the section.
We put $\bold{CH}^n(X\bmod D)=P_\partial$
for the map $\partial:\bigoplus_{y\in X_1}B_y
\to\bigoplus_{x\in X_0}A_x$ in the last section.
By Proposition 1, we have
$Cl(\bold{CH}^n(X\bmod D))=
\text{Coker }(\partial: B\to A)
=CH^n(X\bmod D)$.
We define a pairing
$$P_{k_0,F}(U)\times \bold{CH}^n(X\bmod D)
\to P_{k_0,F}(k).$$
It induces a pairing
$MPic_{k_0,F}(U)\times CH^n(X\bmod D)
\to MPic_{k_0,F}(k)$
of class groups.

Before defining the pairings,
we prepare generalities on the tame symbol and the norm.
First we define tame symbol.
We assume that $X$ is a regular scheme and
$D$ is a divisor with simple normal crossings.
Let $I$ be the index set of the
irreducible components of $D=\bigcup_{i\in I}D_i$.
We put $A_i=
\Gamma(X-D_i,\Bbb G_m)$.
There is an exact sequence
$$0\to
\Gamma(X,\Bbb G_m)
\to A_i @>{\ord_i}>> \ZZ.$$
For a subset $J\subset I$, let $A_J$
denote the fiber product
$\left(\prod_{i\in J}\right)_{\ZZ}A_i$
with respect to the order
$\ord_i:A_i\to \ZZ.$
We write the canonical map
$A_J\to \ZZ$ by $\ord_J$.
We will define pairings
of the form
$$P(U)\times A_J
\to P(D_J^*)$$
called tame symbol
for various Picard categories
$P(U)$ on $U=X-D$
and $P(D_J^*)$
on $D_J^*=
\bigcap_{i\in J}D_i$.

\demo{Remark}
We explain the philosophy behind the definition.
This remark of heuristic argument
has logically independent of the other part of the paper.
Suppose the class group $ClP(U)$ is
$H^1(U,G)$.
We regard $A_I$ as a rough approximation of
$H^{2d}_{D_I}(X,j_!\ZZ(d))$
where $d=\Card I$
and $j:U\to X$ is the open immersion.
Then the pairing
$ClP(U)\times A_I
\to ClP(D_I)$
is induced by the cup-product
$$H^1(U,G)\times
H^{2d}_{D_I}(X,j_!\ZZ(d))
\to
H^{2d+1}_{D_I}(X,G(d))$$
and the Gysin isomorphism
$H^{2d+1}_{D_I}(X,G(d))\simeq
H^1(D_I,G)=P(D_I).$
When $d=1$, we put
$j_!\ZZ(1)=\Ker(\Bbb G_{m,X}
\to \Bbb G_{m,D})[1]$.
Then we get
a map $A\to
H^2_D(X,j_!\ZZ(1))$.
For general $d$,
by cup-product,
we get a map
$\bigotimes_{i\in I}
A_i\to
H^{2d}_{D_I}(X,j_!\ZZ(d)).$
The fiber product $A_I$
approximates the quotient
of the tensor product
$\bigotimes_{i\in I}
A_i$
by the subgroup
$\langle\otimes a_i;
\ord_i(a_i)=0
\text{ for at least 2 $i$'s}\rangle.$

There is another way to explain.
We consider the case $J=I$.
As in [D5] (1.7.10),
we should have a functor
$P(U)\to P(N_I^*)$.
Here $N_I^*=N_I-\cup_{i\in I}N_i$
is the complement
in the normal bundle $N_I$ of
$D_I$ in $X$ of the normal bundles $N_i$
$D_I$ in $D_i$ for $i\in I$.
An element $f$ in $A_I$
defines a section
$\bar f:D_I\to N_I^*$.
Then the pairing
$(\ ,f):P(U)\to P(D_I)$
is the composition of the pull-back by $\bar f$
with $P(U)\to P(N_I^*)$.
\enddemo

First we define the tame symbol for invertible
modules.
Let $P(U)$ and $P(D_I)$
be the Picard categories of invertible
$\Cal O_U $-modules and $\Cal O _{D_I}$-modules respectively
and $A=A_I$ be the abelian group defined above.
We define a pairing
$$(\ ,\ ):P(U)\times A\to P(D_I)$$
as follows.
Let $\EE$ be an
invertible $\Cal O_U $-module and $f\in A$.
We define an invertible
$\Cal O _{D_I}$-module
$(\EE,f)$ as follows.
For an open subset $S\subset X$
and a basis $e$ of $\EE$ on $U\cap S$,
$(\EE,f)$ is generated by
the symbol $(e,f)$ on $D_I\cap S$.
For another basis $e'=ge$
and $g\in\Gamma(U\cap S,\Cal O )^{\times}$,
we put
$$(g,f)=
(-1)^{\sum_i\ord_ig\cdot\ord f}
g^{\ord f}
\prod_if_i^{-\ord_ig}
|_{D_I}$$
$\in \Gamma (D_I\cap S,\Bbb G_m)$
and we impose $(e',f)=(g,f)(e,f)$.
It is easy to check that the invertible sheaf
$(\EE,f)$
is well-defined and a pairing
$P(U)\times A\to P(D_I)$
is defined.
For an invertible
$\Cal O_X $-module $\EX$ extending $\EE$,
there is a canonical isomorphism
$(\EE,f)\to\EX^{\otimes\ord f}|_{D_I}$
defined by $(e,f)\mapsto
e^{\otimes \ord f}$
for a local basis $e$ of $\EX$.

Next we define a pairing for integrable
connections of rank 1.
Let $k$ be a field of characteristic 0
and $X$
be a smooth scheme over $k$.
Let $D=\bigcup_{i\in I}
D_i$ be a divisor with simple normal crossings and
$U$ be the complement.
Let $P^\nabla_X(U)$
be the Picard category
of integrable connections
$(\EE,\nabla:\EE\to \EE\otimes\OmU)$
of rank 1 regular along $D$.
For a local basis $e$ of
an invertible
$\Cal O_U $-module $\EE$
for an object $(\EE,\nabla)$
of $P^\nabla_X(U)$,
let $\nabla \log e$
be the 1-form satisfying
$\nabla e=e\otimes\nabla\log e$.
It has at most logarithmic pole along $D$ and
its residue
$\res_i(\nabla\log e)\in \Cal O_{D_i}$ is defined
for each $i\in I$.
The connection $\nabla$ on $\EE$
is determined by $\nabla\log e$.
Let $J\subset I$
be a subset and
write $T=D_J=\bigcap_{i\in J}D_i$.
Let $D_T=\bigcup_{i\in I-J}(T\cap D_i)$
be a divisor with simple normal crossings
of $T$ and write $T^*=T-D_T$.
We define a pairing
$$(\ ,\ )_T:
P^\nabla_X(U)\times A_J
\to
P^\nabla_T(T^*).$$
For an object $(\EE,\nabla)$
of $P^\nabla_X(U)$
and $f=(f_i)_{i\in J}\in A=A_J$,
we define a logarithmic integrable
connection
$\nabla:
(\EE,f)\to
(\EE,f)\otimes
\Omega^1_{T^*}$
on the invertible
$\Cal O _{T^*}$-module
$(\EE,f)$ regular along
$D_T$ as follows.
For a basis $e$
of $\EE$ on $S\cap U$
for an open subset $S\subset X$
as above, we put
$$\nabla \log(e,f)=
\ord f\cdot\nabla \log e-
\sum_{i\in J}
\res_i(\nabla\log e)\cdot
d\log f_i.$$
The right hand side is a section
of
$\Omega^1_X(log D)|_T$
on $S\cap T$ and is in
$\Ker(
\bigoplus\res_i:
\Omega^1_X(log D)|_T
\to\bigoplus_{i\in J}\Cal O _T)
=\Omega^1_T(\log D_T)$.
It is easy to see that the
connection on $(\EE,f)$ is well-defined
and is regular along $D_T$
and that it defines a pairing
$
P^\nabla_X(U)\times A_J
\to
P^\nabla_T(T^*).$
Let $\EX$ be an invertible $\Cal O_X $-module
extending $\EE$ and
assume $\res_i\nabla\in
\Gamma(D_i,\End_{\Cal O_{D_i}}(\EX|_{D_i}))
=\Gamma(D_i,\Cal O_{D_i})$
is in the image of a constant
$\nabla_i\in\Gamma(X,\Cal O_X )$
such that $d(\nabla_i)=0$ for $i\in J$.
We define a logarithmic integrable connection
$\nabla^\sim:\EX^{\otimes \ord f}\to
\EX^{\otimes\ord f}\otimes\Omega^1_X(\log D')$
where $D'=\bigcup_{i\in I-J}D_i$
by
$$\nabla^\sim\log (e^{\otimes\ord f})
=\ord f\cdot\nabla \log e-
\sum_{i\in J}
\nabla_i\cdot
d\log f_i$$
for a local basis $e$ of $\EX$.
It is easily seen to be well-defined
and the tame symbol
$((\EE,f),\nabla)$ is isomorphic to the restriction
$(\EX^{\otimes \ord f},\nabla^\sim)|_T.$

Finally we define a pairing
$$
P_{k_0,F,X}(U)\times A_J
\to
P_{k_0,F,T}(T^*).$$
Let $k_0$ and $F$
are subfields of $\CC$,
let $X$ be a smooth
scheme over a finite extension $k$ of $k_0$
and $U$ be the complement of a divisor $D$
with simple normal crossings
as in the last paragraph.
We do not assume $X$ is proper.
Let $M_{k_0,F,X}(U)$ be the category
defined quite similary as
$M_{k_0,F}(U)$ defined in Section 1
except that we replace the condition
that the connection
$(\EE,\nabla)$ is regular along the boundary
by that it is regular along $D$.
The rank 1 objects of
$M_{k_0,F,X}(U)$
form a commutative Picard category
$P_{k_0,F,X}(U)$.

Let $\MM=((\EE,\nabla),V,\rho)$
be an object of $P_{k_0,F,X}(U)$.
For $f=(f_i)_{i\in J}\in A=A_J$,
we define the pairing
$(\MM,f)=(((\EE,f),\nabla),(V,f),(\rho,f))$
which is an object of $P_{k_0,F,T}(T^*)$
as follows.
The integrable connection
$((\EE,f),\nabla)$
has been defined above.
We define a local system
$(V,f)$ of $F$-vector spaces
on $(T^*)^{an}$ as follows.
Let $S\subset X$
be an open subset and
$v$ be a multivalued section of
$V$ on $(U\cap S)^{an}$.
We define a horizontal section
$(v,f)$ of
$(\EE,f)^{an}$
on $(T^*\cap S)^{an}$
well-defined up to a multiplicative constant
in the image of
$\pi_1(U\an\cap(X\an)_i)$ in $F^{\times}$
on $(T\cap S)\an\cap (X\an)_i$
for each connected component
$(X\an)_i$ of $X\an.$
Shrinking $S$, if necessary, we take a basis
$e$ of $\EE$
on $S$ and let $\varphi=\rho(v)/e$
be a (multivalued) analytic function on
$(U\cap S)\an$.
We put
$$\align
&(-1)^{-\ord f\cdot\res\nabla\log e}f^{\res\nabla\log  e}\\
=&
\exp\left(\sum_{i\in J}\res_i(\nabla \log e)\cdot
(-\ord f\cdot\pi\sqrt{-1}+\log f_i)\right)
\endalign$$
and
$$(\varphi,f)^\sim=
(-1)^{-\ord f\cdot\res\nabla\log e}f^{\res\nabla\log  e}\times
\varphi^{\ord f}.$$
Here $\res_i(\nabla\log e)$
is regarded as a locally constant function on
$S\an$.
They are well-defined upto a multiplicative constant
in the image of
$\pi_1(U\an\cap(X\an)_i)$ in $F^{\times}$
since
$\exp(2\pi\sqrt{-1}\cdot\nabla_i)$
is the inverse of the eigenvalue
of the monodromy operater along
$D_i^{an}$.
We show that the function
$(\varphi,f)^\sim$
is an invertible holomorphic function
on $(S-\bigcup_{i\in I-J}D_i)\an$.
Since $\rho(v)$ is a horizontal section,
we have
$\nabla\log e=-d\log \varphi$.
Hence the logarithmic differential
$$d\log(\varphi,f)^\sim=
\sum_{i\in J}\res_i(\nabla \log e)\cdot
d\log f_i+\ord f\cdot d\log \varphi$$
does not have pole along $D_i$ for $i\in J$
and is holomorphic
on $(S-\bigcup_{i\in I-J}D_i)\an$.
Therefore
$\log(\varphi,f)^\sim$
is a holomorphic function and
$(\varphi,f)^\sim$ is
invertible and holomorphic.

We define the tame symbol
$(\varphi,f)$
to be the restriction of
$(\varphi,f)^\sim$ on $(T^*)\an$.
It is an analogue of the tame symbol
$(g,f)=
((-1)^{\sum_i\ord_i g\cdot \ord f}
g^{\ord f}
\prod_if_i^{-\ord_i g})|_{T^*}$
for a rational function $g$
since $d\log \varphi=-\nabla\log e$
and $\ord_ig=\res_id\log g$.
We define
$(v,f)=(\varphi,f)(e,f)$.
It is easily seen to be independent of
the choice of a basis $e$.
We define a local system
$(V,f)$
of $F$-vector spaces
of rank 1
to be that generated by
$(v,f)$.
Since $(v,f)$
is well-defined upto a multiplicative constant in
$F^\times$ on each connected component,
the local system $(V,f)$ is well-defined.
The comparison
$(\rho,f):(V,f)\to(\EE,f)\an$
is the natural inclusion.
For an object
$\MM=((\EE,\nabla),V,\rho)$
of $P_{k_0,F,X}(U)$,
the tame symbol
$(\MM,f)$ is defined as the triple
$(\MM,f)=(((\EE,f),\nabla),(V,f),(\rho,f))$
in $P_{k_0,F,T}(T^*)$,
It is easy to check that
the tame symbol defines a pairing
$P_{k_0,F,X}(U)
\times A_J\to
P_{k_0,F,T}(T^*).$

If we take an invertible $\Cal O_X $-module $\EX$ extending $\EE$,
the definition of the tame symbol is rephrased as follows.
We put
$$(-1)^{-\ord f\cdot\nabla}f^\nabla=
\prod_{i\in J}
(-1)^{-\ord f\cdot\nabla_i}
f_i^{\nabla_i}=
\exp(\sum_{i\in J}\nabla_i\cdot
(-\ord f\cdot\pi\sqrt{-1}+\log f_i))$$
where $\nabla_i=\res_i(\nabla:\EX)\in
End_{\Cal O_{D_i}}(\EX|_{D_i})=
\Gamma(D_i,\Cal O_{D_i})$ similarly as above.
The subsheaf
$$(V,f)^\sim=
(-1)^{-\ord f\cdot\nabla}f^\nabla
\rho(V)^{\otimes\ord f}
\subset
\EE^{an\otimes\ord f}$$
is well-defined and is a local system of
$F$-vector spaces on
$U^{an}$ of rank 1.
Since
$\CC\otimes_F(V,f)$
is the sheaf $\text{Ker }(\nabla^\sim)\an$
of horizontal sections of $(\EX^{\otimes\ord f},\nabla^\sim)$,
it is unramified along
$D_i$ for $i\in J$
and extended to a local system
on $(X-\bigcup_{i\in I-J}D_i)\an$
also denoted by $(V,f)^\sim$.
The isomorphism
$(\EE,f)\to\EX^{\otimes\ord f}$
induces an isomorphism
$(V,f)$ to the restriction
$(V,f)^\sim|_{T^*}$.

There is a variant of the tame symbol.
Let $\alpha:Y\to T$
be a morphism
from a smooth $k$-scheme $Y$.
We assume that the reduced inverse image
$D_Y=(\alpha^*D_T)_{red}$
is
a divisor with simple normal crossing.
We put $Y^*=Y-D_Y$.
For $i\in J$,
let $A_{Y,i}$
be the disjoint union
$\amalg_{m\in \ZZ}
\{\text{basis of }N_i|^{\otimes m}_Y\}$
with a structure of abelian group
defined by tensor product.
We have an exact sequence
$$0\to
\Gamma(Y,\Gm)
\to A_{Y,i}
\overset{\ord_i}\to\to
\ZZ.$$
Let $A_{Y,J}$
be the fiber product
${\prod_{i\in J}}_{\ZZ}
A_{Y,i}$
with respect to
$\ord_i:A_{Y,i}\to\ZZ$.
There is a natural map
$\alpha^*:A_J\to
A_{Y,J}$.
Locally on $Y$,
the group $A_{Y,J}$
is the amalgamated sum of
$A_J$ and
$\Gamma(Y,\Gm)^J
=\Ker(\ord:A_{Y,J}\to \ZZ)$
over
$\Gamma(X,\Gm)^J
=\Ker(\ord:A_J\to \ZZ)$.

We define a pairing
$$(\ ,\ )_Y:
P_{k_0,F,X}(U)
\times A_{Y,J}\to
P_{k_0,F,Y}(Y^*)$$
also called the tame symbol
as follows.
The composition of $\alpha^*$ and the
pairing $(\ ,\ )_T$ defines a pairing
$$P_{k_0,F,X}(U)
\times A_J\to
P_{k_0,F,Y}(Y^*).$$
For $u=(u_i)_i
\in \Gamma(Y,\Gm)^J$
and an object
$\MM=((\EE,\nabla),
V,\rho)$
of $P_{k_0,F,X}(U)$,
we imitate the definition above as follows.
For a basis $e$ of $\EE$
on an open subset $S\subset X$,
let $(\EE,u)$
be the invertible $\Cal O_Y $-module
generated by the symbol $(e,u)$
on $\alpha^{-1}(S)$.
For an invertible function
$g\in \Gamma (U\cap S,\Gm)$
and another basis $e'=ge$, we define the tame symbol
$(g,u)=\prod_{i\in J}
u_i^{-\ord_i g}$
and require $(e',u)=(g,u)(e,u)$.
We define
a connection
$\nabla$ on
$(\EE,u)$
by
$\nabla\log(e,u)=
-\sum_{i\in J}
\alpha^*(\res_i\nabla)u_i^{-1}du_i$.
We define a local system
$(V,u)$ of $F$-vector spaces
to be that generated by
$(v,u)=
\prod_{i\in J}
u_i^{\res_i\nabla}\cdot
(e,u)
\subset
(\EE,u)^{an}$
where
$u_i^{\res_i\nabla}
=
\exp(\res_i\nabla\cdot\log u_i)$.
The comparison is the inclusion.
Then we see that a pairing
$$P_{k_0,F,X}(U)
\times \Gamma(Y,\Gm)^J\to
P_{k_0,F,Y}(Y^*)$$
$(\MM,u)\mapsto(((\EE,u),\nabla),(V,u),(\rho,u))$
is well-defined. Since these pairings coincide on
$\Gamma(X,\Gm)^J$,
they induce a pairing
$$(\ ,\ )_Y:
P_{k_0,F,X}(U)
\times A_{Y,J}\to
P_{k_0,F,Y}(Y^*)$$
by Corollary of Lemma 3.1.

We define the norm.
Let $X$ and $Y$
be smooth schemes
over $k$ and $D$ and $E$ be divisors with
simple normal crossings
of $X$ and $Y$ respectively.
Let $\alpha:X\to Y$
be a finite flat morphism such that
the reduced inverse image
$(\alpha^*E)_{red}
=D$ and $\alpha:U\to W$ is finite etale on
$U=X-D$ and $W=Y-E$.
Then we define the norm functor
$$N_{X/Y}:
P_{k_0,F,X}(U)
\to
P_{k_0,F,Y}(W)$$
of the commutative Picard categories.
Let $\MM=
((\EE,\nabla),V,\rho)$
be an object of
$P_{k_0,F,X}(U)$.
The invertible $\Cal O _W$-module
$N(\EE)$ is given in [D7] 7.1.
For a basis $e$ of $\EE$ on the inverse
image $\alpha^{-1}(S)$ of an open subset $S\subset U$,
the invertible $\Cal O_Y $-module $N(\EE)$
is generated by the symbol
$N(e)$ on $S$. For another
local basis $e'=ge$,
we require
$N(e')=N_{X/Y}(g)N(e)$.
We define a connection
$\nabla$ on $N(\EE)$
by
$\nabla\log Ne=
\Tr_{X/Y}\nabla\log e$.
It is clear that
$(N(\EE),\nabla)$ is regular along $E$.
The local system $N(V)$ is defined similarly.
For a basis $v$ of $V$ on the inverse
image $\alpha^{-1}(S)$ of an open subset $S\subset U\an$,
the local system $N(V)$
is locally generated by the symbol
$N(v)$ on $S$. For another
local basis $v'=gv$,
we impose
$N(v')=N_{X/Y}(g)N(v)$.
The comparison map $N(\rho)$ is defined by
$N(\rho)N(v)=N_{X/Y}(\varphi)N(e)$
for an analytic function $\varphi$ such that
$\rho(v)=\varphi e$.
It is easily checked that
the norm is an additive functor.

We study a compatibility of the tame symbol and the norm.
Let $\alpha:X\to Y$
be a flat morphism of
smooth $k$-schemes
and $D=\cup_{i\in I}D_i$
and $E=\cup_{i\in I'}E_i$
be divisors with simple normal crossings
of $X$ and $Y$ respectively.
Assume that
$D=(\alpha^*E)_{red}$ and
$(\alpha^*E_i)_{red}$
is smooth for $i\in I'$
and $\alpha|_U$
is smooth.
Let $\alpha_*:I\to I'$
be the map defined by the condition
$\alpha(D_i)\subset E_{\alpha_*(i)}$.
For $J\subset I$
such that
$\alpha_*|_J$
is injective,
we put $J'=\alpha(J)$,
$X_J^*=
X-\cup_{i\in I-J,\alpha(i)\in J'}D_i$
and define a map
$\alpha^*:
A_{Y,J'}\to
A_{X_J^*,J}$
as follows.
For ${i\in J}$,
let $m_i$ be the coefficient of $D_i$ in
$\alpha^*E_{\alpha_*(i)}$,
$m=\prod_{i\in J}m_i$
and $\check m_i=m\cdot m_i^{-1}$.
We define a map
$\alpha^*:
A_{Y,J'}\to
A_{X_J^*,J}$
by
$(f_{i'})_{i'\in J'}
\mapsto
(\alpha^*f_{\alpha_*(i)}^{\check m_i})_{i\in J}$.
It is well-defined and induces the
multiplication by $m$ on the quotient $\ZZ$.

\proclaim{Lemma 3.2}
Assume
$\alpha:X\to Y$
is finite flat.
Let $J'\subset I'$
and $\alpha^*J'=
\{J\in I; \alpha|_J \text{ maps $J$ bijectively to $J'$}\}$.
Then there is an isomorphism between the following two pairings
$P_{k_0,F,X}(U)
\times
A_{Y,J'}\to
P_{k_0,F,E_{J'}}(E_{J'}^*):$
$$
\align
P_{k_0,F,X}(U)
\times
A_{Y,J'}
@>{id\times \alpha^*}>>&
P_{k_0,F,X}(U)
\times
\prod_{J\in \alpha^*J'}A_{X_J^*,J}\\
@>{(\ ,\ )_J}>>&
\prod_{J\in \alpha^*J'}
P_{k_0,F,D_J}(D_J^*)
@>>{\bigotimes_JN_{D_J/E_{J'}}}>
P_{k_0,F,E_{J'}}(E_{J'}^*),
\endalign
$$
$$P_{k_0,F,X}(U)
\times
A_{Y,J'}
@>{N_{X/Y}\times id}>>
P_{k_0,F,Y}(W)
\times
A_{Y,J'}
@>>>
P_{k_0,F,E_{J'}}(E_{J'}^*).
$$
\endproclaim

\demo{Proof}
For
$f\in A_{Y,J'}$
and an object $\MM=((\EE,\nabla),V,\rho)$ of
$P_{k_0,F,X}(U)$,
we define an isomorphism
$$(N_{X/Y}\MM,f)_Y
\simeq
\bigotimes_{J\in\alpha^*J'}
N_{D_J/E_{J'}}(\MM,\alpha^*f).$$
We define an isomorphism of invertible
$\Cal O _{E_{J'}}$-modules
$$\beta:(N_{X/Y}\EE,f)_Y
\simeq
\bigotimes_{J\in\alpha^*J'}
N_{D_J/E_{J'}}(\EE,\alpha^*f)$$
by $(Ne,f)\mapsto
\otimes_JN(e,\alpha^*f)$
on $S\cap E_{J'}$
for a basis $e$ of
$\EE$ on
$\alpha^{-1}(S)\cap U$
for an open subset $S\subset Y$.
We show that it is well-defined,
that it preserves the integrable connections
and that the induced map
$\beta^{an}$
maps the subsheaf
$(N_{X/Y}V,f)$
onto
$\otimes_{J\in \alpha^*J'}
N_{D_J/E_{J'}}(V,\alpha^*f)$.
This will prove that $\beta$
induces an isomorphism
$(N_{X/Y}\MM,f)_Y
\simeq
\bigotimes_{J\in\alpha^*J'}
N_{D_J/E_{J'}}(\MM,\alpha^*f).$
It is sufficient to check projection formulae
$$\gather
(N_{X/Y}g,f)
=
\prod_{J\in\alpha^*J'}
N_{D_J/E_{J'}}(g,\alpha^*f)\\
\nabla\log(N_{X/Y}e,f)=
\sum_{J\in\alpha^*J'}
\nabla\log N_{D_J/E_{J'}}(e,\alpha^*f)\\
(N_{X^{an}/Y^{an}}\varphi,f)=
\prod_{J\in\alpha^*J'}
N_{D_J^{an}/E_{J'}^{an}}(\varphi,\alpha^*f)
\endgather$$
for $g\in \Gamma(U\cap S,\Gm)$,
for a basis $e$ of $\EE$
on $U\cap\alpha^{-1}(S)$ as above and
for a basis $v$ of $V$ and the analytic function
$\varphi=\rho(v)/e$
on $(U\cap\alpha^{-1}(S))\an$ respectively.
Since it is checked etale locally,
we may assume that
$X=Y[\pi'_j,_{j\in J'}]/(\pi_i^{\prime m_i}-\pi_i)$
for some $m_i$, that
$I=J=I'=J'$
and that $D_I\simeq E_I$.
We put $u_i=f_i\pi_i^{-\ord f}\in\Gamma(X,\Gm)$,
$v=g\prod_i\pi_i^{\prime -\ord_i g}\in \Gamma(Y,\Gm)$,
an invertible holomorphic function
$\psi=\varphi\prod_i\pi_i^{\prime \res_i \nabla\log e}$,
and the regular 1-form $\omega=\nabla\log e -
\sum_i\res_i \nabla\log e\cdot d\log\pi'_i$.
Using
$$Ng=
(-1)^{m\sum_i\ord_ig}
Nv\prod_i(-\pi_i)^{\check m_i\ord_i g}$$
with $N(v)|_{E_I}=(v|_{D_I})^m$,
$$\nabla\log Ne=\text{Tr }\nabla\log e=
\text{Tr }\omega+
\sum_i\check m_i\ \res_i \nabla\log e\cdot d\log\pi_i$$
with
$\text{Tr }\omega|_{E_I}=
m\cdot \omega|_{D_I}$,
$$N\varphi=
(-1)^{-m\sum_i\res_i\nabla\log e}
N\psi\prod_i(-\pi_i)^{-\check m_i\res_i\nabla\log e}$$
with $N\psi|_{E_I}=(\psi|_{D_I})^m$
and the definition
$\alpha^*f=
(\alpha^*f_i^{\check m_i})_{i\in J}$,
we obtain
$$
\align
&(Ng,f)
=(g,\alpha^*f)\\
=&
(-1)^{m\cdot\ord f\sum_i\ord_i g}
(v|_{D_I})^{m\cdot\ord f}
\prod_iu_i|_{E_I}^{-\check m_i\cdot\ord_ig},
\endalign
$$$$\align
&\nabla\log(Ne,f)
=\nabla\log(e,\alpha^*f)\\
=&m\cdot\ord f\cdot\omega|_{D_I}-
\sum_i\res_i(\nabla\log e)
\cdot\check m_i\frac{du_i}{u_i}|_{E_I},
\endalign
$$$$\align
&(N\varphi,f)
=(\varphi,\alpha^*f)\\
=&(-1)^{m\cdot\ord f\sum_i\res_i(\nabla\log e)}
(\psi|_{D_I})^{m\cdot\ord f}
\prod_i(u_i|_{E_I})^{-\check m_i\cdot\res_i(\nabla\log e)}.
\endalign
$$
Thus an isomorphism
$$(N_{X/Y}\MM,f)_Y
\simeq
\bigotimes_{J\in\alpha^*J'}
N_{D_J/E_{J'}}(\MM,\alpha^*f)$$
is defined.
It is easy to check that
it defines an isomorphism of pairings.
\enddemo

\proclaim{Corollary}
Assume further that $\dim X=1$ and
$E=\{y\}$.
Then there is an isomorphism between the following two pairings
$P_{k_0,F,X}(U)
\times
k(Y)^\times\to
P_{k_0,F}(y)$
$$
\align
P_{k_0,F,X}(U)
\times
k(Y)^\times
@>{id\times \partial_x\circ\alpha^*}>>
&P_{k_0,F,X}(U)
\times
\prod_{x \in D}A_x\\
@>{((\ ,\ )_x)_x}>>&
\prod_{x\in D}
P_{k_0,F}(x)
@>{\bigotimes_xN_{x/y}}>>
P_{k_0,F}(y)
\endalign
$$$$P_{k_0,F,X}(U)
\times
k(Y)^\times
@>{N_{X/Y}\times \partial_y}>>
P_{k_0,F,Y}(W)
\times
A_y@>{(\ ,\ )_y}>>
P_{k_0,F}(y).
$$
\endproclaim


Now we define the pairing.
For a closed point $x\in X_0$,
We put
$I_x=\{i\in I; x\in D_i\}$.
Since the abelian group
$A_x$ in the last section is the same as
$A_{x,I_x}$
defined above for the inclusion $\{ x \} \to D_{I_x}$,
we have defined a pairing
$$(\ ,\ )_x:
P_{k_0,F,X}(U)
\times
A_x\to
P_{k_0,F}(x)$$
By compositing with the norm
$N_{x/k}:P_{k_0,F}(x)\to
P_{k_0,F}(k),$
we obtain a pairing
$$(\ ,\ )_X=\bigotimes_{x\in X_0}N_{x/k}\circ(\ ,\ )_x:
P_{k_0,F,X}(U)
\times
\bigoplus_{x\in X_0}A_x\to
P_{k_0,F}(k).$$

\proclaim{Proposition 2}
Assume further that $X$ is projective and $n=\dim X$.
Then there exists a pairing
$$(\ ,\ )_X:
P_{k_0,F,X}(U)
\times
\bold{CH}^n(X\bmod D)\to
P_{k_0,F}(k)$$
inducing the pairing
$(\ ,\ )_X$ above
by $\bigoplus_{x\in X_0}A_x\to
\bold{CH}^n(X\bmod D)$.
It induces a pairing
$(\ ,\ )_X:
MPic_{k_0,F,X}(U)
\times
{CH}^n(X\bmod D)\to
MPic_{k_0,F}(k)$
on the class groups.
\endproclaim

\demo{Proof}
The second assertion follows from
the adelic presentation, Proposition 1
$$CH^n(X\bmod D)=
\Coker(\bigoplus_{y\in X_1}B_y\to
\bigoplus_{x\in X_0}A_x).$$
By Lemma 3.1,
it is enough to define a trivialization
of the pairing
$(\ ,\partial \ )_X:
P_{k_0,F,X}(U)
\times B_y\to
P_{k_0,F}(k):
(\MM,f)\mapsto\bigotimes_{x\in \overline{\{y\}}}
N_{x/k}(\MM,\partial_{x,y}f)_x$
for $y\in X_1$.

First we consider the de Rham component.
Let $P(U)$ be the commutative Picard category
of invertible $\Cal O_U $-modules
and $P(k)$ be that of invertible $k$-vector spaces.
We define a trivialization of the pairing
$(\ , \partial \ )_X:P(U)\times B_y\to P(k)$.
We recall that for an invertible
$\Cal O_U $-module $\EE$ and $f\in A_x$
for $x\in X_0$,
the invertible
$\kx$-vector space
$(\EE,f)_x$
is generated by the symbol
$(e,f)$ for a local basis
$e$ of $\EE$ on
$X_x\cap U$.
For another basis
$e'=ge$ and
$g\in \Gamma(X_x\cap U,\Gm)$,
we have
$(e',f)=(g,f)(e,f)$
where
$$(g,f)=(-1)^{\ord_xf\cdot\sum_{i\in I_x}\ord_ig}
\bar g^{\otimes \ord_xf}\otimes
\bigotimes_{i\in I_x}f_i^{-\ord_i g}$$
and $\bar g$ is the basis of
the invertible sheaf $\Cal O (-\text{div }g)$
defined by $g$ at $x$.

For an invertible $\Cal O_U $-module
$\EE$ and for $f\in B_y$, $y\in X_1$,
we define a trivialization
$$k\to (\EE,\partial f)=
\bigotimes_{x\in X_0}
N_{x/k}(\EE,\partial_{x,y}f)_x$$
as follows.
Let $Y$ be the closure of $\{y\}$
and $\Sigma
=\{x\in Y;\partial_{x,y}f\ne0,
x\in D_i\text{ for some }\allowmathbreak i\notin I_y,
\text{ or $\bar f$ is not invertible at $x$}\}$
be a finite closed subset of $Y$,
where $\bar f$ is the image of $f$ by $B_y\to \kyx$.
Let $X_\Sigma$
be the semi-localization of $X$ at $\Sigma$.
It is the spectrum of a semi-local ring
since the finite set
$\Sigma$ is contained in an affine open subset
of $X$. The restriction of $\EE$ on $X_\Sigma\cap U
$ is free of rank 1.
We define the trivialization by
$$1\mapsto (e,\partial f)_X:=
\prod_{x\in Y-\Sigma}N_{x/k}\bar f(x)^{-\ord_x(e)}\times
\bigotimes_{x\in \Sigma}
N_{x/k}(e,\partial_{x,y}f)_x$$
for a basis $e$ of $\EE$ on $X_\Sigma\cap U$
if $\Sigma\neq\emptyset$ and by the identity if $\Sigma=\emptyset$.
Here for a basis $e$ of $X_\Sigma\cap U$,
let an invertible $\Cal O_Y $-module $\EE_Y(e)$
be the restriction of an invertible
$\Cal O_X $-module extending $\EE$ on $U$ and
$\Cal O \cdot e$ on $X_\Sigma$.
It is uniquely determined if $\Sigma\neq\emptyset$
and $e$ defines a non-zero rational section of
$\EE_Y(e)$. We write $\ord_x(e)$ its order at a closed point $x\in Y$.

We show $(e,\partial
f)_X$ is independent of the choice of a basis $e$
and the trivialization is well-defined.
For $g\in \Gamma(X_\Sigma\cap U,\Gm)$,
let $\bar g$
be the section of
$\Cal O_Y (-\sum_{i\in I_y}\ord_ig\cdot D_i)$,
$\Sigma_g=
\{x\in Y;\text{$\bar g$ is not a basis}\}$
and $\ord_x(\bar g)$ be its order at $x\in Y$.
We put
$$(g,\partial f)_x=\cases
(g,\partial_{x,y}f)_x &\quad x\notin\Sigma_g\\
\bar f(x)^{-\ord_x\bar g} &\quad x\in \Sigma_g
\endcases$$
$\in \kxx$ for $x\in Y$
and
$(g,\partial f)_X=\prod_{x\in Y}
N_{x/k}(g,\partial f)_x\in k^\times$.
Since
$(ge,\partial f)=
(g,\partial f)
(e,\partial f)$,
it is enough to prove
$(g,\partial f)_X=1$.

Let $\pi=(\pi_i)_{i\in I_y}$
be a family of bases of
$\Cal O_Y (-D_i)$ on a neighbourhood of $\Sigma\amalg\Sigma_g$.
We take a lift of $\pi_i$ to a basis of
$\Cal O _{X_\Sigma}(-D_i)$ and write it also $\pi_i$.
Let $\{\pi,\bar f\}\in B_y$
be the image of $\pi\otimes \bar f
\in A_y\otimes \kyx$.
Since $\Sigma$
for
$\{\pi,\bar f\}$ and for
$f -\{\pi,\bar f\}$ are disjoint with
$\Sigma_g$,
it is sufficient to show
$(g,\partial\{\pi,\bar f\})_X
=(g,\partial(f-\{\pi,\bar f\}))_X=1$.
In other word, in the proof of
$(g,\partial f)_X=1$, we may assume one of the conditions
\roster
\item
$\bar f=1$.
\item
$f=\{\pi,\bar f\}$ where $\pi=(\pi_i)_{i\in I_y}$ as above.
\endroster
In case (1),
we have
$f=(f_i)_i\in \bigoplus_{i\in I_y}K_2(y)\subset B_y$,
$(\partial_{x,y}f)_x=
(\partial_xf_i)_i\in \bigoplus_{i\in I_y}\kxx\subset A_x$
and
$(g,\partial_{x,y}f)_x=
\prod_{i\in I_y}(\partial_xf_i)^{-\ord_i g}$.
Hence
$$\prod_{x \in Y}N_{x/k} (g,\partial f)_x=
\prod_{i\in I_y}(\prod_{x\in Y}\partial_xf_i)^{-\ord_i g}=1$$
as required by the reciprocity law for $K_2(y)$.

We consider the case (2).
Let $$(g,\pi)=
(-1)^{\sum_{i\in I_y}\ord_ig}
\bar g\otimes
\bigotimes_i\pi_i^{\otimes -\ord_i g}
\in \kyx $$
be a rational function on $Y$.
We show that
$(g,\partial\{\pi,\bar f\})_x=
((g,\pi),\bar f)_x.$
It implies
$(g,\partial\{\pi,\bar f\})_X=1$
by the reciprocity law
$\prod_{x\in Y}N_{x/k}((g,\pi),\bar f)_x=1.$
For a closed point $\tilde x$
of the normalization
$\tilde Y$ of $Y$,
we put
$$(g,\partial f)_{\tilde x}=\cases
(g,\partial_{\tilde x,y}f)_{\tilde x} &
\quad \tilde x \mapsto x\notin\Sigma_g\\
\bar f(x)^{-\ord_{\tilde x}\bar g} &\quad x\in \Sigma_g
\endcases$$
$\in \kxt^\times$ similarly as before.
Since
$(g,\partial f)_x=\prod_{\tilde x\mapsto x}
N_{\tilde x/x}(g,\partial f)_{\tilde x}$,
it is enough to show
$(g,\partial\{\pi,\bar f\})_{\tilde x}=
((g,\pi),\bar f)_{\tilde x}$
for $\tilde x\mapsto x\in Y$.
We extend $\pi=(\pi_i)_{i\in I_y}$
to $\tilde \pi=(\pi_i)_{i\in I_x}$
where $\pi_i$ is a basis of $\Cal O_X (-D_i)$ at $x$.
Then by the explicit computation of $\partial_{\tilde x,y}$
given in section 2,
we have
$$\partial_{\tilde x,y}\{\pi,\bar f\}=
\tilde\pi^{\ord_{\tilde x}f}\times
((-1)^{\ord_{\tilde x}f}
(\pi_i,f)^{-1}_{\tilde x}
)_{i\in I_x-I_y}
$$
where $\tilde\pi$ also denotes the element of $A_{\tilde x}$
defined by $\tilde\pi$.
By putting
$u=$ \linebreak $(-1)^{\sum_{i\in I_y}\ord_ig}g\cdot
\prod_{i\in I_x}
\pi_i^{-\ord_ig}\in\Cal O _{X,x}^{\times}$,
the left hand side
$(g,\partial\{\pi,\bar f\})_{\tilde x}$
 is
$$\gather
\left((-1)^{\sum_{i\in I_x}\ord_ig}g\cdot
\prod_{i\in I_x}
\pi_i^{-\ord_ig}\right)^{\ord_{\tilde x}f}
\times
\prod_{i\in I_x-I_y}
((-1)^{\ord_{\tilde x}f}
(\pi_i,f)_{\tilde x}^{-1})^{-\ord_ig}\\
=u^{\ord_{\tilde x}f}\cdot
\prod_{i\in I_x-I_y}
(\pi_i,f)_{\tilde x}^{\ord_ig}.
\endgather$$
By $(g,\pi)=
(-1)^{\sum_{i\in I_y}\ord_ig}g
\prod_{i\in I_y}\pi_i^{-\ord_ig}=
u\cdot
\prod_{i\in I_x-I_y}\pi_i^{\ord_ig}$,
the right hand side is the same and the equality
$(g,\partial\{\pi,\bar f\})_{\tilde x}=
((g,\pi),\bar f)_{\tilde x}$
is proved.
Thus the equality
$(g,\partial f)_X=1$ is proved in general
and hence the isomorphism
$k\to(\EE,\partial f)_X$ is well-defined.
We call it the canonical isomorphism.
It is easy to check that the canonical isomorphisms
define an isomorphism of pairings
$1\to (\ ,\partial \ )_X:
P(U)\times B_y\to P(k)$.

We consider $P_{k_0,F,X}(U)$.
We define for
$y\in X_1$ and $f\in B_y$,
an isomorphism of the additive functors
$$1\to (\ ,\partial f )_X:
P_{k_0,F,X}(U)\to P_{k_0,F}(k)$$
whose de Rham component is the canonical isomorphism
defined above.
We show that they define an isomorphism of the pairings
$$1\to (\ ,\partial \ )_X:
P_{k_0,F,X}(U)\times B_y\to P_{k_0,F}(k).$$
The isomorphisms of functors
$(\ ,\partial f)\otimes(\ ,\partial f')
\to(\partial ff')$
is defined as the composite
$(\ ,\partial f)\otimes(\ ,\partial f')
\gets1\otimes 1
\to1
\to(\partial ff')$.
The compatibilities in Lemma 3.1.1
follow from those of the de Rham component
since the forgetful functor
$P_{k_0,F}(k)\to P(k)$
is faithful.
We define the isomorphism
$1\to (\ ,\partial f)_X$
first for curves and then reduce the general case to it.

Assume $X$ is a curve.
We may assume $X$ is integral and $y$ is the
generic point of $X$.
The group $B_y$
is then the function field $\kyx$.
First we consider the case where a function $f\in \kyx$
is constant.
By considering the norm,
we may assume that the constant field
of $X$ is $k$.
We take a branch of the logarithm
$\log f\in \CC (1)^{Hom_{k_0}(k,\CC )}$
of $f\in k^\times$.
Let $\MM=((\EE,\nabla),V,\rho)$
be an object of $P_{k_0,F}(U)$.
Let $e$ be a non-zero rational section of $\EE$.
For $x\in D$,
let $\nabla_x=\res_x\nabla\log e\in \kx$
and $[f^{\nabla_x}]$ be the object
$(\kx,F,\text{ multiplication by }
f^{\nabla_x})$
of $P_{k_0,F}(x)$ where
$f^{\nabla_x}=
\exp(\nabla_x\cdot\log f)
\in (\kx\otimes_{k_0}\CC)^\times)$.
Note that $\partial_{x,y}f\neq0$
if and only if $x\in D$
and then $\partial_{x,y}f=f\in \kxx\subset A_x$.
For $x\in D$,
the isomorphism
$\kx\to
(\EE,\partial_{x,y}f)
:1\mapsto(e,\partial_{x,y}f)$
induces an isomorphism
$[f^{\nabla_x}]
\to(\MM,\partial_{x,y}f)$
in $P_{k_0,F}(x)$.
By taking the tensor product of the norms,
we obtain the isomorphism
$[f^{\sum_{x\in D}\Tr_{x/k}\nabla_x}]
\to(\MM,\partial f)$
whose de Rham component is
$1\to\bigotimes_{x\in D}N_{x/k}(e,\partial_{x,y}f)$.
We have
$\sum_{x\in D}\Tr_{x/k}\nabla_x=
-\sum_{x\notin D}[\kx:k]\ord_xe$
by the residue formula
$\sum_{x\in X}\Tr_{x/k}\res_x\nabla\log e=0$
and by $\ord_xe=
\res_x\nabla\log e$
for $x\notin D$.
Therefore the canonical isomorphism
$k\to(\EE,\partial f)$
defined previously induces an isomorphism
$1\mapsto (\MM,\partial f)$
in $P_{k_0,F}(k)$
also called canonical.
It is easy to see that the canonical isomorphisms
define an isomorphism
$1\to (\ ,\partial f )_X:
P_{k_0,F,X}(U)\to P_{k_0,F}(k)$
of the additive functors.

We assume that a function
$f\in \kyx$ is not constant. It
defines a finite flat morphism
$\alpha:X\to Y=\Bbb P^1_k$.
By shrinking $U$ if necessary,
we may assume that
$U=\alpha^*\alpha(U)$
and $\alpha|_U$
is etale.
We put
$E=\Bbb P^1_k-\alpha(U)$
and $t$ be the coordinate of $\Bbb P^1_k$.
Then by Lemma 3.2,
there is a canonical isomorphism
$$(\ ,\partial f)_X\simeq
(\ ,\partial t)_Y\circ N_{X/Y}$$
of additive functors
$P_{k_0,F,X}(U)\to
P_{k_0,F}(k)$.
Since the canonical isomorphism
$(\EE,f)_X\to(N(\EE),t)_Y$ sends
$(e,f)$ to $(N(e),t)$,
we are reduced to the case
$X=\Bbb P^1_k$
and $f=t$ is the coordinate.

Now we assume
$X=\Bbb P^1_k$
and $f=t$ is the coordinate
and define an isomorphism of functors
$$(\ ,\partial t)_X\simeq1 : P_{k_0,F,X}(U)\to
P_{k_0,F}(k)$$
of the commutative Picard categories.
Let $\MM=
((\EE,\nabla),V,\rho)$
be an object of $P_{k_0,F,X}(U)$.
By shrinking $U$ if necessary,
we may assume $0,\infty\in D$
and we may take a basis $e$ of $\EE$ on $U$.
We identify
$k=(\EE,\partial t)_X$
by the canonical map
$1\mapsto\bigotimes_{x\in D}
N_{x/k}(e,\partial_xt)$.
It is enough to show that
the local system
$(V,\partial t)_X
\subset
(\EE,\partial t)_X\an=\CC$
of $F$-vector spaces
coincides with
$F\subset\CC$
on $k\an$.
Since it suffices to show it for each embedding
$k\to \CC$ over $k_0$,
we may assume $k=\CC$.

We take a simply connected path $\gamma$
connecting $0$ to $\infty$
which does not pass other points in $D$.
Let $|\gamma|$
denote the underlying set and
put $|\gamma|^*=|\gamma|-\{0,\infty\}$.
We take a small tubular neighborhood
$T$ of $|\gamma|^*$.
There are 2 connected component
$T^+$ and $T^-$ of $T-|\gamma|^*$.
We assume that a small circle
around 0 with positive direction
cuts $\gamma$ from $T^-$ to $T^+$.
We take a branch $\log t$ of the logarithm
of the coordinate $t$ on the simply connected region
$\PP_{\CC}-|\gamma|$.
Let $\log^+t$ and $\log^-t=\log^+t+2\pi\sqrt{-1}$
be the branch of the logarithm
on $|\gamma|^*$
continuous on
$|\gamma|^*\cup T^+$ and
$|\gamma|^*\cup T^-$
respectively.
Take a basis $v$ of
the $F$-vector space
$V(|\gamma|^*)$
of dimension 1.
We define a basis
$(v,\partial_xt)$
of
$(V,\partial_xt)$
for $x\in D$.
For $x\ne 0,\infty$,
let $(v,t)=
t^{\nabla_x}$
be the section
$\exp(\nabla_x\cdot\log t)$
of $(V,t)^\sim
=\Cal O_X \an$
and
$(v,\partial_xt)$
be its fiber at $x$.
For $x=0$,
let $(v,t)=
(-1)^{-\nabla_0}t^{\nabla_0}v$
be the section
$\exp(\nabla_0(-\pi\sqrt{-1}+\log^- t))\cdot v$
of $(V,t)^\sim|_{|\gamma|^*}$
and
$(v,\partial_0t)$
be its fiber at $0$.
Similarly for $x=\infty$,
let $(v,t)=
(-1)^{\nabla_\infty}t^{\nabla_\infty}v^{\otimes-1}$
be the section
$\exp(\nabla_\infty(\pi\sqrt{-1}+\log^+t))\cdot v^{\otimes-1}$
of $(V,t)^\sim|_{|\gamma|^*}$
and
$(v,\partial_\infty t)$
be its fiber at $\infty$.

It is enough to show that
$$\prod_{x\in D}
\frac{(\rho_,\partial_xt)
(v,\partial_xt)}
{(e,\partial_xt)}
=1.$$
We compute each term for $x\in D$.
For $x\ne0,\infty$,
it is
$t(x)^{\nabla_x}=
\exp(\nabla_x\log t(x))$.
Here $\log$ is the branch chosen above
and
$\nabla_x=\res_x(\nabla\log e)$.
We take a base point $b$ on $|\gamma|^*$.
Then for $x=0$ and $\infty$
they are
$$\gather
\exp\left(\nabla_0(-\pi\sqrt{-1}+\log^- t(b))
+\int_b^0(\nabla_0d\log t-
\nabla\log e)\right)
\times
\frac{\rho(v)(b)}{e(b)}\\
\exp\left(\nabla_\infty(\pi\sqrt{-1}+\log^+ t(b))
+\int_b^\infty(\nabla_\infty d\log t+
\nabla\log e)\right)
\times
(\frac{\rho(v)(b)}{e(b)})^{-1}.
\endgather$$
Therefore it is sufficient to show
$$\gather
\nabla_0(\pi\sqrt{-1}+\log^+ t(b))
+\int_b^0(\nabla_0d\log t-
\nabla\log e)\\
+\nabla_\infty(\pi\sqrt{-1}+\log^+ t(b))
+\int_b^\infty(\nabla_\infty d\log t+
\nabla\log e)\\
+\sum_{x\in D,\ne 0,\infty}
\nabla_x\cdot\log t(x)=0.
\endgather$$

For $x\in D,\ne0,\infty$,
let $\gamma_x$ be a small circle with positive direction
around $x$.
For $x=0$,
let $\gamma_0$
be a loop starting at $b$,
going along $\gamma$ in $T^+$,
turning around 0 in the positive direction
and going back to $b$ along $\gamma$ in $T^-$.
Similarly for $x=\infty$,
let $\gamma_\infty$
be a loop starting at $b$,
going along $\gamma$ in $T^-$,
turning around $\infty$ in the positive direction
and going back to $b$ along $\gamma$ in $T^+$.
Then by Cauchy's integral theorem,
the sum of the integral of the
holomorphic 1-form
$$\sum_{x\in D}\left(\frac1{2\pi\sqrt{-1}}
\int_{\gamma_x}
\log t\nabla\log e\right)=0.$$
Hence it is enough to show that the integrals
are equal to the corresponding terms in the equality above.
For $x\in D,\ne0,\infty$,
it is the theorem of residue.
For 0 and $\infty$,
we have
$$
\frac1{2\pi\sqrt{-1}}
\int_{\gamma_0}
\log td\log t=
-\frac1{2\pi\sqrt{-1}}
\int_{\gamma_\infty}
\log td\log t=
\log^+ t(b)+
\pi\sqrt{-1}$$
Since
$\nabla\log e-\nabla_0d \log t$ is holomorphic
at 0, we have
$$
\frac1{2\pi\sqrt{-1}}
\int_{\gamma_0}
\log t
(\nabla\log e
-\nabla_0d \log t)=
\int_b^0
\nabla_0d \log t
-\nabla\log e.$$
Similarly for $\infty$,
we have
$$
\frac1{2\pi\sqrt{-1}}
\int_{\gamma_\infty}
\log t
(\nabla\log e
+\nabla_\infty d \log t)=
\int_b^\infty
\nabla\log e
+\nabla_0d \log t.$$
Thus Proposition 2 is proved for $X=\PP_k$
and hence for curves.

We prove Proposition 2 for higher dimension.
We define an isomorphism
$1\to(\ ,\partial f)_X:
P_{k_0,F,X}(U)\to
P_{k_0,F}(k)$
for $f\in B_y$ and $y\in X_1$.
First we consider the case $f\in
\bigoplus_{i\in I_y}K_2(y)=
\text{ Ker }(B_y\to\kyx)$.
Let $D_i$ be an irreducible component of $D$
and $k_i$ be the constant field of $D_i$.
We define a pairing
$$(\ ,\ )_i:P_{k_0,F,X}(U)\times k_i^\times\to
P_{k_0,F}(k_i)$$
as follows.
For an invertible
$\Cal O_U $-module $\EE$ and $u\in k_i^\times$,
let $(\EE,u)$
be the invertible $k_i$-vector space generated by the
symbol $(e,u)$ for a non-zero rational section $e$ of $\EE$.
For another section $e'=ge$,
we impose $(e',u)=(g,u)(e,u)$
where $(g,u)=u^{-\ord_ig}\in k_i^\times$.
For an object
$\MM=((\EE,\nabla),V,\rho)$ of
$P_{k_0,F,X}(U)$,
we define the local system of $F$-vector space
$(V,u)$ on $\Spec k_i\an$ to be that generated by
$(v,u)=u^{\res_i\nabla\log e}(e,u)$
and the comparison $(\rho,u)$
to be the natural inclusion to $(\EE,u)\an$.
It is easy to see that they define a pairing as above.

For $x\in X_0$,
the restriction of the pairing
$N_{x/k}\circ(\ ,\ )_x:
P_{k_0,F,X}(U)
\times A_x\to
P_{k_0,F}(k)$
to the
$i$-component
$\kxx\subset A_x$
is identified with the pairing
$N_{k_i/k}\circ(\ ,N_{x/k_i}(\ ))_i$
for $i\in I_x$.
On $\bigoplus_{i\in I_y}K_2(y)=
\text{ Ker }(B_y\to\kyx)$,
the map $\partial_{x,y}:B_y\to A_x$
is the direct sum of the tame symbol
$\partial_x:K_2(y)\to \kxx$.
Hence the restriction of the pairing
$(\ ,\partial f)_X:
P_{k_0,F,X}(U)\times B_y\to
P_{k_0,F}(k)$
to the $i$-component $K_2(y)\subset B_y$
for $i\in I_y$ is
identified with the pairing
$N_{k_i/k}\circ(\ ,\prod_{x\in Y}N_{x/k_i}(\partial_x\ ))_i$
where $Y$ is the closure of $\{y\}$.
Hence by the reciprocity law
$1=\prod_{x\in Y}N_{x/k_i}(\partial_x \ ):
K_2(y)\to k_i^\times$,
a trivialization of
the restriction of the pairing
$(\ ,\partial \ )_X$
to $\text{Ker }(B_y\to \kyx)$
is defined.
It is easy to see that the de Rham component is the canonical one.

We prove the general case by reducing to a curve.
The assertion has been proved for an element in the kernel of
$B_y\to \kyx$.
Hence it is sufficient to show that, for
$f\in\kyx$,
there exist an element $\tilde f\in B_y$
lifting $f$ and a trivialization of the functor
$1\to(\ ,\partial \tilde f)_X$
whose de Rham component is the canonical one.
Let $Y$ be the closure of $\{y\}$ for a point $y\in X_1$.
Let $\pi=(\pi_i)_{i\in I_y}$
be a family of bases
$\pi_i$ of
$\Cal O _Y(-D_i)|_W$
on a dense open subset $W\subset Y$.
The element $\pi\in A_y$
defined by $\pi$ has
$\ord_y\pi=1$
and the image of
$\tilde f=\{\pi,f\}\in B_y$
in $\kyx$ is $f$.
The tame symbol
$(\ ,\pi)$
defines an additive functor
$$
P_{k_0,F,X}(U)
\to
P_{k_0,F,\tilde W}
(\tilde W^*)$$
where $\tilde W$ is the inverse image of
$W$ in the normalization $\tilde Y$ of $Y$
and $\tilde W^*$ is the inverse image of
$W-\bigcup_{i\notin I_y}D_i$.

\proclaim{Claim} 1.
For a closed point $\tilde x \in\tilde W$,
there is an isomorphism of additive functors
$$(\ ,\partial_{\tilde x}\tilde f)_{\tilde x}
\to((\ ,\pi),f)_{\tilde x}:
P_{k_0,F,X}(U)
\to
P_{k_0,F}(\tilde x)$$
whose de Rham component is defined by
$(e,\partial_{\tilde x,y}\{\pi,f\})\mapsto
((e,\pi),f)$ for a basis $e$ of $\EE$
on the intersection of $U$ with a neighborhood of $x$.
The first functor is defined by the pairing
$(\ ,\ )_{\tilde x}:
P_{k_0,F,X}(U)\times A_{\tilde x}
\to
P_{k_0,F}(\tilde x)$ and the second is defined by
$(\ ,\ ):P_{k_0,F,\tilde W}(\tilde W^*)
\times A_{\tilde x,\tilde Y}
\to
P_{k_0,F}(\tilde x).$

2. Let $x\in Y-W$.
Assume that $Y$ is regular at $x$ and
that $x\notin D_i$ for $i\notin I_y$.
Then the image of the functor
$(\ ,\pi)$ is in
$P_{k_0,F,\tilde W\cup\{x\}}
(\tilde W^*)$.
Further assume that $f$ is invertible at $x$.
Then there is an isomorphism of additive functors
$$(\ ,\partial_{\tilde x}\tilde f)_{\tilde x}
\to((\ ,\pi),f)_{\tilde x}:
P_{k_0,F,X}(U)
\to
P_{k_0,F}(\tilde x)$$
whose de Rham component is defined in the same way as in 1.
\endproclaim

We complete the proof of Proposition 2 admitting Claim.
For $f\in\kyx$,
we may take a family
$\pi=(\pi_i)_{i\in I_y}$
of bases
$\pi_i$ of
$\Cal O _Y(-D_i)|_W$
on an open subset $W\subset Y$
containing
the singular points of $Y$,
the intersection $D_i\cap Y$
for $i\notin I_y$
and the zeroes and poles of $f$.
Then by the first assertion of Claim 2,
the tame symbol $(\ ,\pi)$ is a functor
$P_{k_0,F,X}(U)\to
P_{k_0,F,\tilde Y}
(\tilde W^*)$.
Further the isomorphisms in Claim 1 and 2 induce
an isomorphism of functors
$$(\ ,\partial\tilde f)_X
\to((\ ,\pi),f)_{\tilde Y}:
P_{k_0,F,X}(U)
\to
P_{k_0,F}(k).$$
Since we have defined the trivialization of the second functor
for a curve $\tilde Y$,
we obtain a trivialization of the first functor
$(\ ,\partial\tilde f)_X$.
Its de Rham component is the canonical one and
hence the assertion is proved for $\tilde f$.
Thus Claim completes the proof of Proposition 2.
\enddemo
\demo{Proof of Claim} 1.
Let $\MM=
((\EE,\nabla),V,\rho)$
be an object of
$P_{k_0,F,X}(U)$.
It is sufficient to show that the isomorphism
$(\EE,\partial_{\tilde x}\tilde f)_{\tilde x}
\to((\EE ,\pi),f)_{\tilde x}$
of the de Rham component
maps the $F$-vector space
$(V,\partial_{\tilde x}\tilde f)_{\tilde x}
\to((V ,\pi),f)_{\tilde x}$.
Take bases $e$ of $\EE$ and $v$ of $V$
in the intersection of $U$ with a neighborhood of $x$
and let
$\varphi=\rho(v)/e$
be an analytic function.
It is enough to show that
$$(\varphi ,\partial_{\tilde x}\tilde f)_{\tilde x}
=((\varphi ,\pi),f)_{\tilde x}$$
modulo $F^\times$.
Extend the family
$\pi=(\pi_i)_{i\in I_y}$
to the family $\tilde \pi=(\pi_i)_{i\in I_x}$
of bases of $\Cal O _{X_\Sigma}(-D_i)$ at $x$ as before.
By the explicit computation
$$\partial_{\tilde x,y}\{\pi,\bar f\}=
\tilde\pi^{\ord_{\tilde x}f}\times
((-1)^{\ord_{\tilde x}f}
(\pi_i,f)^{-1}_{\tilde x}
)_{i\in I_x-I_y},
$$
the left hand side is
$$
((-1)^{-\sum_{i\in I_x}\res_i\nabla\log e}\varphi\cdot
\prod_{i\in I_x}
\pi_i^{\res_i\nabla\log e})^{\ord_{\tilde x}f}
\cdot
\prod_{i\in I_x-I_y}
((-1)^{\ord_{\tilde x}f}
(\pi_i,f)_{\tilde x}^{-1})^{\res_i\nabla\log e}.$$
By $\res_{\tilde x}\nabla\log(e,\pi)=
\sum_{i\in I_x-I_y}
m_i\cdot\res_i\nabla\log e$
where $m_i=\ord_{\tilde x}\pi_i$
for $i\in I_x-I_y$, the right hand side is
$$\gather
(-1)^{-\sum_{i\in I_x-I_y}
m_i\res_i\nabla\log e\cdot
\ord_{\tilde x}f}
\times
\left((-1)^{-\sum_{i\in I_y}\res_i\nabla\log e}\varphi\cdot
\prod_{i\in I_y}
\pi_i^{\res_i\nabla\log e}\right)^{\ord_{\tilde x}f}\\
\times
f^{\sum_{i\in I_x-I_y}
m_i\res_i\nabla\log e}.
\endgather$$
Cancelling the functors of $\varphi$
and $i\in I_y$,
the equality follows from the definition of the tame symbol
$\pi_i^{\ord_{\tilde x}f}(\pi_i,f)^{-1}_{\tilde x}
=(-1)^{-m_i\ord_{\tilde x}f}f^{m_i}$
for $i\in I_x-I_y$.

2.
Let $\pi'=(\pi'_i)_{i\in I_y}$
be a family of bases of $\Cal O_Y (-D_i)$ at $x$
and put $u_i=\pi_i/\pi'_i\in \kyx$.
By the definition of the tame symbol,
we have
$$\nabla\log(e,\pi)=
\nabla\log(e,\pi')-
\sum_{i\in I_y}\res_i(\nabla\log e)\cdot
d\log u_i$$
for a basis $e$ of $\EE$
on the intersection of $U$ with a neighborhood of $x$.
Since
$\nabla\log(e,\pi')$
is regular at $x$ and
$d\log u_i$
has at most logarithmic pole at $x$,
the first assertion is proved.

To show the second assertion,
by the same argument as in 1,
it is enough to show
$$(\varphi ,\partial_{x}\tilde f)_{x}
=((\varphi ,\pi),f)_{x}$$
for $\varphi=\rho(v)/e$ as above.
By $\partial_{x,y}f=
\pi^{\prime\ord_xf}\times(\{u_i,f\}_x)_i=
(f(x)^{-\ord_x\pi_i})_i\in
A_x$,
the left hand side is
$f(x)^{-\sum_i\ord_x\pi_i\cdot\res_i(\nabla\log e)}
=\exp(\sum_i-\ord_x\pi_i\cdot\res_i(\nabla\log e)\log f(x)).$
By $\res_x\nabla\log (e,\pi)=
-\sum_{i\in I_y}\res_i(\nabla\log e)\cdot
\res_xd\log u_i=
-\sum_{i\in I_y}\res_i(\nabla\log e)\cdot
\ord_x\pi_i$,
the right hand side is the same.
Thus Claim 2 and hence Proposition 2 is proved.
\enddemo

\heading
4. Main results.
\endheading

Let $k_0$ and $F$ be subfields of $\CC$
and $U$ be a smooth separated scheme
over a finite extension $k$ of $k_0$.
Let $X$ be a proper smooth scheme over
$k$ including $U$ as the complement
of a divisor $D$ with simple normal crossings.
For an object $\MM$
of $M_{k_0,F}(U)$,
the determinant of the cohomology
$\det_\Gamma R\Gamma_c(U,\MM)$
is defined in section 1.
Let $(\det\MM,c_{X\bmod D})$
denote an object of $P_{k_0,F}(k)$
whose class is the pairing
$([\det\MM],c_{X\bmod D})$
defined in section 3 where
$c_{X\bmod D}$ is
the relative canonical class
defined in previous section 2.

\proclaim{Theorem 1}
Assume $X$ is projective.
Then for an object $\MM$ of
$M_{k_0,F}(U)$,
there exists an isomorphism
$$
\det{}^\Gamma R\Gamma_c(U,\MM)
\otimes
\det{}^\Gamma R\Gamma_c(U,1)^{\otimes-\rank\MM}
\simeq
(\det\MM,c_{X\bmod D})
$$
in
$P_{k_0,F}(k)$.
In other word, we have an equality
$$per_c(\MM)\cdot
per_c(1)^{-\rank\MM}
=
\Gamma(\nabla:\MM)^{-1}\times
([\det\MM],c_{X\bmod D})$$
in
$MPic_{k_0,F}(k)
=k^\times\backslash
(k\otimes_{k_0}\CC)^\times/
(F^\times)^{Hom_{k_0}(k,\CC)}$.
\endproclaim
\proclaim{Corollary}
Under the same assumption,
we have an equality
$$per(\MM)\cdot
per(1)^{-\rank\MM}
=
\Gamma(-\nabla:\MM)^\times
[(\det\MM,c_{X\bmod D})]$$
$\in
k^\times\backslash
(k\otimes_{k_0}\CC)^\times/
(F^\times)^{Hom_{k_0}(k,\CC)}$.
\endproclaim

\demo{Proof of Corollary}
By Theorem 1 for the dual
$\MM^*$ and by Lemma 1.7,
we have
$per(\MM)\cdot
per(1)^{-\rank\MM}
=
\Gamma(\nabla:\MM^*)\times
[(\det\MM^*,c_{X\bmod D})]^{-1}.$
Since
$\Gamma(\nabla:\MM^*)=
\Gamma(-\nabla:\MM)$
and
$\det\MM^*=(\det\MM)^{\otimes-1}$,
Corollary follows.
\enddemo

When $\dim X\le1$,
we have more precise results.
\proclaim{Lemma 4.1}
Assume $\dim X=0$.
Then the canonical isomorphism
$$
\det\Gamma(X,\EE)\otimes
\det\Gamma(X,\Cal O_X )^{\otimes-\rank\MM}
\simeq
N_{X/k}(\det\EE)$$
induces an isomorphism
$$
\det\Gamma(X,\MM)\otimes
\det\Gamma(X,1)^{\otimes-\rank\MM}
\simeq
N_{X/k}(\det\MM)
$$
in $P_{k_0,F}(k)$.
\endproclaim

Assume $\dim X=1$.
We define an object
$(\det\MM,c_{X\bmod D})$
not only as an isomorphism class.
For a reduced divisor $D\subset X$,
let $\bold{Pic}(X,D)$
be the commutative Picard category of
invertible $\Cal O_X $-modules
$\EE$ with a partial trivialization
$\EE|_D\to \Cal O _D$.
There is an equivalence
$\bold{CH}^1(X,D)\to
\bold{Pic}(X,D)$.
For an object
$f=(f_x)_x\in
\bigoplus_xA_x$
of
$\bold{CH}^1(X,D)$,
its image is the invertible
$\Cal O_X $-module
$\EE=\Cal O_X (\sum_x\ord_xf\cdot[x])$
with the trivialization
$\otimes f_x:\EE(x)\to\kx$
for $x\in D$.
Hence the pairings
$P(U)\times
\bold{CH}^1(X,D)\to
P(k)$ and
$P_{k_0,F,X}(U)\times
\bold{CH}^1(X,D)\to
P_{k_0,F}(k)$ defined in section 3
induce pairings
$P(U)\times
\bold{Pic}(X,D)\to
P(k)$ and
$P_{k_0,F,X}(U)\times
\bold{Pic}(X,D)\to
P_{k_0,F}(k)$ respectively.
In [D ](6.1),
a pairing
$P(X)\times
\bold{Pic}(X)\to
P(k)$ attaching
an invertible $k$-vector space $(\Cal L,\Cal M)$
to a pair of invertible $\Cal O_X $-module $\Cal L,\Cal M$
is defined.
There is a canonical isomorphism
of two pairings
$P(X)\times
\bold{Pic}(X,D)\to
P(k)$ where one is induced by
$P(U)\times
\bold{Pic}(X,D)\to
P(k)$ and
the other is induced by
$P(X)\times
\bold{Pic}(X)\to
P(k)$.
We define the additive functor
$(\ ,c_{X\bmod D})$
to be the pairing with the dual object
of $(\Omega^1_X(log D),\res)$ of $\bold{Pic}(X,D)$.

For an integrable connection
$(\EE,\nabla)$ of rank $r$
on $U$ regular along $D$,
we define a canonical isomorphism
$$
\det{}^\Gamma R\Gamma_c(U,DR(\EE))
\otimes
\det{}^\Gamma R\Gamma_c(U,DR(\Cal O_U ))^{\otimes-r}
\simeq
(\det\EE,c_{X\bmod D}).
$$
Take an extension $(\EX,\nabla)$
of $(\EE,\nabla)$.
By the isomorphism of the pairings above,
we identify
$(\det\EE,c_{X\bmod D})=
(\det\EX(-rD),\Omega^1_X(log D)^{\otimes-1})$.
We define a canonical isomorphism
$$\align
\det R\Gamma(X,DR(\EX))
\otimes&
\det R\Gamma(X,DR(\Cal O_X (-D)))^{\otimes-r}\\
\simeq&
(\det\EX(-rD),\Omega^1_X(log D)^{\otimes-1})
\endalign$$
Here
$DR(\Cal O_X (-D))=
[\Cal O_X (-D)\overset{d}\to\to
\Cal O_X (-D)\otimes\Omega^1_X(log D)]$.
In the notation of [D7] , the left hand side is
$$
\det R\Gamma(X,([\EX]-r[\Cal O_X (-D)])\otimes
([\Cal O_X ]-[\Omega^1_X(log D)])).$$
Therefore by Deligne-Riemann-Roch
(Construction 7.2 loc.cit.),
it is canonically isomorphic to
$(\det\EX(-rD),\Omega^1_X(log D)^{\otimes-1})$.
By the following Lemma,
the isomorphisms for small extensions
$\EX$ define an isomorphism
$
\det{}^\Gamma R\Gamma_c(U,DR(\EE))
\otimes
\det{}^\Gamma R\Gamma_c(U,DR(\Cal O_U ))^{\otimes-r}
\simeq
(\det\EE,c_{X\bmod D}).
$
We write
$\lambda(\EX)=$\linebreak
$\det R\Gamma(X,DR(\EX))$
for short.

\proclaim{Lemma 4.2} Let
$\EX'\subset\EX$ be extensions.
Assume the canonical map
$DR(\EX')\to DR(\EX)$ is a quasi-isomorphism.
Then the diagram
$$\CD
\lambda(\EX')\otimes
\lambda(\Cal O_X (-D))^{\otimes-r}
@>>>
\lambda(\EX)\otimes
\lambda(\Cal O_X (-D))^{\otimes-r}\\
@VVV @VVV\\
(\det\EX'(-rD),\Omega^1_X(log D)^{\otimes-1})
@>>>
(\det\EX(-rD),\Omega^1_X(log D)^{\otimes-1})
\endCD$$
is commutative.
Here the upper horizontal map is
$\det(\res\nabla;\EX/\EX')$-times
the isomorphism induced by the quasi-isomorphism $DR(\EX')\to DR(\EX)$
and the lower one is defined by the isomorphism $\EX'|_U=\EX|_U$.
The vertical arrows are defined by Deligne-Riemann-Roch.
\endproclaim
\demo{Proof}
By induction on the length of
$\EX/\EX'$, we may assume
$\FF=\EX/\EX'$ is an $\kx$-vector space for some $
x\in D$
and $\det(\nabla:\EX/\EX')=
N_{x/k}(\det(\res_x\nabla:\FF))$.
By definition of the isomorphisms,
it is sufficient to
show the commutativity of the diagram
$$\CD
k@>>>
\det R\Gamma(([\EX]-[\EX'])\otimes
([\Cal O_X ]-[\Omega^1_X(log D)]))\\
@| @VV{DRR}V\\
k@>>>
(\det \EX\otimes\det \EX^{\prime\otimes-1},
\Omega^1_X(log D)^{\otimes-1}).
\endCD$$
Here the upper horizontal arrow is
$\det(\nabla:\EX/\EX')$-times
the isomorphism
induced by the isomorphism
$\nabla:\EX/\EX'\to
\EX/\EX'\otimes\Omega^1_X(log D)$
and the lower one is defined by the isomorphism
$\EX|_U=\EX'|_U$.
By the definition of the isomorphism of Deligne-Riemann-Roch
7.1 and 7.2 loc. cit,
the diagram is commutative if we replace the
upper horizontal arrow by the isomorphism
induced by
$id\otimes\res_x^{-1}:\EX/\EX'\to
\EX/\EX'\otimes\Omega^1_X(log D)$.
The assertion follows from this and the equality
$\det(\nabla:\EX/\EX')=
N_{x/k}(\det(\res_x\nabla:\FF))$.
\enddemo

\proclaim{Theorem 2}
Assume $\dim X=1.$
Then the canonical isomorphism of the de Rham component above
induces an isomorphism
$$
\det{}^\Gamma R\Gamma_c(U,\MM)
\otimes
\det{}^\Gamma R\Gamma_c(U,1)^{\otimes-\rank\MM}
\simeq
(\det\MM,c_{X\bmod D}).
$$
in $P_{k_0,F}(k)$.
\endproclaim

We give an outline of proof of Theorems.
First, we deduce Theorem 2 for $X=\PP_k$
from [T] Theorem 1.2.
Then we study the vanishing cycles of a
fibration to a curve in Lemma 4.5.
Using this, we prove Theorem 2 for curves.
Finally using it and taking a Lefshetz pencil,
we prove Theorem 1 by induction on the dimension of $X$.

\demo{Remark}
For an arbitrary dimension,
it seems possible to define a {\it canonical\/} isomorphism
of Deligne-Riemann-Roch
$$\gather
\det R\Gamma(X,
DR(\EX))
\otimes
(\det R\Gamma(X,
DR(\Cal O_X (-D))))^{\otimes-r}\\
\simeq
(\det \EE,c_{X\bmod D}).
\endgather$$
Note that
$ch([DR(\EX)]-r[DR(\Cal O_X (-D))])
=c_1(\EX(-D))\cdot c_n(\Omega^1_X(log D))$
modulo higher terms.
It should induce an isomorphism
$$
\det{}^\Gamma R\Gamma_c(U,\MM)
\otimes
\det{}^\Gamma R\Gamma_c(U,1)^{\otimes-r}
\simeq
(\det\MM,c_{X\bmod D}).
$$
as in Theorem 2 for curves
by the same argument as in the proof of Theorem 1 given later.
\enddemo

The following elementary lemma is repeatedly used in the proof.
\proclaim{Lemma 4.3}
Let $E\subset X$ be a smooth divisor such that
$D'=D\cup E$ has simple normal crossings.
Put $U'=X-D'$ and $E^*=E-D$.
Let $\MM$ be an object of $M_{k_0,F,X}(U)$.
Then

(1) Theorem 1 for two of
$\MM, \MM|_{U'}$ and $\MM|_{E^*}$
implies that for the rest.

(2) If $\dim X=1$,
Theorem 2 for $\MM$ and $\MM|_{U'}$
are equivalent.
\endproclaim
\demo{Proof}
(1). We define a canonical isomorphism
$${\det}^\Gamma R\Gamma_c
(U,\MM)\simeq
{\det}^\Gamma R\Gamma_c
(U',\MM|_{U'})\otimes
{\det}^\Gamma R\Gamma_c
(E^*,\MM|_{E^*}).$$
Let $\MM=((\EE,\nabla),V,\rho)$ and
$(\EX,\nabla)$ be a small extension of
$(\EE,\nabla)$.
Then $\EX(-E)$
is a small extension of
$(\EE|_{U'},\nabla)$ and
$\EX|_E$
is a small extension of
$(\EE|_{E^*},\nabla)$.
Let $DR'_X(\EX(-E))$
be the de Rham complex
$(\EX(-E)\otimes\Omega^\bullet_X(\log D'))$.
Then by the exact sequence
$$
0\to
\Omega^1_X(\log D')(-E)\to
\Omega^1_X(log D)\to
\Omega^1_E(\log (D\cap E))\to 0,$$
we have an exact sequence of de Rham complexes
$$
0\to
DR'_X(\EX(-E))\to
DR_X(\EX)\to
DR_E(\EX|_{E})\to
0.$$
By [K-M] from this we obtain an isomorphism
$$\det R\Gamma_c
(U,\MM)\simeq
\det R\Gamma_c
(U',\MM|_{U'})\otimes
\det R\Gamma_c
(E^*,\MM|_{E^*}).$$
The equality
$\Gamma(\nabla:\EX)
=\Gamma(\nabla:\EX')
\Gamma(\nabla:\EX|_E)$
is easily checked and we obtain a canonical isomorphism
$$\det^\Gamma R\Gamma_c
(U,\MM)\simeq
\det^\Gamma R\Gamma_c
(U',\MM|_{U'})\otimes
\det^\Gamma R\Gamma_c
(E^*,\MM|_{E^*}).$$

On the other hand,
applying Corollary 2 of Lemma 2.3 to
$(\EE,\rho)=$\linebreak
$(\Omega^1_X(log D),\res)$ and
$(\FF,\sigma)=
(\Omega^1_X(\log D'),\res)$,
we  obtain
$c_{X\bmod D}=
c_{X\bmod D'}+
c_{E\bmod D\cap E}$
in $CH^n(X\bmod D)$.
Hence we obtain a canonical isomorphism
$$(\det \MM,c_{X\bmod D})
\simeq(\det \MM|_{U'},c_{X\bmod D'})
\otimes
(\det \MM|_{E^*},c_{E\bmod D\cap E})
.$$
Thus the assertion (1) is proved.

(2).
It is enough to check the commutativity
of the diagram of isomorphisms
$$\matrix
\lambda_X(\EX)\otimes
\lambda_X(\Cal O_X (-D))^{\otimes-r}
&@>>>&
\lambda'_X(\EX(-E))\otimes
\lambda'_X(\Cal O_X (-D'))^{\otimes-r}\\
& &\otimes
\lambda_E(\EE|_E)\otimes
\lambda_E(\Cal O _E)^{\otimes-r}\\
@VVV@VVV\\
(\det \EE,c_{X\bmod D})
&@>>>&
(\det \EE|_{U'},c_{X\bmod D'})\\
& &\otimes
N_{E/k}(\det \EE|_{E^*}).
\endmatrix$$
The vertical isomorphisms in the diagram are
given by the isomorphisms
$$\gather
\det R\Gamma(X,
([\EX]-r[\Cal O_X (-D)])
\otimes
([\Cal O_X ]-[\Omega^1_X(log D)]))\\
\to
(\det\EX(-rD),
\Omega^1_X(log D)^{\otimes-1}),\\
\det R\Gamma(X,
([\EX(-E)]-r[\Cal O_X (-D')])
\otimes
([\Cal O_X ]-[\Omega^1_X(\log D')]))\\
\to
(\det\EX(-rD),
\Omega^1_X(\log D')^{\otimes-1}),\\
\det R\Gamma(X,
([\EX]-r[\Cal O_X (-D)])
\otimes
([\Omega^1_X(\log D')]-[\Omega^1_X(log D)]))\\
\to
(\det\EX(-rD),
\Cal O_X (E))
=N_{E/k}(\EE|_E).
\endgather$$
Since the second isomorphism is the same as
$$ \align \det R\Gamma(X,
([\EX]-r[\Cal O_X(-D) ])
& \otimes
([\Cal O_X ]-[\Omega^1_X(\log D')])) \\
& \to
(\det\EX(-rD),
\Omega^1_X(\log D')^{\otimes-1}),\\ \endalign $$
the diagram is commutative.
\enddemo

We deduce Theorem 2 for
$X=\PP_k$ from [T] Theorem 1.2.
It is reformulated as follows.
Let $D$ be a reduced divisor of $X$.
We may assume
$\infty\notin D$ and
put $D'=D\cup\{\infty\}$ and
$U=X-D'$.
Let
$\MM=((\EE,\nabla),V,\rho)$
be an object of $P_{k_0,F,X}(U)$
and $(\EX,\nabla)$ be an extension of
$(\EE,\nabla)$ to $X$.
We compute the cohomology and the tame symbol of
$\MM$ explicitly.
First we compute the
de Rham cohomology assuming that
$\EX$ is a free $\Cal O_X $-module
$E\otimes_k\Cal O_X $ for $E=\Gamma(X,\EX).$
By
$\Gamma(X,\Omega^1_X(\log D'))\simeq\Cal O _D$
and
$H^1(X,\Cal O_X )=H^1(X,\Omega^1_X(\log D'))=0$,
we have a quasi-isomorphism
$$\beta:R\Gamma(X,DR(\EX))\simeq
[\nabla=\bigoplus_{x\in D}\res_x\nabla:
E\to E\otimes \Cal O _D].$$
By identifying $E=\EX(\infty)$,
we obtain an isomorphism also denoted by $\beta$
$$\beta:\det R\Gamma(X,DR(\EX))\simeq
(\det(E\otimes \Cal O _D))^{\otimes-1}
\otimes
\det(\EX(\infty)).$$

Next we consider the singular cohomology.
Let $A=\Bbb A^1_k\subset X$ and $j:U\to A$
be the open immersion.
We compute
$R\Gamma(A\an,j_!V)$
explicitly by choosing some paths.
For each embedding $k\to \CC$ over $k_0$,
we take a base point $b\in U\an$,
a small disc $\Delta_x\subset X\an$
with center at $x$
for each $x\in D^{\prime an}$,
a base point $b_x$ of the punctured disc
$\Delta^*_x=
\Delta_x-\{x\}$
and  a path $\gamma_x$
connecting $b$ and $b_x$.
We assume that the union
$|\gamma|=
\bigcup_{x\in D^{\prime an}}
|\gamma_x|$
of the underlying sets
$|\gamma_x|$
of paths $\gamma_x$
is contractible.
Then we have a distinguished triangle
$$\to
R\Gamma(A\an,j_!V)\to
R\Gamma(U\an,V)\to
\bigoplus_{x\in D\an}R\Gamma(\Delta^*_x,V)\to.$$
Since the topological space
$U\an$ is homotopically equivalent
to
$\bigcup_{x\in D\an}
(|\gamma_x|\cup\Delta^*_x)$,
we obtain a distinguished triangle
$$\to
R\Gamma(A\an,j_!V)\to
R\Gamma(
\bigcup_{x\in D\an}|\gamma_x|,V)\to
\bigoplus_{x\in D\an}R\Gamma(b_x,V)\to.$$
In other word, we have a quasi-isomorphism
$$\alpha:R\Gamma(A\an,j_!V)\simeq
[\bigoplus\gamma_{x*}:V_b\to
\bigoplus_{x\in D\an}V_{b_x}].$$
Here the isomorphism
$\gamma_{x*}:V_b
\to V_{b_x}$
is defined by the path $\gamma_x$.
By identifying by the isomorphism
$\gamma_{\infty*}:V_b\to V_{b_\infty}$,
we obtain an isomorphism also denoted by $\alpha$
$$\alpha:\det R\Gamma(A\an,j_!V)\simeq
(\det\bigoplus_{x\in D\an}V_{b_x})^{\otimes-1}\otimes
\det V_{b_\infty}.$$

We consider tame symbols.
We fix an embedding $k\to \CC$
and assume $k=\CC$.
By the rational section $dt$
of $\Omega^1_X(\log D')$
for the coordinate $t$ of $X=\PP_k$,
we have
$c_{X\bmod D}=
-(\sum_{x\in D}
(t-t(x))[x]+(-t)[\infty])
\in CH^1(X\bmod D').$
Hence there is an isomorphism
$$(\det V,c_{X\bmod D})\simeq
\bigotimes_{x\in D}
(\det V,t-t(x))_x^{\otimes-1}
\otimes
(\det V,-t^{-1})_\infty.$$
Keeping the notation
$\gamma_x$ etc. above,
we define isomorphisms
$(\det V,t-t(x))_x
\simeq
\det V_{b_x}$
for $x\in D$
and
$(\det V,-t^{-1})_\infty
\simeq
\det V_{b_\infty}$.
To simplify the notation, we replace $V$
by $\det V$ and assume $\rank V=1$.
We take a branch  of the logarithms
$\log (t-t(x))$
for $x\in D$
and $\log t$ on $|\gamma|$
so that
$\log(t(b_\infty)-t(x))-\log(t(b_\infty))=
-\sum_{n=1}^\infty
(t(x)/t(b_\infty))^n/n$.
For $x\in D$,
we put $\nabla_x
=\res_x\nabla\in\CC$.
For $v\in V_{b_x}$,
we take the branch of
$(-1)^{-\nabla_x}
(t-t(x))^{\nabla_x}
v$
of $(V,t-t(x))^\sim$
such that the stalk at $b_x$
is
$\exp((-\pi \sqrt{-1}+
\log (t(b_x)-t(x)))\nabla_x)v.$
We write
$(v,t-t(x))_x\in
(V,t-t(x))_x$
 the stalk of this branch at $x$.
Then the map
$v\mapsto
(v,t-t(x))_x$
defines an isomorphism
$(\ ,t-t(x))_x:V_{b_x}\to
(V,t-t(x))_x.$
Similarly for $\infty$,
an isomorphism
$(\ ,-t^{-1})_{\infty}:
V_{b_\infty}\to
(V,-t^{-1})_{\infty}$
is defined by using the branch with stalk
$\exp(\log t(b_\infty)\nabla_\infty)v$
at $b_\infty.$

Assuming $k_0=k=\CC$,
we define an injective homomorphism
$$(\rho,dt):
((\det\bigoplus_{x\in D}V_{b_x})^{\otimes -1}
\otimes
\det V_{b_\infty})
\to
\det(E\otimes\Cal O _D)^{\otimes -1}
\otimes
\det \EX(\infty)$$
using tame symbols above.
Since
$\det\bigoplus_{x\in D}V_{b_x}
\simeq
\bigotimes_{x\in D}\det V_{b_x}$
and
$\det(E\otimes \Cal O _D)
\simeq
\bigotimes_{x\in D}\det \EX(x)$,
it is enough to define isomorphisms
$$\gather
(\rho ,t-t(x))_x:
\det V_{b_x}
\to
\det \EX(x)\\
(\rho ,-t^{-1})_\infty:
\det V_{b_\infty}
\to
\det \EX(\infty)
\endgather$$
for $x\in D$
and $\infty$ and put
$(\rho,dt)=
\bigotimes_x
(\rho ,-t^{-1})_\infty$.
We define it for $x\in D$
by the commutativity of the diagram
$$\CD
\det V_{b_x}
@>{(\rho ,t-t(x))_x}>>
\det \EX(x)\\
@V{(\ ,t-t(x))_x}VV
@VV{(\ ,t-t(x))_x}V\\
(\det V,t-t(x))_x
@>>>
(\det \EE, t-t(x))_x.
\endCD$$
The lower horizontal arrow is the map induced by the
comparison map $\rho$ defined in section 3,
the right vertical is also defined in the last section
and the left vertical arrow is
the map defined in the previous paragraph.
It is defined similarly for $\infty$.

If the extention
$\EX$ is small at $D$ and big at $\infty$,
the map $\rho$
induces an quasi-isomorphism
$$
R\Gamma(\rho):
R\Gamma(A\an,j_!V)
\to
R\Gamma(X,DR(\EX)).
$$

\proclaim{Theorem T}
Assume
$k_0=k=\CC$,
$\EX\simeq\Cal O_X ^r$
and the real parts of eigenvalues
of $\nabla_x=\res_x\nabla$ are positive for $x\in D$
and negative for $x=\infty$.
Put
$(-1)^{\Tr(\nabla_\infty)}=
\exp(\Tr(\nabla_\infty)\pi\sqrt{-1})$
and
$\Gamma=
\prod_{x\in D}
\Gamma(\res_x\nabla)\cdot
(-1)^{\Tr(\nabla_\infty)}\cdot
\Gamma(1-\res_\infty\nabla)^{-1}$.
Then the diagram
$$\CD
\det R\Gamma(A\an,j_!V)
@>>>
\det R\Gamma(X,DR(\EX))\\
@V{\alpha}VV @VV{\beta}V\\
(\det\bigoplus_{x\in D}V_{b_x})^{\otimes -1}
\otimes
\det V_{b_\infty}
@>>{(\rho,dt)}>
\det(E\otimes\Cal O _D)^{\otimes -1}
\otimes
\det \EX(\infty)
\endCD$$
is commutative.
Here the upper horizontal arrow is
$(-1)^{r(d-1)}\Gamma^{-1}$-times
$\det R\Gamma(\rho)$
and $r=\rank \MM$ and $d=\deg D$.
\endproclaim

\demo{Remark}
The isomorphisms $\alpha$
and $(\rho,dt)$
may depend on the paths
$\{\gamma_i\}_i$,
but the composition is independent.
\enddemo

\demo{Proof}
We may assume $F=\CC$.
We show that the statement is
equivalent to Theorem 1.2 [T].
We define isomorphisms
$$\gather
(\bigotimes_{x\in D}
\det E)^{\otimes -1}
\otimes
\det E
\simeq
(\det\bigoplus_{x\in D}V_{b_x})^{\otimes -1}
\otimes
\det V_{b_\infty}\\
(\bigotimes_{x\in D}
\det E)^{\otimes -1}
\otimes
\det E
\simeq
\det(E\otimes\Cal O _D)^{\otimes -1}
\otimes
\det \EX(\infty)
\endgather$$ as follows.
The second isomorphism is defined
by the canonical isomorphism
$E\to\EX(x)$ for $x\in D'$.
The first one is induced by
the composition of the canonical isomorphisms
$E\to \EX(b)$,
$\rho(b):\EX(b)\to V_b$
and
$\gamma_{x*}:V_b\to V_{b_x}$ for $x\in D'$.
By the isomorphisms,
the maps
$\beta\circ\det R\Gamma(\rho)\circ\alpha^{-1}$
and
$(\rho,dt)$
are regarded as the multiplications by numbers.
Let
$\det R\Gamma(\rho)$
and $(\rho,dt)$
also denote the corresponding numbers respectively.
We prove
$$\align
(-1)^{r(\deg D-1)}
\det R\Gamma(\rho)
\times
\det (-\res_\infty\nabla)=&
 \det H_c\\
(-1)^{-\Tr(\nabla_\infty)}
\Gamma
\times
\det (-\res_\infty\nabla)=&
 \prod_{i=1}^n
\Gamma_{\lambda_i}(P)
\Gamma_{\infty}(-P)^{-1}\\
(-1)^{\Tr(\nabla_\infty)}(\rho,dt)=&
 \prod_{i=1}^n
(P,x-\lambda_i)_{\gamma_i}
\times(P,1/x)_{\gamma_\infty}^{-1}.
\endalign$$
Here the notation in the right hand sides
is that in Theorem 1.2 loc.cit.
The equalities imply that the statements are equivalent.

We give a brief dictionary between the notation here and
loc.cit.
The logarithmic integrable connection
$\nabla:\EX\to \EX\otimes\Omega^1_X(log D)$
here corresponds to
$\partial_P:\Cal O (V^*)_{\log}\to
\Omega^1(V^*)_{\log}$
defined by
$\partial_Pf=df+fP$
in Section 1.2 loc.cit.
The global section
$E=\Gamma(X,\EX)$ here is $V^*$ there.
The coordinate $t$ here is denoted by
$x$ there and the value $t(x_i)$ for
$x_i\in D$ is $\lambda_i\in \CC$ there.
The residues
$\nabla_i=\res_{x_i}\nabla$
at $x_i\in D$ here is the right multiplication
by the matrix $B^{(i)}:V^*\to V^*$.

We reformulate the definition of $\det H_c$.
Fix a numbering
$D=\{x_1,\cdots,x_n\}$.
We define isomorphisms
$E^{n-1}\simeq
H^1(X,DR(\EX))$
and
$E^{*n-1}\simeq
H^1(A\an,j_!V)^*$.
Here ${}^*$
denotes the dual.
For de Rham cohomology,
we have defined a quasi-isomorphism
$$\beta:R\Gamma(X,DR(\EX))
\simeq
[\bigoplus_x\res_x\nabla:E\to\bigoplus_{x\in D}E].$$
Since
$\sum_{x\in D}\res_x\nabla=
-\res_\infty\nabla:E\to E$
is an isomorphism,
$\ssum :\bigoplus_{x\in D}E\to E$
induces a retraction of
$\bigoplus_x\res_x\nabla:E\to\bigoplus_{x\in D}E$.
Hence
$H^q(X,DR(\EX))=0$ for $q\neq1$ and an isomorphism
$\beta':H^1(X,DR(\EX))\to
\Ker(\ssum :\bigoplus_{x\in D}E\to E)$.
For singular cohomology, we have
$$
\alpha:R\Gamma(A\an,j_!V)\simeq
[\bigoplus_x\gamma_{x*}:V_b\to \bigoplus_{x\in D}V_{b_x}]
\simeq
[\text{diagonal}:E\to\bigoplus_{x\in D}E].
$$
Hence we have
$H^q(A\an,j_!V)=0$ for $q\neq1$ and
$H^1(A\an,j_!V)=
\Coker(\text{diagonal}:E\to\bigoplus_{x\in D}E)$.
Taking the dual, we obtain an isomorphism
$\alpha^*:H^1(A\an,j_!V)^*=
\Ker(\ssum :\bigoplus_{x\in D}E^*\to E^*)$.

We identify
$E^{n-1}
\simeq\Ker(\ssum :E^n\to E)$
by $id\times 0-0\times id$.
Here
$id\times 0:
E^{n-1}\to
E^{n-1}\times E=E^n$
and
$0\times id:
E^{n-1}\to
E\times E^{n-1}=E^n$.
Similarly we identify
$E^{*n-1}
\simeq\Ker(\ssum :E^{*n}\to E^*)$
by $id\times 0-0\times id$.
Thus we have defined isomorphisms
$E^{n-1}\simeq
H^1(X,DR(\EX))$ and
$E^{*n-1}\simeq
H^1(A\an,j_!V)^*.$
The isomorphism
$H^1(\rho):
H^1(A\an,j_!V)^*\to
H^1(X,DR(\EX))$
induces a perfect pairing
$$\gather
E^{*n-1}\times
E^{n-1}
\simeq
H^1(A\an,j_!V)^*\times
H^1(X,DR(\EX))\\
@>{id\times H^1(\rho)^{-1}}>>
H^1(A\an,j_!V)^*\times
H^1(A\an,j_!V)
@>{can}>>\CC.
\endgather$$
Then the determinant
$\det H_c$
is characterized  by the condition that the induced pairing
$\det (E^{*n-1})\times
\det (E^{n-1})\to\CC$
is $\det H_c$-times the canonical pairing.

We prove the first equality.
The definition of $\det H_c$ above is rephrased as follows.
We have quasi-isomorphisms
$$\align
\beta':R\Gamma(X,DR(\EX))\simeq&
[\ssum:\bigoplus_{x\in D}E\to E]\\
\alpha^*:R\Gamma(A\an,j_!V)^*\simeq&
[\ssum:\bigoplus_{x\in D}E^*\to E^*].
\endalign$$
They induce isomorphisms
$$\align
\beta':\det R\Gamma(X,DR(\EX))\simeq&
\det E^{\otimes 1-n}\\
\alpha^*:\det R\Gamma(A\an,j_!V)^*\simeq&
\det E^{*\otimes 1-n}.
\endalign$$
The isomorphism
$R\Gamma(\rho):
R\Gamma(A\an,j_!V)\to
R\Gamma(X,DR(\EX))$
induces a perfect pairing
$$\gather
(\det E^*)^{\otimes 1-n}\times
(\det E)^{\otimes 1-n}
\simeq
\det R\Gamma(A\an,j_!V)^*\times
\det R\Gamma(X,DR(\EX))\\
@>{id\times \det R\Gamma(\rho)^{-1}}>>
\det R\Gamma(A\an,j_!V)^*\times
\det R\Gamma(A\an,j_!V)
@>{can}>>\CC.
\endgather$$
Then
$\det H_c$
is characterized by the condition that the
pairing is
$(\det H_c)^{-1}$-times the canonical one.
Since the isomorphism
$\alpha^*:\det R\Gamma(A\an,j_!V)^*\simeq
\det E^{*\otimes 1-n}$
is the dual of
$\alpha:
\det R\Gamma(A\an,j_!V)\simeq
\det E^{\otimes 1-n}$,
we have
$$\det H_c=
\beta'\circ
\det R\Gamma(\rho)\circ
\alpha^{-1}=
(\beta'\circ\beta^{-1})
\det R\Gamma(\rho).$$
Hence the first equality will follow from
$\beta'\circ\beta^{-1}=
(-1)^{r(n-1)}\det(-\res_{\infty}\nabla).$
This is the consequence of Lemma below
applied to $E'=E^n$,
$r=\bigoplus_{x\in D}\res_x\nabla$
and $s=\ssum$.

\proclaim{Lemma 4.4}
Let $r:E\to E'$
and $s:E'\to E$
be homomorphisms of
vector spaces $E$ and $E'$ of finite dimensions
such that $s\circ r$ is an automorphism of $E$.
Consider the quasi-isomorphism $b$
$$\CD
    @.  [E'  @>s>> E]\\
@.       @VV{id}V  @.\\
[E @>r>> E'] @.
\endCD$$
of complexes where $E'$ is put on degree 1.
Then the automorphism $\det b$
of $\det E\otimes (\det E')^{\otimes-1}$
is the multiplication by
$(-1)^{\dim E(\dim E'-\dim E)}\times \det (s\circ r)^{-1}$.
\endproclaim
\demo{Proof}
Since the determinant of $(s\circ r, id)
:[E\overset {r}\to \to E']\to
[E\overset {r(s\circ r)^{-1}}\to \to E']$
is $\det s\circ r$, we may assume $s\circ r =1$.
Then the assertion follows from Koszul's rule on the sign.
\enddemo

The second equality follows immediately from
the definition, the dictionary
$\res_{x_i}\nabla=B^{(i)}$
and from the formula
$(-s)\Gamma(1-s)=\Gamma(-s).$

We prove the third equality.
By taking the determinant,
we may assume the rank $r$ is 1.
By identification
$E\simeq \EX(b)
\overset{\rho_b}\to\to V_b
\overset{\gamma_{x*}}\to\to V_{b_x}$
and $E\simeq \EX(x)$
for $x\in D'$ above,
we regard
$(\rho,t-t(x))_x$
for $x\in D$,
$(\rho,-t^{-1})_\infty$
etc. as numbers.
By the definition
$(\rho,dt)=
\prod_x(\rho,t-t(x))^{-1}_x
\times(\rho,-t^{-1})_\infty$
and by the residue formula
$\sum_{x\in D,\infty}\nabla_x=0$,
it is sufficient to show that
$$\gather
(\rho,t-t(x_i))_{x_i}
=(-1)^{-\nabla_x}(P,x-\lambda_i)_{\gamma_i}^{-1}\\
(\rho,-t^{-1})_\infty
=(P,1/x)^{-1}_{\gamma_\infty}
\endgather$$
for $x_i\in D$ and $\infty$.
We compute the left
hand sides following the definition.
For $x\in D,=\infty$ and $e\in E$,
let $v$ be the (multivalued) section
of $V$ on $\Delta^*_x$ such that stalk
$v$ at $b_x$ is $\gamma_{x*}(\rho_b(e))$.
Let $\varphi$ be the analytic function
$\rho(v)/e$ on $\Delta^*_x$.
Then for $x\in D$
by the definition in Section 3,
$(\rho,t-t(x))_x=(\varphi,t-t(x))_x$
is the value at $x$ of the branch of the invertible function
$(-1)^{-\nabla_x}(t-t(x))^{\nabla_x}h$
taking the value
$\exp(\nabla_x(-\pi\sqrt{-1}+
\log(t(b_x)-t(x))))
\rho(v)/e$
at $b_x$.
On the other hand,
the notation
$$(P,x-\lambda_i)_{\gamma_i}=
\lim_{x\to\lambda_i}
D(x)\cdot(x-\lambda_i)^{-tr B^{(i)}}$$
loc.cit Section 1.1 Definition
is translated here as follows.
We have
$D(x)=F(x)$ there is $\varphi(x)^{-1}$ here,
$tr B^{(i)}=B^{(i)}$ there is $\nabla_{x_i}$ here and
$(x-\lambda_i)^{-tr B^{(i)}}$ there is
$\exp(-\nabla_{x_i}\log (t-t(x_i)))$.
Hence we have
$$(\rho,t-t(x_i))_{x_i}
=\exp(-\nabla_{x_i}\pi\sqrt{-1})(P,x-\lambda_i)_{\gamma_i}^{-1}$$
Similarly
$(\rho,-t^{-1})_\infty=(\varphi,-t^{-1})_\infty$
is the value at $\infty$ of the branch of the invertible function
$t^{-\nabla_\infty}\varphi$
taking the value
$\exp(-\nabla_\infty
\log(t(b_\infty)))\rho(v)/e$
at $b_\infty$
and
$(P,1/x)_{\gamma_\infty}=
\lim_{x\to\infty}
D(x)x^{trB^{(\infty)}}$.
Hence by a similar dictionary, we have
$$
(\rho,-t^{-1})_\infty
=(P,1/x)^{-1}_{\gamma_\infty}.$$
Thus the third equality and hence Theorem T is proved.
\enddemo

\demo{Proof of Theorem 2 for $X=\PP_k$}
Changing the coordinate if necessary,
we always assume $\infty\in U$.
Put $U'=U-\{\infty\}$.
Let
$\MM
=((\EE,\nabla),V,\rho)$
be an object of
$M_{k_0,F,X}(U')$.
For an extension $\EX$ of $\EE$,
we put
$\Gamma'(\nabla:\EX)=
\prod_{x\in D}N_{x/k}\Gamma(\res_x\nabla)
\times(-1)^{\Tr(\nabla_\infty)}
\Gamma(1-\res_\infty\nabla)^{-1}
\in(k\otimes_{k_0}\CC)^\times$.
We define a variant
$\det^\Gamma R\Gamma'
(U',\MM)$
to be the triple
$$
\projlim_{(\EX,\nabla)}
(\det R\Gamma
(X,DR(\EX)),
\det R\Gamma
(A\an,j_!V),
\Gamma'(\nabla:\EX)^{-1}\times
\det R\Gamma(\rho))$$
$\in P_{k_0,F}(k)$
similarly as in the definition of
$\det^\Gamma R\Gamma_c
(U',\MM)$.
Here the extensions
$(\EX,\nabla)$ of
$(\EE,\nabla)$
are small at $D$ and big at $\infty$
and the transition morphism
for $\EX'\subset \EX$ is
$(\det(\nabla:\EX/\EX'))$-times the isomorphism
induced by the inclusion.
By the equality
$s(-1)^{-(s-1)}\Gamma(1-(s-1))=
(-1)^{-s}\Gamma(1-s)$,
they form a projective system.

We show Theorem T implies

\proclaim{Claim 1}
Assume $\infty\in U$.
Let
$\MM
=((\EE,\nabla),V,\rho)$ be an object of
$M_{k_0,F,X}(U)$.
Assume that there exists an
extension $(\EX,\nabla)$ of
$(\EE,\nabla)$
such that
$\EX\simeq\Cal O_X (-m)^r$ for some $m\in \ZZ$.
Then the isomorphism of Deligne-Riemann-Roch
$$\align
\det R\Gamma(X,DR(\EX))
\otimes&
\det R\Gamma(X,DR(\Cal O_X (-D)))^{\otimes-r}\\
\simeq&
(\det\EX(-rD),\Omega^1_X(log D)^{\otimes-1})
\endalign$$
induces an isomorphism
$$
\det{}^\Gamma R\Gamma'(U',\MM|_{U'})
\otimes
\det{}^\Gamma R\Gamma'(U',1)^{\otimes-r}
\simeq
(\det\MM|_{U'},c_{X\bmod D'})
$$
in $P_{k_0,F}(k)$.
\endproclaim

We complete the proof
of Theorem 2 for $X=\PP_k$
admitting Claim 1.
We show that the assumption of Claim 1
is not restrictive.

\proclaim{Claim 2}
Shrinking $U$
if necessary,
there exists an extension $(\EX,\nabla)$
such that $\EX\simeq\Cal O_X (m)^r$ for some $m\in\ZZ$.
\endproclaim
\demo{Proof of Claim 2}
By a theorem of Grothendieck,
a locally free $\Cal O_X $-module $\EX$ is
isomorphic to $\bigoplus_j\Cal O_X (m_j)$
for some $(m_j)_{j=1,\cdots,r}$.
We prove by induction on $M=
\sum_j(m_j-m)$
where $m=\min_jm_j$.
If it is 0, there is nothing to prove.
Otherwise, take $j_0$
such that $m_{j_0}\neq m$
and a $k$-rational point $x\in U$.
Then the kernel $\EX'$ of the composite
$\EX\to\EX(x)\overset{pr_{j_0}}\to\to
\Cal O_X (m_{j_0})(x)$
is an extension of
$\EE|_{U-x}$
and $M$ for $\EX'$ is $M-1$.
Thus Claim is proved.
\enddemo

We show that the conclusion of Claim 1 implies
that of Theorem 2 for $X=\Bbb P^1$.
Let
$\MM
=((\EE,\nabla),V,\rho)$ be an object of
$M_{k_0,F,X}(U)$
and $\EX$ be an extension of $\EE$ as in the assumption of Claim 1.
It is sufficient to define canonical isomorphisms
$$\gather
{\det}^\Gamma R\Gamma_c
(U,\MM)
\simeq
{\det}^\Gamma R\Gamma'
(U',\MM|_{U'})
\otimes
\det (\MM|_\infty(-1))\\
(\det\MM,c_{X\bmod D})
\simeq
(\det\MM|_{U'},c_{X\bmod D'})
\otimes
\det (\MM|_\infty)
\endgather$$
satisfying the following condition.
Write
$$
\align
\lambda(\EX) &=
\det R\Gamma(X,
[\EX\overset\nabla\to\to
\EX\otimes\Omega^1_X(log D)])
\\
\lambda'(\EX) &=
\det R\Gamma(X,
[\EX\overset\nabla\to\to
\EX\otimes\Omega_X(\log D')])
\endalign $$
for short.
The condition is that de Rham components
of the isomorphisms form a commutative diagram
$$\CD
\lambda(\EX)
\otimes
\lambda(\Cal O_X (-D))^{\otimes-r}
&{\to}&
\lambda'(\EX)
\otimes
\lambda'(\Cal O_X (-D))^{\otimes-r}
 \otimes\det\EX(\infty)
\\
@V{DRR}VV   @VV{DRR}V\\
(\det \EX(-rD),
\Omega^1_X(log D)^{\otimes-1})
&{\to}&
(\det \EX(-rD),
\Omega_X(\log D')^{\otimes-1})
\otimes\det\EX(\infty)
\endCD$$
for small extensions $(\EX,\nabla)$.
Note that
$\Gamma(\nabla:\EX)=
\Gamma'(\nabla:\EX)$
since $\res_\infty\nabla=0$.
We define the isomorphisms.
The first one is defined
by the distinguished triangles
$$\gather
\to
R\Gamma
(X,DR(\EX))
\to
R\Gamma
(X,DR'(\EX))
\to
\EX(\infty)[-1]\to\\
\to
R\Gamma_c(U\an,V)
\to
R\Gamma_c(U^{\prime an},V|_{U^{\prime an}})
\to
V_\infty(-1)[-1]\to
\endgather$$
The second one is a consequence of the equality
$c_{X\bmod D}=c_{X\bmod D'}+[\infty]$ in Lemma 2.1.
The commutativity is
proved in the same way as in Lemma 4.2.
Thus the proof of Theorem 2 for $X=\PP_k$
is reduced to the proof of Claim 1.

\demo{Proof of Claim 1}
We assume
$\EX\simeq
\Cal O_X (-m)^r$
for some $m$.
Replacing $\EX$ by
$\EX(-m'D)$ for large $m'\in \Bbb N$
if necessary,
we may assume that the following condition is satisfied.
We put $\EX'=\EX(m[\infty])$
so that $\EX'\simeq\Cal O_X ^r$.
Then the real parts of all the conjugates in $\CC$
over $k_0$ of the eigenvalues of
$\res_x\nabla\in\kx$ for $x\in D$ are positive
and those for $\res_\infty\nabla_{\EX'}$
are negative.
The condition is also satisfied
for the extension $\Cal O_X (-D)$ of $(\Cal O_U ,d)$
and $\Cal O_X '=\Cal O_X (-D+\deg D[\infty])$.
It is sufficient to show that for each embedding
$k\to \CC$ over $k_0$,
there is an isomorphism
$$\det R\Gamma(A\an,j_!V)
\otimes
\det R\Gamma(A\an,j_!F)^{\otimes-r}
\to
(\det V\otimes F^{\otimes -r},c_{X\bmod D'})
$$
of $F$-vector spaces such that the diagram
$$\matrix
\det R\Gamma(A\an,j_!V)
\otimes
&@>>>&
\det R\Gamma(X,DR(\EX'))
\otimes
\\
\quad\det R\Gamma(A\an,j_!F)^{\otimes-r}
& &
\quad\det R\Gamma(X,DR(\Cal O_X '))^{\otimes-r}\\
@VVV
@VV{DRR}V
\\
(\det V\otimes F^{\otimes -r},c_{X\bmod D'})
&@>>>&
(\det\EX'\otimes
\det\Cal O_X ^{\prime\otimes-r},
\Omega^1_X(\log D')^{\otimes-1}).
\endmatrix
$$
is commutative.
Here the upper horizontal arrow is
$\Gamma'(\nabla:\EX)^{-1}\times$\linebreak
$\Gamma'(\nabla:\Cal O_X (-D))^r\times
\det R\Gamma (\rho)\otimes
\det R\Gamma (\text{canonical})^{\otimes -r}$
and the lower one is the comparison map of
$(\det \MM|_{U'},c_{X\bmod D'})$.

By the definition of the tame symbol given in section 3,
we have a commutative diagram
$$\CD
(\det\bigoplus_{x\in D\an}V_{b_x})^{\otimes-1}
\otimes
\det V_{b_\infty}
\otimes
&@>>>&
\det(\EX'|_D)^{\otimes-1}
\otimes
\det(\EX'(\infty))
\otimes
\\
\quad((\det\bigoplus_{x\in D\an}F)^{\otimes-1}
\otimes F)^{\otimes -r}
& @. &
\quad(\det(\Cal O_X '|_D)^{\otimes-1}
\otimes
\det(\Cal O_X '(\infty)))^{\otimes-r}
\\
@VVV @. @VVV
\\
(\det V\otimes F^{\otimes -r},c_{X\bmod D'})
&@>>>&
(\det\EX'\otimes
\det\Cal O_X ^{\prime\otimes-r},
\Omega^1_X(\log D')^{\otimes-1}).
\endCD
$$
Here the upper horizontal arrow is
$(\det\rho,dt)\otimes(\det \text{ can },dt)^{\otimes-r}$
and the lower one is
$(\det\rho\otimes \text{ can}^{\otimes-r},c_{X\bmod D'})$.
The left vertical arrow is induced by
the tame symbols
$(\ ,t-t(x))_x$
for $x\in D\an$ and
$(\ ,-t^{-1})_\infty$
and the right one is defined by
a property of the norm
$
\det(\EX'|_D)
\otimes
\det(\Cal O_X '|_D)^{\otimes-r}
\simeq
N_{D/k}(\det(\EX'|_D))
$
and the isomorphism
$$\gather
N_{D/k}(\det(\EX'|_D))
\otimes
(\det(\EX'(\infty))
\otimes
\det(\Cal O_X '(\infty))^{\otimes-r})\\
\to
(\det\EX'\otimes
\det\Cal O_X ^{\prime\otimes-r},
\Omega^1_X(\log D')^{\otimes-1}).
\endgather
$$
The isomorphism of Deligne-Riemann-Roch is
the composition of this and
$$\gather
\det R\Gamma(X,DR(\EX'))
\otimes
\det R\Gamma(X,DR(\Cal O_X '))^{\otimes-r}\\
@>{\beta_{\EX'}\otimes\beta_{\Cal O_X '}^{\otimes-r}}>>
\det(\EX'|_D)^{\otimes-1}
\otimes
\det(\EX'(\infty))
\otimes
(\det(\Cal O_X '|_D)^{\otimes-1}
\otimes
\det(\Cal O_X '(\infty)))^{\otimes-r}.
\endgather$$
Therefore by Theorem T for $\EX'$ and $\Cal O_X '$,
the composite of
$\alpha_V\otimes\alpha_F^{\otimes-r}$
and the left vertical arrow of the diagram above
gives the required isomorphism of the $F$-vector spaces
$\det R\Gamma(A\an,j_!V)
\otimes
\det R\Gamma(A\an,j_!F)^{\otimes-r}
\to
(\det V\otimes F^{\otimes -r},c_{X\bmod D'}).
$
Thus Claim 1 and hence Theorem 2
for $X=\PP_k$ are proved.
\enddemo
\enddemo

To prove the general case,
we study a fibration $f:X\to Y$
to a curve $Y$.
We consider a local situation below.
Let $Y$ be a smooth curve over $k$ and
$E$ be a reduced divisor containing a closed point $y$.
Let $X$ be a smooth scheme over $k$,
$D$ be a divisor with simple normal crossings
and $f:X\to Y$ be a proper and flat morphism.
We assume that
\roster
\item"$\star$"
The underlying set of
$f^*(E)$ is in $D$.
On the complement
$W=Y-E$,
the restriction
$f_W:X_W\to W$
is smooth and $D_W=D\cap X_W$
is a divisor with relative normal crossings.
\endroster
Let $(\EE,\nabla)$
be an integrable connection on
$U=X-D$
and $(\EX,\nabla)$
be a small extension to $X$.
Let
$\Omega^1_{X/Y}(\log D/\log E)$
be the cokernel of the natural morphism
$f^*\Omega^1_Y(\log E)\to
\Omega^1_X(log D)$.
It is a locally free $\Cal O_X $-module of rank
$\dim X-1$.
Let
$\Omega^q_{X/Y}(\log D/\log E)=
\bigwedge^q\Omega^1_{X/Y}(\log D/\log E)$.
For a logarithmic integrable connection
$(\EX,\nabla)$,
we define the relative de Rham complex
$$\align
&DR_{X/Y}(\EX)=
(\EX\otimes
\Omega^\bullet_{X/Y}(\log D/\log E))\\
=&
[\EX
\overset\nabla\to\to
\EX\otimes
\Omega^1_{X/Y}(\log D/\log E)
\overset\nabla\to\to
\cdots].
\endalign$$
On the higher direct image
$\FY^q=
R^qf_*DR_{X/Y}(\EX)$,
the Gauss-Manin connection
$\nabla:
\FY^q\to
\FY^q\otimes
\OmYE$
is defined as follows.
By the definition of
$DR_{X/Y}(\EX)$,
we have an exact sequence of complexes
$$
0\to
DR_{X/Y}(\EX)\otimes f^*\OmYE[-1]
\to
DR_X(\EX)\to
DR_{X/Y}(\EX)\to0.
$$
The connection $\nabla$ is defined as the
boundary map of the induced long exact sequence
$$
\to
R^qf_*DR_X(\EX)\to
\FY^q\overset\nabla\to\to
\FY^q\otimes\OmYE\to.
$$

We generalize the terminology for logarithmic integrable connections
on locally free $\Cal O _Y$-modules defined in section 1
to those on coherent modules.
Let $(\FF,\nabla)$ be an integrable connection on $W$ .
Let $(\FY,\nabla)$ be the pair of a coherent $\Cal O _Y$-module
$\FY$  and a logarithmic integrable connection
$\nabla:\FY\to\FY\otimes\OmYE$
extending $(\FF,\nabla)$.
The connection $\nabla$
induces a connection on the torsion part
$\FF_{Y,tors}\to
\FF_{Y,tors}\otimes\OmYE$.
We say an extension $\FY$
is small if $\FY/\FF_{Y,tors}$
is small in the sense of section 1
and if
$\FF_{Y,tors}\to
\FF_{Y,tors}\otimes\OmYE$
is an isomorphism.
If $\FY$ and $\FY'$ are small extensions
of $\FF$,
a morphism
$(\FY,\nabla)\to(\FY',\nabla)$
extending the identity on $\FF$ induces a
quasi-isomorphism
$DR(\FY)\to DR(\FY')$.
An extension $\FY$
is small if and only if
$DR(\FY)\an=j_{W!}\an DR(\FF)\an$.

We define the characteristic rational function
$\Phi_{\FF_{Y,y}}(T)\in
\ky(T)^\times$ as follows.
The $\ky$-vector spaces
$T_i=Tor^{\Cal O _{Y,y}}_i(\FF_{Y,y},\ky)$
for $i=0,1$
are the cokernel and the kernel of the map
$\FY(-E)_y\to\FF_{Y,y}$
respectively.
Hence the connection
$\nabla$ induces an
$\ky$-linear endomorphisms
$$\res_{i,y}\nabla:
T_i\to T_i\otimes
\OmYE_y
@>{id\otimes\res_y}>>
T_i.
$$
We define
$$\Phi_{\FF_{Y,y}}(T)=
\prod_{i=0,1}
\det(T-\res_{i,y}\nabla:
Tor^{\Cal O _{Y,y}}_i(\FF_{Y,y},\ky))^{(-1)^i}.
$$
When $Y$ is proper, we  put
$\Phi_{\FY}(T)=
\prod_{y\in E}
N_{y/k}
(\Phi_{\FF_{Y,y}}(T))
\in k(T)^\times$
and
$\Gamma(\nabla:\FY)=
(\Gamma(\Phi_{\FY}(T)))
\in (k\otimes_{k_0}\CC)^\times$
as in section 1.
For an object $\MM=((\FF,\nabla),V,\rho)$ of
$M_{k_0,F,Y}(W)$,
the determinant
$\det^\Gamma
R\Gamma_c
(W,\MM)$
is canonically isomorphic to the triple
$$(\det R\Gamma(Y,DR(\FY)),
\det R\Gamma_c(U,V),
\Gamma(\nabla:\FY)^{-1}\cdot\det R\Gamma_c(\rho))$$
for a small extension $\FY$ of $\FF$.

\proclaim{Lemma 4.5}
Let $f:X\to Y$ and $(\EX,\nabla)$
be as above.
Then
\roster
\item
For each $q$,
the Gauss-Manin connection
$(\FY^q,\nabla)$
is small.
\item
For $y\in E$,
the characteristic rational functions
$\Phi^q(T)=
\Phi_{\FY^q,y}(T)$
satisfy
$$\prod_{D_i\subset f^{-1}(y)}
N_{k_i/\ky}
\left(
\prod_{\ell=0}^{m_i-1}
\frac{\Phi_{\EX,i}(m_iT-\ell)}{m_i^{\rank\EE}}
\right)^{c_i}=
\prod_q\Phi^q(T)^{(-1)^q},$$
where $k_i$ is the constant field of $D_i$ and $m_i$ is the multiplicity
of $D_i$ in $f^{-1}(y)$.
\endroster
\endproclaim
\demo{Proof}
(1).
It is sufficient to show that
$DR(\FY^q)\an\simeq
j_{W!}\an DR(\FF^q)\an.$
By the long exact sequence
$$
\to
R^qf_*DR_X(\EX)\to
\FY^q\to
\FY^q\otimes\OmYE\to
$$
and
$$R^qf_*DR_X(\EX)\an
\simeq
R^qf_*j_!DR_U(\EE)\an
\simeq
j_{W!}DR_W(\FF^q)\an,$$
we see
$DR(\FY^q)\an|_{E\an}=0$.
Thus the assertion is proved.

(2).
By the exact sequence
$$0\to
\FY^q\otimes\ky\to
H^q(X_y,DR_{X/Y}(\EX))\to
Tor_1^{\Cal O _Y}(\FY^{q+1},\ky)
\to0,$$
the right hand side is equal to
$$\align
&\det(T-\res_y\nabla:
R\Gamma(X_y,DR_{X/Y}(\EX)))\\:=&
\prod_q
\det(T-\res_y\nabla:
H^q(X_y,DR_{X/Y}(\EX)))^{(-1)^q}.
\endalign$$
We prove
$\det(T-\res_y\nabla:
R\Gamma(X_y,DR_{X/Y}(\EX)))$
is equal to the left hand side.
We take a sequence of divisors
$0=D_0<D_1<\cdots<D_\ell=\sum_i\ell_iD_i<\cdots<D_M=X_y$
such that $D_{\ell+1}-D_\ell=D_{i(\ell)}$ is irreducible.
The map $[0,M)\to\amalg_{\in I_y}[0,m_i)$
defined by $\ell\mapsto
\ell_{i(\ell)}\in[0,m_{i(\ell)})$
is a bijection.
We have
$$\gather
\det(T-\res_y\nabla:
R\Gamma(X_y,DR_{X/Y}(\EX)))\\=
\prod_{\ell=0}^{N-1}
\det(T-\res_y\nabla:
R\Gamma(D_{i(\ell)},DR_{X/Y}(\EX(-D_\ell))|_{D_{i(\ell)}})).
\endgather$$
Since
$\Phi_{i,\EX(-D_\ell)}(T)=
\Phi_{i,\EX}(T-\ell_i),$
it is sufficient to show that
$$\det(T-\res_y\nabla:
R\Gamma(D_i,DR_{X/Y}(\EX(-D_\ell))|_{D_i}))=
N_{k_i/\ky}
\left(
\frac{\Phi_{\EX(-D_\ell),i}(m_iT)}{m_i^{\rank\EE}}
\right)^{c_i}$$
for $D_i\subset f^{-1}(y)$.
Changing the notation, we write $\EX$ for $\EX(-D_\ell)$.

We show that
$\res_y\nabla$ on
$R\Gamma(D_i,DR_{X/Y}(\EX)|_{D_i})$
is induced by
$m_i^{-1}\res_i\nabla$
on $\EX|_{D_i}$.
\proclaim{Lemma 4.6}
Let $0\to A@>a>> B@>b>> C\to 0$
be an exact sequence of complexes.
Let $r:B\to A$
be a componentwise retraction and $s:C\to B$
be the corresponding section: $1_B=ar+sb$.
Then the canonical map
$\tilde b:K=\roman{cone }[A\to B]\to C$
admits a homotopy inverse
$\tilde s=(s,rds):C^p\to B^p\oplus A^{p+1}$.
The composite
$C\to K\to A[+1]$
is given by $rds:C^p\to A^{p+1}$.
\endproclaim

\demo{Proof}
A homotopy connecting
$1_B-\tilde s\circ\tilde b$
is given by $K^p\to B^p\overset{r}\to\to A^p\to K^{p-1}$.
The rest is straightforward to check.
\enddemo

We apply Lemma to the exact sequence
$$0\to
f^*\OmYE\otimes
DR_{X/Y}(\EX)|_{D_i}[-1]
\to
DR_X(\EX)|_{D_i}
\to
DR_{X/Y}(\EX)|_{D_i}
\to
0$$
and the retraction
$$\align
&m_i^{-1}d\log\pi\otimes id\otimes \res_i:
\EX\otimes
\Omega^q_X(\log D)|_{D_i}\\
\to&\
f^*\OmYE\otimes
\EX\otimes
\Omega^{q-1}_{X/Y}(\log D/\log E)|_{D_i}
\endalign$$
for a prime element $\pi$ of $\Cal O _{Y,y}$.
The section $s$ is induced by the canonical map
$$
\Omega^q_{X/Y}(\log D/\log E)|_{D_i}
\simeq
\Omega^q_{D_i}(\log D|_{D_i})
\to
\Omega^q_X(\log D)|_{D_i}
$$
where
$D|_{D_i}$
is the divisor $\bigcup_{j\ne i}
D_j\cap D_i$ with simple normal crossings of $D_i$.
By an elementary computation,
the map $rds$ in Lemma 4.6 is
$$\gather
m_i^{-1}d\log\pi\otimes\res_i\nabla\otimes id:
\EX\otimes
\Omega^\bullet_{X/Y}(\log D/\log E)|_{D_i}\\
\to f^*\OmYE\otimes
\EX\otimes
\Omega^\bullet_{X/Y}(\log D/\log E)|_{D_i}.
\endgather$$
Therefore by Lemma 4.6,
$\res_y\nabla$
on $R\Gamma(D_i,DR_{X/Y}(\EX)|_{D_i})$
is induced by
$m_i^{-1}\res_i\nabla$ on
$\EX|_{D_i}$.

To complete the proof of (2),
it is sufficient to show
$$\det(T-\res_i\nabla
:R\Gamma(D_i,DR_{X/Y}(\EX)|_{D_i}))
=N_{k_i/k}
\Phi_{\EX,i}(T)^{c_i}.$$
This follows from Lemma below applied for
$X=D_i, \EE=\EX|_{D_i}, f=\res_i\nabla$
and $\FF=\Omega_{X/Y}^1(\log D/\log E)|_{D_i}
\simeq\Omega^1_{D_i}(\log D|_{D_i})$.

\proclaim{Lemma 4.7}
Let $X$ be a proper smooth scheme
of dimension $n$ over a field $k$
and let $k'$ be the constant field of $X$.
Let $\EE$ be a torsion-free coherent
$\Cal O_X $-module and $f$ be an endomorphism of $\EE$.
Then
\roster
\item
The eigenpolynomial
$\Phi(T)=
\det (T-f:\EE(\xi))\in \kappa(\xi)[T]$
of $f$ at the generic point
$\xi$ has coefficient in $k'$.
\endroster
Further let $\FF$
be a locally free $\Cal O_X $-module
of rank $n$
and $c$ be the degree
$\deg :CH^n(X)\to CH^0(k')=\ZZ$
of the top chern class $c_n(\FF)$.
Then
\roster
\item"{(2)}"
We have
$$\prod_{q=0}^n
\det (T-f\otimes id:R\Gamma(X,\EE\otimes
\wedge^q\FF))^{(-1)^q}
=N_{k'/k}\Phi(T)^c$$
in $k(T)^\times$.
\endroster
\endproclaim

\demo{Proof}
We may assume $k'=k$.
To show  (1),
it is enough to prove the minimal polynomial
$\Phi_0(T)\in\kappa(\xi)[T]$
of $f$ at $\xi$
divides a non-zero polynomial in $k[T]$.
The minimal polynomial of
$f\in End_{\Cal O_X }(\EE)$
over $k$
is in $k[T]$ and divisible by $\Phi_0(T)$.
We show (2).
By taking a direct sum decomposition,
we may assume $\Phi(T)$
is a power of an irreducible polynomial
$\Phi_1(T)\in k[T]$.
Then the both sides are powers of
$\Phi_1(t)$
and it is enough to compare the degrees.
By Riemann-Roch
$$\gather
\sum_{q=0}^n
(-1)^q\chi(X,\EE\otimes
\wedge^q\FF)
=\left((ch(\EE)\cdot\sum_q(-1)^q
ch(\wedge^q\FF))\cdot td\Omega_X^{1*}
\right)_{\dim 0}\\
=\rank \EE\cdot c_n(\FF),
\endgather$$
Lemma 4.7 is proved.
\enddemo
Thus Lemma 4.5 is proved.
\enddemo

\proclaim{Corollary}
For $D_i\in f^{-1}(y)$,
let $\nabla_i=\Tr \ \res_i\nabla
\in k_i$
and $c_i$ be the Euler number
of $D_i^*\otimes_{k_i}{\bar k_i}$.
We put $m_i^{\nabla_i}=
\exp(\nabla_i\log m_i)
\in (k_i\otimes_{k_0}\CC)^\times$ where
$\log m_i\in \Bbb R$.
Then we have
$$\gather
\prod_{D_i\subset
f^{-1}(y)}
N_{k_i/\ky}
\left(\frac
{\Gamma_{D_i}(\nabla:\EX)}
{\Gamma_{D_i}(\nabla:\Cal O_X (-D))^r}
\right)^{c_i}\\
=\prod_q
\left(\frac
{\Gamma_y(\nabla:R^qf_*DR_{X/Y}(\EX))}
{\Gamma_y(\nabla:R^qf_*DR_{X/Y}(\Cal O_X (-D)))^r}
\right)^{(-1)^q}
\times
\prod_{D_i\subset
f^{-1}(y)}
N_{k_i/\ky}m_i^{\nabla_ic_ir}
\endgather$$
in
$(\ky\otimes_{k_0}\CC)^\times$.
\endproclaim
\demo{Proof}
We may assume $k=\CC$.
By the formula
$$
\frac
{\Gamma(s)}
{\Gamma(1)}
=
\prod_{\ell=0}^{m-1}
\frac
{\Gamma((s+\ell)/m)}
{\Gamma((1+\ell)/m)}
\cdot m^s,$$
we have
$$
\Gamma\left(\frac
{\Phi(T)}{(T-1)^r}\right)
=
\Gamma\left(
\prod_{\ell=0}^{m-1}
\frac
{\Phi(mT-\ell)}{(mT-(1+\ell))^r}\right)
\cdot m^{\Sigma(\Phi(T))}
$$
for a polynomial $\Phi(T)$ of degree $r$.
For $\Phi(T)=
\sum_{i=0}^ra_iT^i=
a_r\prod_{i=1}^r(T-\alpha_i)$,
we put
$\Sigma(\Phi(T))=-a_{r-1}/a_r=
\sum_{i=1}^r\alpha_i.$

Corollary follows immediately from the definitions,
the formula and Lemma 4.5 (2).
\enddemo


\demo{Proof of Theorem 2}
Let $f:X\to Y$
be a finite flat morphism of proper smooth curves over $k$.
Let $W$ be an open subscheme of $Y$ such that
on the inverse image $U=f^*(W)$ the restriction
$f|_U:U\to W$ is etale.
We put $D=X-U$ and $E=Y-W$.
We have a canonical isomorphism
$f^*\OmYE\simeq\Omega^1_X(log D)$.
For an extension $(\EX,\nabla)$ to $X$
of an integrable connection
$(\EE,\nabla)$ on $U$ of rank $r$,
the canonical isomorphism
$f_*DR_X(\EX)\simeq DR_Y(f_*\EX)$
induces an isomorphisms
$$
\det R\Gamma(X,DR(\EX))
\simeq
\det R\Gamma(Y,DR(f_*\EX)).
$$
The canonical isomorphisms
$$N_{X/Y}(\det\EX\otimes
\Cal O_X (-D)^{\otimes-r})
\simeq
\det f_*\EX
\otimes
\det f_*\Cal O_X (-D)^{\otimes -r}$$
and
$f^*\OmYE\simeq\Omega^1_X(log D)$
induce an isomorphism
$$
\align
(\det \EX(-rD),\Omega^1_X(log D))\simeq&
(N_{X/Y}\det \EX(-rD),\OmYE)\\
\simeq&
(\det f_*\EX
\otimes
\det f_*\Cal O_X (-D)^{\otimes -r},
\OmYE).
\endalign$$
They are compatible with the isomorphisms
of Deligne-Riemann-Roch and the following diagram is commutative
$$\CD
\det\Gamma(X,DR(\EX))\otimes
& @.  &
\det\Gamma(Y,DR(f_*\EX))\otimes \\
\det\Gamma(X,DR(\Cal O_X (-D)))^{\otimes-r}
&@>>>&
\det\Gamma(Y,DR(f_*\Cal O_X (-D)))^{\otimes-r}\\
@V{DRR}VV @. @VV{DRR}V\\
(\det \EX(-rD),\Omega^1_X(log D))
&@>>>&
(\det f_*\EX
\otimes
\det f_*\Cal O_X (-D)^{\otimes -r},
\OmYE).
\endCD$$

Let $\MM=((\EE,\nabla),V,\rho)$
be an object of
$M_{k_0,F,X}(U)$
and $(\EX,\nabla)$ be an extension of $(\EE,\nabla)$
as above.
Let
$m_{\EX}^\nabla
=(\exp(\sum_{x\in D}\sigma
(\Tr_{x/k}(\Tr(\res_x\nabla)))\times\log m_x))_\sigma
\in (k\otimes_{k_i}\CC)^\times$ where
$m_x$ is the ramification index at $x$ and $\log m_x\in \Bbb R$
and let
$[m_{\EX}^\nabla]$
be the corresponding object of
$P_{k_0,F}(k)$.
For extensions $\EX\subset \EX'$,
the multiplication by
the inverse of
$\prod_{x\in D}
m_x^{\dim_k(\EE'_{X,x}/\EE_{X,x})}$
on $k$ induces an isomorphism
$[m_{\EX}^\nabla]
\simeq [m_{\EX'}^\nabla]$.
Define an object
$[m_{\MM}^\nabla]$
of $P_{k_0,F}(k)$ to
be the projective limit of
$[m_{\EX}^\nabla]$
for small extensions $\EX$.

\proclaim{Claim}
The canonical isomorphisms above
induce isomorphisms
$$\gather
{\det}^\Gamma
R\Gamma_c(U,\MM)
\otimes
{\det}^\Gamma
R\Gamma_c(U,1)^{\otimes-r}\\
\simeq
{\det}^\Gamma
R\Gamma_c(W,f_*\MM)
\otimes
{\det}^\Gamma
R\Gamma_c(W,f_*1)^{\otimes-r}
\otimes
[m^\nabla_{\MM}],\\
(\det \MM,c_{X\bmod D})
\simeq
(\det f_*\MM\otimes
\det f_*1^{\otimes-r},
c_{Y\bmod E})
\otimes
[m_{\MM}^\nabla].
\endgather$$
\endproclaim
We show Claim implies Theorem 2.
Take a finite flat morphism $f:X\to Y=\PP_k$.
Shrinking $U$ if necessary,
we assume that $f|_U$ is etale and
$U=f^*f(U)$.
By Theorem 2 for $Y=\PP_k$
already proved,
the isomorphism of Deligne-Riemann-Roch
defines an isomorphism on
the right hand sides of the isomorphisms in Claim.
Therefore by the commutativity above,
the isomorphism of Deligne-Riemann-Roch
also defines an isomorphism
on the left hand sides of the isomorphisms.
Thus Theorem 2 is proved.

\demo{Proof of Claim}
We show the first isomorphism.
By Lemma 4.5 (1), the isomorphism above
for a small extension $(\EX,\nabla)$
induces an isomorphism
$$
\det R\Gamma_c(U,\MM)
\simeq
\det R\Gamma_c(W,f_*\MM).
$$
Therefore it is sufficient to show an equality
$$
\Gamma_X(\nabla:\EX)
\Gamma_X(\nabla:\Cal O_X (-D))^{-r}
=
\Gamma_Y(\nabla:f_*\EX)
\Gamma_Y(\nabla:f_*\Cal O_X (-D))^{-r}
\times m_{\EX}^\nabla$$
in $(\CC^\times)^{Hom_{k_0}(k,\CC)}$.
It is an immediate consequence of Corollary of Lemma 4.5.

We show the second isomorphism.
By Lemma 2.1 applied to
$f^*\OmYE\simeq\Omega^1_X(log D)$
and
$f^*\res_y=m_x^{-1}\res_x$
for $x\in D$ and
$y=f(x)\in E$,
we have
$$c_{X\bmod D}=
f^*c_{Y\bmod E}-
\sum_{x\in D}m_x[x].$$
Hence there is a canonical isomorphism
$$(\det \MM,c_{X\bmod D})\simeq
(N_{X/Y}\det \MM
,c_{Y\bmod E})\otimes
\bigotimes_{x\in D}
(\det \MM,m_x)_x^{\otimes-1}.$$
By definition, there is a natural isomorphism
$\bigotimes_{x\in D}
(\det \MM,m_x)_x^{\otimes -1}
\simeq
[m_{\MM}^\nabla]$.
Further by the canonical isomorphism
$N_{X/Y}\det \MM
\simeq
\det f_*\MM\otimes
\det f_*1^{\otimes-r},
$
the right hand sides of the above isomorphism
and of the Claim are canonically isomorphic.
It is easy to check that
the de Rham component of the canonical isomorphism
defined here is the same as that given there.
Thus Claim is proved.
\enddemo
\enddemo

\demo{Proof of Theorem 1}
We prove it by induction on $\dim X$.
It is already proved for $\dim X\le 1$.
First we consider the case
\roster
\item"$\star$"
There is a proper flat morphism
$f:X\to Y$
to a proper smooth curve satisfying the following condition.
For $W=f(U)$, the restriction $f_W:f^{-1}W\to W$ is smooth and
$D\cap f^{-1}(W)$ has relative normal crossings.
\endroster
This is equivalent to the condition $\star$ before.
We put $E=Y-W$.
Let $(\EX,\nabla)$ be an extension
of an integrable connection
$(\EE,\nabla)$ on $U$.
Then we have a spectral sequence
$$E_2^{p,q}=
H^p(Y,DR_YR^qf_*DR_{X/Y}(\EX))
\Rightarrow
H^{p+q}(X,DR_X(\EX)).$$
In fact, the exact sequence of complexes
$$
0\to
DR_{X/Y}(\EX)
\otimes f^*\OmYE[-1]
\to
DR_X(\EX)\to
DR_{X/Y}(\EX)
\to0$$
induces a quasi-isomorphism
$$
DR_X(\EX)\to
[DR_{X/Y}(\EX)\to
DR_{X/Y}(\EX)
\otimes f^*\OmYE].$$
Taking the direct image, we get
$$
Rf_*DR_X(\EX)
\to
[Rf_*DR_{X/Y}(\EX)\to
Rf_*DR_{X/Y}(\EX)
\otimes \OmYE].$$
The canonical filtration on
$Rf_*DR_{X/Y}(\EX)$
induces an increasing filtration
$\tau$ on
$Rf_*DR_X(\EX)$
such that
$$\align
Gr^\tau_qRf_*DR_X(\EX)
\simeq&
[R^qf_*DR_{X/Y}(\EX)\to
R^qf_*DR_{X/Y}(\EX)
\otimes \OmYE]\\
=&
DR_YR^qf_*DR_{X/Y}(\EX)[-q].
\endalign$$
{}From this, we obtain the spectral sequence above
by [D2] (1.4.8) and taking its decalage.
For a small extension $\EX$,
by Lemma (1),
the spectral sequence gives an isomorphism
$$
\det R\Gamma_c(U,DR(\EE))
\simeq
\det R\Gamma_c(W,Rf_!DR_{U/W}(\EE)).
$$
Here we put
$Rf_!DR_{U/W}(\EE)=
Rf_*DR_{X/Y}(\EX)|_W$
for a small extension $\EX$.

Let $I_1\subset I$
be the index set of the components
$f(D_i)\subset E$
and $I_2$ be the complement $I-I_1$.
Similarly as in the proof of Theorem 2,
we put
$m_{\EX}^\nabla
=\prod_{i\in I_1}(c_i\cdot
\exp(\sigma
(\Tr_{k_i/k}(\Tr(\res_i\nabla)))\cdot\log m_i))_\sigma
\in (k\otimes_{k_0}\CC)^\times$ where
$m_i$ is the multiplicity of
$D_i$ in the fiber $f^*{f(D_i)}$
and $\log m_i\in \Bbb R$.
We also put
$\Gamma_2(\nabla:\EX)=
\prod_{i\in I_2}N_{k_i/k}
\Gamma(\res_i\nabla:\EX)^{c_i}
\in (k\otimes_{k_0}\CC)^\times$.
The objects
$[m_{\MM}^\nabla]=
\projlim[m_{\EX}^\nabla]$
and
$[\Gamma_2(\nabla:\MM)]=
\projlim[\Gamma_2(\nabla:\EX)]$
are defined as the limit
for small extensions $\EX$ as before.

\proclaim{Claim}
The canonical isomorphisms above
induce isomorphism
$$\gather
{\det}^\Gamma
R\Gamma_c(U,\MM)
\otimes
{\det}^\Gamma
R\Gamma_c(U,1)^{\otimes-r}\\
\simeq
{\det}^\Gamma
R\Gamma_c(W,Rf_!\MM)
\otimes
{\det}^\Gamma
R\Gamma_c(W,Rf_!1)^{\otimes-r}
\otimes
[m_{\MM}^\nabla]
\otimes
[\Gamma_2(\nabla:\EX)].
\endgather$$
There is an isomorphism
$$
(\det \MM,c_{X\bmod D})
\simeq
(\det Rf_!\MM\otimes
\det Rf_!1^{\otimes-r},
c_{Y\bmod E})
\otimes
[m_{\MM}^\nabla]
\otimes
[\Gamma_2(\nabla:\EX)].
$$
\endproclaim
By Theorem 2,
the right hand sides of isomorphisms
in Claim are isomorphic to each other.
Hence Claim implies Theorem 1
in this case.

\demo{Proof of Claim}
We show the first isomorphism.
The isomorphism above
for a small extension $(\EX,\nabla)$
induces an isomorphism
$
\det R\Gamma_c(U,\MM)
\simeq
\det R\Gamma_c(W,f_*\MM).
$
Therefore it is sufficient to show an equality
$$\gather
\Gamma_X(\nabla:\EX)
\Gamma_X(\nabla:\Cal O_X (-D))^{-r}\\
=
\Gamma_Y(\nabla:Rf_*DR_{X/Y}(\EX))
\Gamma_Y(\nabla:Rf_*DR_{X/Y}(\Cal O_X (-D)))^{-r}
\times m_{\EX}^\nabla
\times
\Gamma_2(\nabla:\EX)
\endgather$$
in $(k\otimes_{k_0}\CC)^\times$.
It is an immediate consequence of Corollary of Lemma 4.5.

We show the second isomorphism.
By Lemma 2.1
for partial trivializations
$\res$ and $m_i^{-1}\res_i$ of
$\Omega^1_X(log D)$ and by Lemma 2.2
applied to the exact sequence
$$0\to f^*\OmYE\to\Omega^1_X(log D)\to
\Omega^1_{X/Y}(\log D/\log E)\to 0,$$
we have
$$c_{X\bmod D}=
f^*c_{Y\bmod E}\cup
c_{X/Y\bmod D/E}
-
\sum_{i\in I_1}m_i
i_{i*}c_{D_i}.$$
Since $(\det \MM,
\sum_{i\in I_1}m_i
i_{i*}c_{D_i})=
[m_{\MM}^\nabla]$
by the definition of the tame symbol,
there is a canonical isomorphism
$$(\det \MM,c_{X\bmod D})\simeq
(\det \MM,
f^*c_{Y\bmod E}\cup
c_{X/Y\bmod D/E})
\otimes
[m_{\MM}^\nabla].$$
Therefore it is sufficient to define an isomorphism
$$\gather
(\det Rf_!\MM\otimes
\det Rf_!\bold1^{\otimes-r},
c_{Y\bmod E})
\otimes
[\Gamma_2(\nabla:\MM)]
\\
\simeq
(\det \MM,
f^*c_{Y\bmod E}\cup
c_{X/Y\bmod D/E}).
\endgather$$

\proclaim{Lemma 4.8}
There is a dense open subset $W'\subset W$
such that for $y\in W'$
the restriction
$\EX|_{X_y}
=\EX\otimes_{\Cal O_X }
\Cal O _{X_y}$
is also reflexive.
\endproclaim

We finish the proof of Claim and
Theorem 1 under the assumption $\star$
admitting Lemma 4.8.
By Lemma 4.8 and Corollary of Proposition 1 in section 2,
we may write
$c_{Y\bmod E}=
\sum_{y\in W'}n_y[y]$.
Then by proper base change and by the definition of
$[\Gamma_2(\nabla:\MM)]$,
the left hand side is
$$\bigotimes_{y\in W'}
N_{y/k}\left(
\det{}^\Gamma R\Gamma_c(U_y,\MM|_{U_y})\otimes
\det{}^\Gamma R\Gamma_c(U_y,1)^{\otimes-r}
\right)^{n_y}.$$
In fact,
by
$c_i=\sum_yn_y\deg y\cdot c'_i$
where $c'_i$ is the Euler number of a geometic fiber
$D_i^*|_W\to W$ for $i\in I_2$,
we have
$\Gamma_2(\nabla:\EX)=
\prod_{y\in W'}
N_{y/k}
\Gamma(\nabla:\EX|_{X_y})^{n_y}$
for a small extension $\EX$.
On the other hand,
$$\gather
(\det \MM,
f^*c_{Y\bmod E}\cup
c_{X/Y\bmod D/E})\\
\simeq\bigotimes_{y\in W'}
N_{y/k}(\det
\MM|_{U_y},c_{X_y\bmod D_y})^{n_y}.
\endgather$$
In fact
for the Gysin map
$i_{y*}:
CH^{n-1}(X_y\bmod D_y)\to
CH^n(X\bmod D)$,
we have
$f^*c_{Y\bmod E}\cup
c_{X/Y\bmod D/E}=
\sum_{y\in W'}
n_y i_{y*}(c_{X_y\bmod D_y})$
and
$(\det\MM,i_{y*}c)=
(\det\MM|_{U_y},c)$.

Therefore by the assumption of the induction
applied for $X_y$,
Theorem 1 is proved under $\star$ and Lemma 4.8.

\demo{Proof of Lemma 4.8}
Let ${}^*$ denote the dual.
Let $W'$ be a dense open subset such that
for any closed point $y\in W'$
the fiber $X_y$
does not contain any component of
the associated cycles
$Ass(\Cal Ext^1_{\Cal O_X }
(\EX,\Cal O_X ))$
and
$Ass(\Cal Ext^1_{\Cal O_X }
(\EX^*,\Cal O_X ))$.
We show for $y\in W'$
the canonical morphisms
$\EX^{*}|_{X_y}
\to
(\EX|_{X_y})^*$
and
$\EX^{**}|_{X_y}
\to
(\EX^*|_{X_y})^*$
are isomorphisms.
This proves Lemma since then the canonical map
$\EX|_{X_y}
\to
(\EX|_{X_y})^{**}$
is also an isomorphism and
$\EX|_{X_y}$
is reflexive.
We show the assertion for $\EX$
and the same argument works for $\EX^*$.
For a closed point $y\in W'$, by the exact sequence
$$0\to
\EX^{*}|_{X_y}
\to
(\EX|_{X_y})^*
\to
\Cal Ext^1_{\Cal O_X }
(\EX,\Cal O (-X_y)),$$
and by the assumption that
$\Cal Ext^1_{\Cal O_X }
(\EX,\Cal O (-X_y))$
has no submodule
supported on the fiber $X_y$,
the canonical morphism
$\EX^{*}|_{X_y}
\to
(\EX|_{X_y})^*$
is an isomorphism.
Thus Lemma is proved.
\enddemo
\enddemo

We prove the general case by
taking a Lefshetz pencil.
We follow Step 2 and Step 3 in the
proof of [S] Theorem 1.
By Lemma 3 loc.cit,
there exists a Lefschetz pencil $(H_y)_{y\in \PP_k}$
of $X$
satisfying the following properties (1)-(4).
\roster
\item
The axis of the pencil
meets transversally $D_J=\bigcap_{i\in J} D_i$
for all $J\subset I$.
\item
There exists a dense open subscheme
$W\subset Y=\PP_k$
such that for all $y\in W$,
the hyperplane $H_y$
meets transversally
$D_J$ for all $J\subset I$.
\item
For all $y\in Y$
and all $J\subset I$,
the hyperplane $H_y$
meets transversally $D_J$
except at isolated ordinary quadratic singularities.
\endroster
For $J\subset I$,
let $S_J=\{x\in D_J;
\text{ $x$ is a singular point of
$H_y\subset D_J$ for some $y\in Y$\}}$.
Then
\roster
\item"{(4)}"
$S_J$ and $S_{J'}$ are disjoint
for all $J\ne J',\subset I$.
\endroster
We define a blowing-up $X_3$ of $X$,
a morphism $f_3:X_3\to Y=\PP_k$
and a divisor $D_3$ with simple normal crossings
satisfying the condition $\star$ above as follows.
Let $X_1$ be the blowing-up of $X$ with the center the intersection
$Z$ of $X$ with the axis of the pencil
and $f_1:X_1\to Y=\PP_k$
be the map defined by the pencil.
Let $S\subset X_1$ be the disjoint union
$S=\amalg_{J}S_J$.
Let $X_2$ be the blowing up of
$X_1$ at $S$.
For $x\in S_J$
and $J\ne\emptyset$,
let $B_x$
be the intersection of the
exceptional divisor $E_x$
and the proper transform of
$D_J$ in $X_2$.
Let $X_3$
be the blowing-up of $X_2$
at all $B_x$
for $x\in S-S_\emptyset$.
We put $E=f_1(S)$,
$W=Y-f_1(S)$
and let $D_3$
be the divisor
$\varphi_3^*D\cup f_3^*E$
where
$f_3:X_3\to Y$
and $\varphi_3:X_3\to X$.
Then it is easy to check that $D_3$
has simple normal crossings and
$f_3:X_3\to Y$
satisfies the condition $\star$.

We have proved Theorem for $(X_3,D_3)$ above.
Therefore together with Lemma,
the following Lemma completes the proof of Theorem 1.

\proclaim{Lemma 4.9}
Let $X,D$ and
$\MM$ be as in the statement of Theorem 1.
Assume Theorem 1 holds for lower dimensions.
Let $Z$ be a regular closed subscheme
of $X$ satisfying one of the following conditions
(1)-(3) and let
$\varphi:X'\to X$ be the blowing-up at $Z$.
Then Theorem 1 for
$\MM$ on $U$
is equivalent to that for
$\MM'=\varphi^*\MM$ on $U'=\varphi^*(U)$.
\roster
\item
$Z$ is a regular closed subscheme of
codimension 2 meeting transversally
$D_J$ for all $J\subset I$.
\item
$Z=\{x\}$ for a closed point $x$ of $X$.
\item
$Z=D_J$ for some $J\subset I$
\endroster
\endproclaim

\demo{Proof} Let $D'$ be the divisor $\varphi^*D=X'-U'$
with simple normal crossing.

(1). Let $D_Z=D\cap Z$ and $Z^*=Z-D_Z$.
It is easy to see that there is an isomorphism
$${\det}^\Gamma R\Gamma_c(U',\varphi^*\MM)
\simeq
{\det}^\Gamma R\Gamma_c(U,\MM)
\otimes
{\det}^\Gamma R\Gamma_c(Z^*,\MM|_{Z^*}(-1)).$$
We show
$(\det \varphi^*\MM,c_{X'\bmod D'})_{X'}
\simeq
(\det \MM,c_{X\bmod D})_X
\otimes(\det \MM|_{Z^*},c_{Z\bmod Z^*})_Z$.
We apply Corollary 2 of Lemma 2.3,
to the exact sequence
$$0\to
\varphi^*\Omega^1_X(log D)\to
\Omega^1_{X'}(\log D')\to
\Omega_{E/Z}^1\to 0.$$
By the exact sequence
$$\align
0\to
N_{E/X}(=\Cal O _E(1))\to
\text{Image }&
(\varphi^*\Omega^1_X(log D)|_E\to
\Omega^1_{X'}(\log D')|_E)
\\
&\to
\varphi^*
\Omega^1_Z(\log D_Z)\to 0,
\endalign$$
we obtain
$c_{X'\bmod D'}
=\varphi^*c_{X\bmod D}
+i_{E_*}(\varphi^*
c_{Z\bmod D_Z}\cup c_1(\Cal O _E(1)))$
in $CH^n(X'\bmod D')$.
By writing an element of the relative Chow group
by 0-cycles (Corollary of Proposition 1 in section 2),
we have isomorphisms
$(\det \varphi^*\MM,\varphi^*c)_{X'}=
(\det \MM,c)_X$
 for $c\in CH^n(X\bmod D)$,
$(\det \varphi^*\MM,
i_{E_*}c)=
(\det \MM|_{E\cap U'},c)$
 for $c\in CH^{n-1}(E\bmod D'\cap E)$ and
$(\det \varphi^*\MM,
i_{E_*}(\varphi^*c\cup c_1(\Cal O _E(1))))=
(\det \MM|_{Z^*},c)$
 for $c\in CH^{n-2}(Z\bmod D_Z)$.
By the isomorphisms and by the equality above,
we obtain the isomorphism.
Under Theorem 1 for $Z$,
(1) is proved.

(2) for $x\in U$.
It is easy to see that there is an isomorphism
$${\det}^\Gamma R\Gamma_c(U',\varphi^*\MM)\simeq
{\det}^\Gamma R\Gamma_c(U,\MM)
\otimes
\bigotimes_{q=1}^{n-1}N_{x/k}(\det(\MM(-q))(x)).$$
We show
$(\det \varphi^*\MM,c_{X'\bmod D'})_{X'}
\simeq
(\det \MM,c_{X\bmod D})_X
\otimes N_{x/k}\det \MM(x)^{\otimes(n-1)}$.
By Corollary 1 of Lemma 2.3
applied to
the exact sequence
$$0\to
\varphi^*\Omega^1_X(log D)\to
\Omega^1_{X'}(\log D')\to
\Omega_E^1\to 0$$
and by an elementary computation,
we obtain
$c_{X'\bmod D'}
=\varphi^*c_{X\bmod D}
+i_{E_*}((n-1)[x'])$
in $CH^n(X'\bmod D')$
for a $\kx$-rational point $x'\in E$.
Similarly as in (1),
we obtain an isomorphism.

(2) for $x\in D$ and (3).
In these cases,
$\det^\Gamma R\Gamma_c(U',\varphi^*\MM)=
\det^\Gamma R\Gamma_c(U,\MM)$.
In fact except for the case (2)
where $\Card I_x=1$,
the Euler numbers $c_i$ of $D_i^*$ do not change
and that for the exceptional divisor $E^*$ is 0.
In the exceptional case, the Euler number of
$c_{i_0}$ for $x\in D_{i_0}$ decreases by 1,
that for $E^*$ is 1 and the characteristic polynomials
$\Phi_{\EX}$
for $D_{i_0}$ and for $E$ are the same.
We show
$(\det \varphi^*\MM,c_{X'\bmod D'})_{X'}
\simeq
(\det \MM,c_{X\bmod D})_X$.
Let $D_1$ be a regular family of
subschemes $D\cup\{Z\}$,
$\rho_1$
be the partial trivialization of
$\Omega^1_X(log D)$ defined by
$\res_i$ for $D_i$
and $\sum_{D_i\supset Z}\res_i$ for $Z$
and let
$\rho'_1$
be the partial trivialization of
$\varphi^*\Omega^1_X(log D)$ on $D'$ defined by
pulling-back $\rho_1$.
We consider the canonical maps
$$
H^n(X,\Kn(X\bmod D))
\gets
H^n(X,\Kn(X\bmod D_1))
\to
H^n(X',\Kn(X'\bmod D')).
$$
The first one is an isomorphism since
$H^{n-1}(Z,\Kn(Z))=0$ and
$H^n(Z,\Kn(Z))=0$ by $\dim Z\le n-2$.
The images of
the relative top chern class
$c_n(\Omega^1_X(log D),\rho_1)$ are
$c_n(\Omega^1_X(log D),\res)$
and $c_n(\varphi^*\Omega^1_X(log D),\varphi^*\rho_1)$ respectively.
Further by Lemma 2.4
applied to
$\varphi^*\Omega^1_X(log D)
\to
\Omega^1_{X'}(\log D)$,
we have
$c_n(\varphi^*\Omega^1_X(log D),\varphi^*\rho_1)=$\linebreak
$c_n(\Omega^1_{X'}(\log D),\res)$.
Thus we obtain an isomorphism
$(\det \varphi^*\MM,c_{X'\bmod D'})_{X'}
\simeq
(\det \MM,c_{X\bmod D})_X$
by Corollary of Proposition 1
and Lemma 4.9 is proved.
\enddemo
Now proof of Theorem 1 is completed.
\enddemo

We have variants of main results
in the category defined below.
Let $k_0$ be a subfield of $\CC$ and
$k$ and $F$ be finite extensions of $k_0$.
Let $U$ be a smooth scheme over $k$.
For an integrable connection
$\nabla:\EE\to \EE
\otimes \OmU$,
we call a multiplication  by $F$ on $(\EE,\nabla)$
a $k_0$-homomorphism
$F\to End_{\Cal O_X }^\nabla(\EE)=
\Ker (\nabla\text{ on }End_{\Cal O_X }(\EE)).$
Let
$M'_{k_0,F}(U)$
be the category of triples
\roster
\item
An integrable connection $(\EE,\nabla)$
on $U$ regular along boundary with a multiplication by $F$.
\item
A local system $V$ of $F$-vector spaces on $U\an$.
\item
A morphism $\rho:V\to \EE\an$ on $U\an$
inducing an isomorphism $V\otimes_{k_0}\CC
\to\Ker \nabla\an$ of local systems
of $F\otimes_{k_0}\CC$-modules.
\endroster
For this category,
we have analogues of
Theorems 1 and 2.

We define briefly
the both sides of the isomorphisms there.
Note that an integrable connecion
$(\EE,\nabla)$ with a multiplication by $F$
is naturally identified with an
integrable connection denoted
$(\EE^{(F)},\nabla^{(F)})$
on $U\otimes_{k_0}F$.
Let $X$ be a proper smooth scheme over $k$
including $U$ as the complement of
a divisor $D$ with simple normal crossings.
For an integrable connection
$(\EE,\nabla)$ on $U$
with a multiplication by $F$,
it is regular along $D$
if and only if $(\EE^{(F)},\nabla^{(F)})$
is regular along $D\otimes_{k_0}F$.
Hence if it is regular,
there exists an extension
$(\EX,\nabla)$ to $X$ with
an extended multiplication
$F\to End^\nabla_{\Cal O_X }(\EX)$.
The rank of an object
$\MM=((\EE,\nabla),V,\rho)$
of $M'_{k_0,F}(U)$
is the rank of $\EE^{(F)}$
as an $\Cal O _{U\otimes F}$-module
which is the same as the rank of
$V$ as a local system of $F$-vector spaces.
The determinant ${\det}_F(\EE,\nabla)$
over $F$ of an integrable connection
$(\EE,\nabla)$ with a multiplication by $F$
is defined as that corresponding to
$\det_{\Cal O _{U\otimes F}}\EE^{(F)}$.
The class group
$MPic'_{k_0,F}(U)$
of the category of rank 1 object
$P'_{k_0,F}(U)$ is defined similarly as in Section 1.
For a finite extension $k$ of $k_0$,
there is a canonical isomorphism
$MPic'_{k_0,F}(k)=
(k\otimes_{k_0}F)^\times\backslash
(k\otimes_{k_0}F\otimes_{k_0}\CC)^\times/
(F^\times)^{Hom_{k_0}(k,\CC)}$.

For an object
$\MM=((\EE,\nabla),V,\rho)$
of $M'_{k_0,F}(U)$,
the determinant of cohomology
${\det}_F
R\Gamma_c(U,\MM)
\in
P'_{k_0,F}(k)$
is defined similarly as in Section 1.
It is the alternating tensor product of
the determinant over $F$ of the compact support cohomology
$$
H^q_c(U,\MM)=
(H^q(X,DR(\EX)),
H^q_c(U\an,V),
H^q_c(\rho))$$
$\in
M'_{k_0,F}(k)$
where $\EX$ is a small extension of
$\EE$ with a multiplication by $F$.
We also define a variant
${\det}_F^\Gamma
R\Gamma_c(U,\MM)$
as an object of
$P'_{k_0,F}(k)$.
Take a small extension $(\EX,\nabla)$
with an extended multiplication by $F$.
We define $\Gamma_F(\nabla:\EX)
\in (k\otimes F\otimes \CC)^\times$
as follows.
By $k\otimes F=\prod_j K_j$ for fields $K_j$,
the fiber product $X\otimes_{\QQ}F$
is the disjoint union $\amalg X_j$.
Let $\EE_{X_j}$ be the restriction on
$X_j$ of the corresponding connection $\EX^{(F)}$.
Then $\Phi_{\EE_{X_j}}(T)\in K_j(T)^\times$
and $\Gamma(\nabla:\EE_{X_j})=
\Gamma(\Phi_{\EE_{X_j}}(T))
\in (K_j\otimes_{\QQ}\CC)^\times$
are defined in Section 1.
We put
$\Gamma_F(\nabla:\EX)
=(\Gamma(\nabla:\EE_{X_j}))_j
\in \prod_j(K_j\otimes_{\QQ}\CC)^\times
= (k\otimes F\otimes \CC)^\times$.
Now similarly as in Section 1, we define
${\det}_F^\Gamma
R\Gamma_c(U,\MM)$ to be
the inverse limit of the triples
$$
({\det}_F R\Gamma(X,DR(\EX)),
{\det}_F R\Gamma_c(U\an,V),
\Gamma_F(\nabla:\EX)^{-1}
\det R\Gamma_c(\rho))$$
for small extensions $\EX$.

We also define a pairing
$$
P'_{k_0,F,X}(U)\times CH^n(X\bmod D)
\to P'_{k_0,F}(k).$$
As is shown in Section 3, the essential point is
the definition of the tame symbol.
We consider locally and do not assume $X$ is proper.
Let $T=D_J$
and define a pairing
$$
P'_{k_0,F,X}(U)\times A_J
\to P'_{k_0,F,T}(T^*).$$
The integrable connection
$((\EE,f),(\nabla,f))$
is defined
as that corresponding to
the integrable connection
$((\EE^{(F)},f),(\nabla^{(F)},f))$
on $T^*\otimes_{k_0}F$.
We define the local system
$(V,f)$ on $T^{*an}$
of $F$-vector spaces.
We have a local system
$(V,f)^\sim
=(-1)^{-\ord f\cdot\nabla}f^\nabla
\allowbreak\rho(V)^{\otimes\ord f}$
inside a locally free
$\Cal O_U \an\otimes_{k_0}F$-module
$\EE^{an\otimes\ord f}$ on $U\an$
by fixing an extension $\EX$.
It is unramified along $D_i$ for $i\in J$
and $(V,f)$ is defined as the restriction on
$T^{*an}$
of the extension of $(V,f)^\sim$.
The comparison
$(\rho,f):(V,f)
\to
\Ker(\nabla\an\text{ on }
(\EE,f)\an)$
is the restriction of the inclusion
$(V,f)^\sim\to \EX^{an\otimes \ord f}$.
Using the tame symbol defined this way,
we obtain a pairing
$
P'_{k_0,F,X}(U)\times CH^n(X\bmod D)
\to P'_{k_0,F}(k)$
as in Section 3.

By the same argument as above,
we have
\proclaim{Theorem 1'}
Assume $X$ is projective.
Then for an object $\MM$ of
$M'_{k_0,F}(U)$,
there exists an isomorphism
$$
{\det}_F^\Gamma R\Gamma_c(U,\MM)
\otimes
{\det}_F^\Gamma R\Gamma_c(U,1)^{\otimes-\rank\MM}
\simeq
({\det}_F\MM,c_{X\bmod D})
$$
in
$P'_{k_0,F}(k)$.
\endproclaim
\proclaim{Theorem 2'}
Assume further that $\dim X=1.$
Then the isomorphism of Deligne-Riemann-Roch
induces an isomorphism
$$
{\det}_F^\Gamma R\Gamma_c(U,\MM)
\otimes
{\det}_F^\Gamma R\Gamma_c(U,1)^{\otimes-\rank\MM}
\simeq
({\det}_F\MM,c_{X\bmod D}).
$$
\endproclaim


\subheading{5. $\ell$-adic sheaves. Motives}

In this section, first we define a pairing
$$
CH^n(X\bmod D)\times\PA{k}
\to
\PAT{U}
$$
compatible with the pairings
$$
CH^n(X)\times\PA{k}
\to
\PA{X}
$$
defined in [S2] Proposition 1.
Let $k$ be a field,
$X$
be a proper regular scheme over $k$
and $U$ be an open subscheme
whose complement $D$ is
a divisor with simple normal crossings.
As in [S2], let
$\PAT{U}$
be the quotient of $\PA{U}$
classifying the abelian coverings of
$U$ tamely ramified along $D$.
It is the Pontrjagin dual of
$H^1(U,\QZ)\tame=
H^1(U,\Bbb Q/\ZZ[1/p])\oplus
H^1(X,\Bbb Q_p/\ZZ_p)$
for $p=\chara k$.

\proclaim{Proposition 3}
Let $X$ and $U$ over a field $k$ as above.
Then there is a pairing
$$
CH^n(X\bmod D)\times\PA{k}
\to
\PAT{U}
$$
characterized by the following property.
For a closed point $x\in U$,
the pairing with the class $[x]$
is given by the composition
$$
\PA{k}@>{[\kx:k]_i\cdot tr_{x/k}}>>\PA{x}@>{i_{x*}}>>\PAT{U}.
$$
of the inseparable degree times
the transfer $tr_{x/k}$ and
the map induced by the inclusion $i_x$.
\endproclaim
\demo{Proof}
By Corollary of Proposition 1,
it is unique.
By taking the dual
it is sufficient to define a pairing
$$
H^1(U,\QZ)\tame\times
CH^n(X\bmod D)
\to
H^1(k,\QZ)$$
such that, for a closed point $x\in X_0$,
the pairing with the class $[x]$
is the composition
$H^1(U,\QZ)\tame
\overset {i_x^*}\to\to
H^1(x,\QZ)
@>{[\kx:k]_i\times\Tr_{x/k}}>>
H^1(k,\QZ)$.
For the $p=\chara k$-part,
this is Proposition 1 of
[S2].
There $X$ is assumed to be smooth but
it is enough to assume only $X$ is regular.
Hence it is enough to define a pairing
$$
H^1(U,\Zm)\times
CH^n(X\bmod D)
\to
H^1(k,\Zm)$$
for an integer $m\in \ZZ$ invertible in $k$.
Let $\kappa:\Gamma(U,\Gm)\to
H^1(U,\mu_m)$
be the map induced by the Kummer sequence.

First we assume $X$
is a proper regular integral curve.
Let $\Gm(X\bmod D)$
be the complex
$[\Bbb G_{m,X}\to
\Bbb G_{m,D}]$
on $X^{et}$.
The Kummer sequence
induces a map
$H^1(X,\Gm(X\bmod D))
\to
H^2_c(U,\mu_m)$.
The pairing
$\Tr_{U/k}\circ \cup:
H^1(U,\ZZ/m)\times
H^2_c(U,\mu_m)
\to
H^1(k,\ZZ/m)$
induces a pairing
$H^1(U,\ZZ/m)\times
CH^1(X\bmod D)
\to
H^1(k,\ZZ/m)$.
By the canonical isomorphism
$CH^1(X\bmod D)
\to
H^1(X,\Gm(X\bmod D))$,
it gives the required pairing.

For general $X$, similarly as in Section 3,
we define a local pairing
$$H^1(U,\Zm)\times
A_x
\to H^1(x,\Zm)$$
for a closed point $x\in X_0$
using the tame symbol
and check the reciprocity
for $y\in X_1.$
For a while, we drop the assumption $X$ is proper
and we define the tame symbol in this context.
Let $J\subset I$ be a subset of
the index set of the irreducible components of $D$
and let $T^*\subset T=\bigcap_{i\in J}D_i$
and $A_J=\left({\prod_{i\in J}}\right)_{\ZZ}
\Gamma(X-D_i,\Gm)$
be the fiber product as in section 3.
We define a pairing
$$
H^1(U,\Zm)\times A_J
\to
H^1(T^*,\Zm).$$
Replacing $X$ by
$X-\bigcup_{i\in I-J}D_i$,
we may assume $T=T^*$ and $I=J$ and we put $A=A_J$.
It is enough to define the pairing
Zariski locally on $X$ and patch them afterward.
Hence we may assume
that the ideal sheaf
$\Cal O (-D_i)$ is a trivial
invertible $\Cal O_X $-module and
$\ord_i:A_i\to\ZZ$ is surjective
for $i\in I.$
We show that the heuristic argument
in Remark in Section 3 actually works.
For each $i\in I$,
the Kummer sequence induces a map
$A_i\to
H^1(X-D_i,\ZZ/m(1))
\to
H^2_{D_i}(X,j_{i!}\ZZ/m(1))$
where $j_i:X-D_i\to X$
is the open immersion.
Let $\tilde A$ be the tenser product
$\tilde A=\bigotimes_{i\in I}A_i$
and $r=\Card I$.
By taking the cup-product, we obtain
$\tilde A\to
H^{2r}_T(X,j_!\Zm(r))$.
The cup-product
$$H^1(U,\Zm)\times
H^{2r}_T(X,j_!\Zm(r))
\to
H^{2r+1}_T(X,\Zm(r))
\simeq
H^1(T,\Zm)$$
induces a pairing
$(\ , \ ):H^1(U,\Zm)\times
\tilde A\to
H^1(T,\Zm).$
Let $\tilde A^2$
be the subgroup of $\tilde A$
generated by
$\otimes_ia_i$ such that $\ord_ia_i=0$
for at least two $i\in I$.
Since $\ord_i:A_i\to \ZZ$
are surjective, the fiber product $A$
is isomorphic to
$\tilde A/\tilde A^2$
by $\otimes a_i\mapsto
(a_i^{\prod_{j\ne i}\ord_ja_j})_i$.
We show that the pairing $(\ ,\ )$
annihilates $\tilde A^2$.
We may assume $\Card I=2$ and
$I=\{1,2\}$.
The pairing is factorized as
$$\gather
H^1(U,\Zm)\times
A_1\times A_2
@>{(\ ,\ )_1\times res}>>
H^1(D_1-T,\Zm)\times
\Gamma(D_1-T,\Gm)\\
@>{(\ ,\ )_2}>>
H^1(T,\Zm).
\endgather$$
Further in the first pairing
$(\ ,\ )_1$,
the restriction to
$\Gamma(X,\Gm)\subset A_1$
is
$$\gather
H^1(U,\Zm)
\times
\Gamma(X,\Gm)
@>{\partial\times(\kappa\circ res)}>>
H^0(D_1,\Zm(-1))\times
H^1(D_1,\Zm(1))\\
@>{\cup}>>
H^1(D_1,\Zm).
\endgather$$
Here the first arrow is the product of
the boundary map and the restriction
followed by $\kappa:
\Gamma(D_1,\Gm)
\to
H^1(D_1,\Zm(1))$
induced by the Kummer sequence.
In the second pairing $(\ ,\ )_2$,
the restriction to
$H^1(D_1,\Zm)\times
\Gamma(D_1,\Gm)$
is trivial.
Hence the pairing annihilates
$\tilde A^2$
and the pairing
$(\ , \ ):H^1(U,\Zm)\times
 A\to
H^1(T,\Zm)$
is defined.

Using the tame symbol defined above,
we obtain a local pairing
$$(\ ,\ )_x:
H^1(U,\Zm)\times A_x\to
H^1(x,\Zm)$$
as in Section 3.
To define a pairing
$(\ ,\ )_X:
H^1(U,\Zm)\times CH^n(X\bmod D)\to
H^1(k,\Zm)$
it is sufficient to show the reciprocity law
that the pairing
$\bigoplus_x\Tr_{x/k}\circ(\ ,\ )_x:
H^1(U,\Zm)\times \bigoplus A_x\to
H^1(k,\Zm)$
annihilates the image of the boundary map
$\partial:\bigoplus_{y\in X_1}B_y
\to \bigoplus_{x\in X_0}A_x$
by Proposition 1 in Section 2.

We prove the reciprocity for higher dimension
as in the proof of Proposition 2 in section 3.
Let $D_i$ be an irreducible component of $D$,
$k_i$ be the constant field of $D_i$
and $\partial_i:
H^1(U,\Zm)\to
H^0(D_i,\Zm(-1))=
H^0(k_i,\Zm(-1))$
be the bounary map.
Then for a closed point $x\in D_i$,
the restriction of the pairing to
the $i$-th component
$\kxx\subset A_{x,i}\cap A_x$
is factorized as
$$H^1(U,\Zm)\times\kxx
@>{\partial_i\times(\kappa\circ N_{x/k_i})}>>
H^0(k_i,\Zm(-1))\times H^1(k_i,\mu_m)
@>{\cup}>>H^1(k_i,\Zm).$$
Hence the pairing annihilates
the kernel of the norm
$\bigoplus_{x\in (D_i)_0}\kxx\to k_i^\times.$
Since the composition
$\bigoplus_{y\in (D_i)_1}K_2(y)\to
\bigoplus_{x\in (D_i)_0}\kxx\to k_i^\times$
is the 0-map,
the restriction of the pairing
$(\ ,\partial\ )_X$
to the kernel of
$B_y\to \kyx$
is trivial for $y\in X_1$.
Therefore it is sufficient to show that for
each
$y\in X_1$
and $f\in \kyx$,
there is a lifting $\tilde f\in B_y$
of $f$ such that the pairing
$(\ ,\partial \tilde f)$ is trivial.

Let $Y$ be the closure of $y$
and $\tilde Y$
be the normalization of $Y$
and let $f\in \kyx$.
We take an open subset
$W\subset Y$ containing
the singularities of $Y$,
the intersection $Y\cap D_i$ for $i\notin
I_y=\{i\in I; y\notin D_i\}$
and the zeroes and poles of $f$
and take a family $\pi$
of bases $\pi_i$
of $\Cal O_Y (-D_i)|_W$
for $i\in I_y$.
Let $\tilde W^*$
be the inverse image of
$W-\bigcup_{i\in I-I_y}D_i$
in $\tilde Y$.
The tame symbol $(\ ,\pi)$
defines a map
$H^1(U,\ZZ/m)\to
H^1(\tilde W^*,\ZZ/m)$.
Similarly as loc.cit.,
it is sufficient to show the equality of the map
$$((\ ,\pi),\partial f)_{\tilde Y}=
(\ ,\partial\{\pi,f\})_X:
H^1(U,\ZZ/m)\to
H^1(k,\ZZ/m).$$
It follows from the local equality
$$((\ ,\pi),\partial_{\tilde x} f)_{\tilde x}=
(\ ,\partial_{\tilde x,y}\{\pi,f\})_{\tilde x}:
H^1(U,\ZZ/m)\to
H^1(\tilde x,\ZZ/m)$$
for $\tilde x\in \tilde Y$.

First we assume the image $x$ of $\tilde x$ is in $W$.
By replacing $X$ by an etale neighborhood
$X'$ with $x\in X'\to X$,
we may assume
$H^1(U,\Zm)$
is generated by
$\bigoplus_{i\in I_x}
H^1(X-D_i,\Zm)$.
\comment
We have
$$(\chi,\pi)=
\cases
(\chi,\pi)|_{D_{I_y-\{i\}}}
&\quad \text{ if } i\in I_y\\
\chi|_{D_{I_y}}
&\quad \text{ if } i\notin I_y
\endcases
$$
for $\chi\in
H^1(X-D_i,\Zm)$.
\endcomment
Hence we may assume
$\Card I_x\le 1$.
If $I_x=\emptyset$,
both sides
$((\ ,\pi),\partial_{\tilde x} f)$ and
$(\ ,\partial_{\tilde x,y}\{\pi, f\})$
are simply
$\ord_{\tilde x}f\cdot (\ |_{\tilde x}):
H^1(X,\Zm)\to
H^1(\tilde x,\Zm)$.
We assume $\Card I_x=1$
and let 0 be the unique element
and $\pi_0$ be a prime element of $D=D_0$ at $x$.
First assume
$y\in D=D_0$.
We may assume $\pi=(\pi_0)$.
Then they are equal to
$\ord_{\tilde x}f\cdot (\ ,\pi_0)|_{\tilde x}:
H^1(U,\Zm)\to
H^1(\tilde x,\Zm)$.
Finally assume $y\notin D$.
Then $I_y=\emptyset$, $\pi$ is the empty family
and
$(\chi,\pi)=\chi|_{\tilde W^*}$.
We have an exact sequence
$$0\to
H^1(X,\Zm)\to
H^1(U,\Zm)\to
H^0(D,\Zm(-1))$$
and the assertion is proved for
$H^1(X,\Zm)$.
Hence by localizing if necessary,
we may assume
$\Zm(1)\simeq\Zm$
on $X$ and
$\chi=\kappa(\pi_0)
\in
H^1(U,\Zm)\simeq
H^1(U,\Zm(1))$.
Then
the left hand side is
$(\chi|_{\tilde W^*},\partial_{\tilde x} f)=
\kappa\{\pi_0,f\}_{\tilde x}
\in H^1(\tilde x,\Zm)$.
We compute the right hand side.
By the explicit computation
$\partial_{\tilde x,y}f=
(-\pi_0)^{\ord_{\tilde x}f}\cdot
\{\pi_0^{-1},f\}_{\tilde x},$
we have
$$
(\chi ,\partial_{\tilde x,y}f)=
\kappa(\{\pi_0,-\pi_0\}^{\ord_{\tilde x}f}\cdot
\{\pi_0,\{\pi_0^{-1},f\}_{\tilde x}\})=
\kappa(\{\pi_0,f\}_{\tilde x}).$$
Thus the local equality is proved
if $x\in W$.

We assume $x\notin W$.
We take a local basis $\pi'_i$
of $N_i|_Y$ at $x$
for $i\in I_y=I_x$ and put
$\pi_i=u_i\cdot\pi'_i$, $u_i\in \kyx$.
Then for $\chi\in H^1(U,\Zm)$,
we have
$(\chi,\pi)=(\chi,\pi')+
\sum_i\partial_i\chi\cup\kappa(u_i)$
and $(\chi,\pi')$
is unramified at $x$.
Since $f$ is a unit at $x$,
we have
 $((\chi,\pi'),\partial_xf)=0$
and
$$\align
&((\chi,\pi),\partial_xf)=
\sum_i(\partial_i\chi\cup\kappa(u_i),\partial_xf)\\
=&
\sum_i\partial_i\chi\cup\kappa(\{u_i,f\}_x)=
-\sum_i\ord_xu_i\cdot\partial_i\chi\cup\kappa(f(x)).
\endalign$$
On the other hand, by
$\partial_{x,y}\{\pi, f\}=
(f^{-\ord_xu_i})_i
\in\bigoplus_i\kxx\subset A_x$,
the left hand side
$(\chi,\partial_{x,y}\{\pi, f\})$
is also
$-\sum_i\ord_xu\cdot\partial_i\chi\cup\kappa(f(x)).$
Thus the equality is proved for $x\notin W$.
Thus the proof of Proposition 3 is completed.

\enddemo

In the rest of paper, we assume the Gersten conjecture
for K-theory holds for a discrete valuation ring.
By [G],
it implies the Gersten conjecture
for a smooth scheme over a discrete valuation ring.

\proclaim{Proposition 4}
Let $X$ be a proper smooth
scheme over a regular noetherian connected scheme
$S$ and $D$ be a divisor of $X$
with simple normal crossings
relative to $S$.
Then there exists a unique pairing
$$
\PA{S}\times
CH^n(X\bmod D)\to
\PAT{U}$$
such that for a point $s\in S$,
the diagram
$$\matrix
\PA{s}
&\times&
CH^n(X_s\bmod D_s)
&@>>>&
\PAT{U_s}\\
@VVV@AAA@VVV\\
\PA{S}
&\times&
CH^n(X\bmod D)
&@>>>&
\PAT{U}
\endmatrix$$
is commutative.
\endproclaim

\demo{Proof}
The uniqueness follows from the
surjectivity of
$\PA{\eta}\to\PA{S}$
for the generic point $\eta$ of $S$.
We show the existence.
We show the pairing
$\PA{\eta}\times
CH^n(X\bmod D)\to
\PA{\eta}\times
CH^n(X_\eta\bmod D_\eta)\to
\PAT{U_\eta}\to
\PAT{U}$
factors the quotient $\PA{S}$.
In other word,
the image of  the cup-product
$$
H^1(U,\QZ)\tame
\times
CH^n(X\bmod D)
\to
H^1(\eta,\QZ)$$
is in
$H^1(S,\QZ)$.
By the purity of the branch locus,
we may assume $S$ is the spectrum of a
discrete valuation ring.
Let $K$ be the fraction field
and $F$ be the residue field of $S$.

For the $p$-part where $p=\chara F$,
it is Proposition 1 of [S2].
Let $m$ be an integer invertible in $F$
 and consider the pairing
$$
H^1(U,\ZZ/m)
\times
CH^n(X_\eta\bmod D_\eta)
\to
H^1(\eta,\ZZ/m).$$
By the Gersten resolution, we have an isomorphism
$$H^n_{X_F}(X,\Kn(X\bmod D))\to
\bigoplus_iH^{n-1}(D_{i,F},\KK_{n-1}(D_{i,F}))=
\bigoplus_iCH^{n-1}(D_{i,F})$$
and an exact sequence
$$
CH^n(X\bmod D)
\to
CH^n(X_K\bmod D_K)
\overset\partial\to\to
\bigoplus_iCH^{n-1}(D_{i,F}).$$
The ramification theory gives a map
$\partial_i:H^1(U,\ZZ/m)
\to
H^0(D_i,\ZZ/m(-1))$
for a irreducible component $D_i$ and
$\partial:H^1(K,\ZZ/m)
\to
H^0(F,\ZZ/m(-1))$.
Since $\Ker\partial=H^1(S,\ZZ/m)$,
the following implies a pairing
$$H^1(U,\ZZ/m)
\times
CH^n(X\bmod D)
@>>>
H^1(S,\ZZ/m)$$
is induced.
\proclaim{Claim 1}
The diagram
$$\CD
H^1(U,\ZZ/m)
\times
CH^n(X_K\bmod D_K)
@>>>
H^1(K,\ZZ/m)\\
@V{(\oplus_i\partial_i)\times\partial}VV@VV{\partial}V\\
\bigoplus_iH^0(D_{i,F},\ZZ/m(-1))\times
\bigoplus_iCH^{n-1}(D_{i,F})
@>>>
H^0(F,\ZZ/m(-1))
\endCD$$
is commutative.
Here the lower pairing is the sum of
$$\gather
H^0(D_{i,F},\ZZ/m(-1))\times
CH^{n-1}(D_{i,F})
@>{can\times\deg_F}>>
H^0(F_i,\ZZ/m(-1))\\
@>{\Tr_{F_i/F}}>>
H^0(F,\ZZ/m(-1))
\endgather$$
where $F_i$ is the constant field of $D_{i,F}$.
\endproclaim

We prove Claim 1 later and assume a pairing is defined.
We show the commutativity of the
diagram for an arbitrary point $s$ of $S$.
By considering the blowing-up at
the closure of $s$ and by the same argument
as in the proof of Proposition 1 in [S2],
we may assume that the codimension of $s$ is 1.
Hence we may assume $S$ is the
spectrum of a discrete valuation ring
and $s$ is the closed point.
For $p=\chara F$-part, it is proved in loc.cit.
Let $m$ be an integer invertible in $F$.
We take a prime element $\pi$ of $K$.
Let $sp_\pi$ denote the composition
$$\gather
CH^n(X_K\bmod D_K)
@>{\{\ ,\pi\}}>>
H^n(X_K,\KK_{n+1}(X_K\bmod D_K))\\
\overset\partial\to\to
H^{n+1}_{X_F}(X,\KK_{n+1}(X\bmod D))\\
\simeq
H^n(X_F,\Kn(X_F\bmod D_F))=
CH^n(X_F\bmod D_F)
\endgather$$
and also
$$
H^1(K,\ZZ/m)
@>{\cup\{\pi\}}>>
H^2(K,\ZZ/m(1))
\to
H^1(F,\ZZ/m).$$
The restriction to the closed fiber
$res: CH^n(X\bmod D)\to
CH^n(X_F\bmod D_F)$
is equal to the composition
$$CH^n(X\bmod D)
\to
CH^n(X_K\bmod D_K)
@>{sp_\pi}>>CH^n(X_F\bmod D_F)$$
and the restriction
$
H^1(S,\ZZ/m)\to
H^1(F,\ZZ/m)$
is
$$
H^1(S,\ZZ/m)\to
H^1(K,\ZZ/m)
@>{sp_\pi}>>
H^1(F,\ZZ/m).$$
Hence it is sufficient to show the following.

\proclaim{Claim 2}
The diagram
$$\CD
H^1(U,\ZZ/m)\times
CH^n(X_K\bmod D_K)
@>{(res,sp_\pi)}>>
H^1(U_F,\ZZ/m)\times
CH^n(X_F\bmod D_F)\\
@VVV @VVV\\
H^1(K,\ZZ/m)
@>>{sp_\pi}>
H^1(F,\ZZ/m)
\endCD$$
is commutative.
\endproclaim

We deduce Claims 1 and 2 from
\proclaim{Claim 3}
Let $y$ be a closed point of $U_K$
and $Y$ be the closure of $\{y\}$ in $X$.
Then

(1). The image of the class $[y]$ by
$\partial:
CH^n(X_K\bmod D_K)\to
\bigoplus_iCH^{n-1}(D_{i,F})$
is
$\partial[y]=((D_i,Y))_i\in
\bigoplus_iCH^{n-1}(D_{i,F})$.
Here $(D_i,Y)$ is the intersection product
$\sum_{x\in Y\cap D_i}
\text{length }\Cal O _{Y\cap D_i,x}\cdot [x]$.

(2). For a closed point $x$ of $Y_0=Y\cap X_F$,
let $\partial_{x,y}$ denote the map
$\kyx\to A_x$ defined in Section 2.
Then the diagram
$$\CD
\kyx
@>{\oplus\partial_{x,y}}>>
\bigoplus_{x\in Y_0}A_x\\
@VVV@VVV\\
H^n(X_K,\KK_{n+1}(X_K\bmod D_K))
@>>{\partial}>
CH^n(X_F\bmod D_F)
\endCD$$
is commutative.
\endproclaim

We show Claim 3 implies Claims 1 and 2 and
hence Proposition 4.

\demo{Proof of Claim 1}
Let $\chi\in
H^1(U,\ZZ/m)$
and $y$ be a closed point of $U_K$.
For a closed point $\tilde x$
of the normalization
$\tilde Y$
of $Y$,
the boundary
$\partial_{\tilde x}(\chi|_y)
\in H^0(\tilde x,\ZZ/m(-1))$
is given by
$\sum_i\ord_{\tilde x}(D_i|_{\tilde Y})
\cdot (\partial_i\chi)|_{\tilde x}$.
Hence $$\partial(\chi,[y])=
\partial(\Tr_{y/K}(\chi|_y))=
\sum_{x\in\tilde Y_0}
\Tr_{\tilde x/F}
\partial_{\tilde x}(\chi|_y)$$
is equal to
$\sum_i\sum_{\tilde x}
[\kappa(\tilde x):F_i]\cdot
(\ord_{\tilde x}(D_i|_{\tilde Y}))
\cdot
\Tr_{F_i/F}(\partial_i\chi).$
Since
$\sum_{\tilde x}
[\kappa(\tilde x):F_i]\cdot
\ord_{\tilde x}(D_i|_{\tilde Y})=
\deg_{F_i}(D_i,Y)$,
we have
$\partial(\chi,[y])=
\sum_i(\partial_i\chi,(D_i,Y))$.
Since the right hand side is
$(\sum_i\partial_i\chi,
\partial[y])$
by Claim 3 (1)
and $CH^n(X_K\bmod D_K)$ is generated by
the class of closed points of $U_K$,
Claim 1 is proved.
\enddemo

\demo{Proof of Claim 2}
We use the notation above.
By Claim 3 (2),
we have $(\chi,sp_\pi[y])
=\sum_{x\in Y_0}\Tr_{x/F}
(\chi,\partial_{x,y}(\pi))_x$.
On the other hand,
we have $sp_\pi(\chi,[y])
=\partial_F\Tr_{y/K}(\chi,\pi)
\allowbreak
=\sum_{\tilde x\in\tilde Y_0}
\Tr_{\tilde x/F}(\chi,\pi)_{\tilde x}$.
Hence it is sufficient to show
$(\chi,\partial_{\tilde x,y}(\pi))_{\tilde x}
=(\chi,\pi)_{\tilde x}$.
This is proved by exactly
the same argument as the equality
$(\chi,\partial_{\tilde x,y}(\{\pi,f\}))_{\tilde x}
=((\chi,\pi),\partial_{\tilde x}f)$
in the proof of the reciprocity law in
Proposition 4.
\enddemo

\demo{Proof of Claim 3}
Similarly as in the proof of Proposition 1,
it is sufficient to prove that
the boundary map
$$\gather
\ZZ=H^n_y(X,\Kn(X\bmod D))
\to
H^{n+1}_x(X,\Kn(X\bmod D))
=\bigoplus_{D_i\ni x}\ZZ\\
\kyx=H^n_y(X,\KK_{n+1}(X\bmod D))
\to
H^{n+1}_x(X,\KK_{n+1}(X\bmod D))
=A_x
\endgather$$
are given by
$1\mapsto
((D_i,Y)_x)_i$
and $\partial_{x,y}$
respectively.
The second one is proved by the same argument
as Proposition 1 in Section 2 and the first one
is similar and easier.
The only point to be proved is the crucial Claim
loc.cit.. However,
it is proved by taking etale locally a
smooth projection $X\to D$ over $S$
by the same argument.
The rest is proved in the same way and we
leave the detail to the reader.
\enddemo
\enddemo


Using Propositions 3 and 4 above,
we have a variant of Theorem 1 [S2]
for the relative canonical class in
the relative Chow group.
To state it, we review Jacobi sums
following [A].
Let $k$ be a field of
characteristic $p\ge0$
and $F$ be a field of characteristic 0.
Let $\bar k$ and $\bar F$ be separable closures
and $\hZZ'(1)_{\bar k}=
\projlim_{p\nmid d}\mu_d(\bar k)$
and $\QQ/\ZZ(1)_{\bar F}=
\injlim_{d}\mu_d(\bar F)$.
Let $B_F(k)$
be the free abelian group with basis
$Hom(\hZZ'(1)_{\bar k},\QZ(1)_{\bar F})$
and $B_F^0(k)$ be the kernel of the natural map
$B_F^0(k)\to
Hom(\hZZ'(1)_{\bar k},\QZ(1)_{\bar F})$.
If we choose an isomorphism
$\hZZ'(1)_{\bar k}\simeq
\hZZ'(1)_{\bar F}$,
the group
$Hom(\hZZ'(1)_{\bar k},\QZ(1)_{\bar F})$
is naturally identified with
$(\QZ)'=\injlim_{p\nmid d}\ZZ/d$
and there are inclusions
$B_F(k)\to\Bbb B$ and
$B_F^0(k)\to\Bbb B^0$.
Here the free abelian group
$\Bbb B$ with basis $\QZ$ and
$\Bbb B^0=\Ker\Bbb B\to\QZ$
are as in [A].
The groups
$B_F(k)$ and $B_F^0(k)$
have natural actions of
the absolute Galois groups
$G_k=\Gal(\bar k/k)$
and
$G_F=\Gal(\bar F/F)$.
The automorphism group of
$Hom(\hZZ'(1)_{\bar k},\QZ(1)_{\bar F})$
as an abstract group is
$\hZZ^{\prime\times}=\projlim_{p\nmid d}(\ZZ/d)^\times$.
The action is given by the cyclotomic characters
$G_k\to\hZZ^{\prime\times}$ and
$G_F\to\hZZ^{\prime\times}$.
We define Jacobi sum $J(\aaa)$ for
an element
$\aaa\in B^0_F(k)^{G_k\times G_F}$
fixed by the actions of
$G_k$ and $G_F$.

First we consider the case where
$k$ is a finite field of order $q$.
The subgroup
$B_F(k)^{G_k}$
of $B_F(k)$ consisting of the elements
fixed by $G_k$ is generated
by the elements of the form
$\aaa=\sum_{i=0}^{f-1}
[q^ia]$.
 Here the order $m(a)$ of $a\in$\linebreak
$Hom(\hZZ'(1)_{\bar k},\QZ(1)_{\bar F})$
is prime to $q$ and $f$ is the order of $q$
in $(\ZZ/m(a))^\times$.
For a non-trivial additive character
$\psi_0:k\to \bar F^\times$,
we define a map
$g(\ ,\psi_0):
B_F(k)^{G_k}\to
\bar F^\times$
by
$$g(\aaa ,\psi_0)=
-\sum_{x\in E_f^\times}a^{-1}(x)
\psi_0(\Tr_{E_f/k}x)$$
for an element
$\aaa=\sum_{i=0}^{f-1}
[q^ia]$ as above.
Here $E_f$ is the extension of $k$ of degree $f$
in $\bar k$ and $a$ is regarded as a character
of $E_f^\times$
by $E_f^\times=\mu_{q^f-1}(\bar k)$.
The restriction $J_k$
of
$g(\ ,\psi_0)$ to
$B^0_F(k)^{G_k\times G_F}$
is independent of the choice of $\psi_0$
and the image is in $F^\times$.
Hence $J_k$ is a map
$B^0_F(k)^{G_k\times G_F}
\to F^\times$.

We consider a general field $k$.
Let
$\aaa=\sum_an_a[a]\in B^0_F(k)$
be an element fixed by $G_k$ and $G_F$ and
let $m(\aaa)$ denote the order of the subgroup of
$Hom(\hZZ'(1)_{\bar k},\QZ(1)_{\bar F})$
generated by $\{a;n_a\ne0\}$.
We define a subfield $F_{\aaa}$ of $F$
which is a subfield of $\QQ(\zeta_{m(\aaa)})$.
Let $G_{\aaa}$ be the subgroup of
$(\ZZ/m(\aaa))^\times$
fixing $\aaa$.
We define $F_{\aaa}$ to be the subfield of
$\QQ(\zeta_{m(\aaa)})\subset\bar F$
fixed by the subgroup $G_{\aaa}$ of
$(\ZZ/m(\aaa))^\times=
\Gal(\QQ(\zeta_{m(\aaa)})/\QQ)$.
Since $G_{\aaa}$ includes the image of the
cyclotomic character $G_F\to
(\ZZ/m(\aaa))^\times$, it is a subfield of $F$.
Take an isomorphism
$\hZZ'(1)_{\bar k}\simeq
\hZZ'(1)_{\bar F}$.
Then
the image of $\aaa\in B_F^0(k)$
in $\Bbb B^0$
is fixed by $G_{F_{\aaa}}$
and defines an algebraic Hecke character
$J(\aaa)$ of $F_{\aaa}$ with values in $F_{\aaa}$
with conductor dividing a power of $m(\aaa)$
by [A].
Since $G_{\aaa}$ also includes the image of
$G_k$, the isomorphism
$\hZZ'(1)_{\bar k}\simeq
\hZZ'(1)_{\bar F}$ chosen
induces a homomorphism
$\Cal O _{F_{\aaa}}[1/m(\aaa)]\to k$
and a homomophism
$G_k\ab\to
\PA{\Cal O _{F_{\aaa}}[1/m(\aaa)]}$.
If $F$ is a finite extension of
$\Ql$ for a prime $\ell\ne p$,
the algebraic Hecke character $J(\aaa)$
induces a character
$\PA{\Cal O _{F_{\aaa}}[1/\ell m(\aaa)]}
\to F^\times$.
The composite character
$J(\aaa):G_k\ab
\to F^\times$
is independent of the choice of an isomorphism
$\hZZ'(1)_{\bar k}\simeq
\hZZ'(1)_{\bar F}$.
If $k$ is a finite field,
the value of the geometric Frobenius $Fr_k$
is given by
$J(\aaa)(Fr_k)=J_k(\aaa)\in F^\times$.

We return to a geometric situation.
Let $X$ be a proper smooth scheme
over a field $k$ of characteristic $p\ge0$
and $U$ be an open subscheme of $X$
whose complement $D$ is a divisor with simple
normal crossings.
Let $F_\llll$ be a finite extension
of $\Ql$, for $\ell\ne p$,
and $V_\llll$ be a smooth $F_\llll$-
sheaf on $U\et$ of rank $r$.
We assume that the ramification of $V_\llll$
along $D$ is tame and
the local monodromy at each irreducible
component $D_i$ is quasi-unipotent.
The assumption is automatically satisfied
if $k$ is a number field
by the monodromy theorem of Grothendieck.
We define an element
$\aaa_{D,V_\llll}
\in B_{F_\llll}^0(k)$
fixed by $G_k$ and $G_F$
and a character
$J_{D,V_\llll}=
J(\aaa_{D,V_\llll}):
G_k\ab
\to F_\llll^\times$.
Let $k_i$
be the constant field of an irreducible component
$D_i$ and $c_i$
be the Euler number of $D_i^*\otimes_{k_i}\bar k_i$
where $D_i^*=D_i-\cup_{j\ne i}D_j$ as before.
Since the ramification of $V_\llll$
is tame, the local monodromy defines a
continuous representation
$\rho_i:
\hZZ'(1)_{\bar k_i}\to GL_r(F_\llll)$.
Since it is assumed to be quasi-unipotent,
its semi-simplification
$\rho_i^{ss}$
is the direct sum $\bigoplus_{j=1}^ra_{ij}$ of characters
$a_{ij}:
\hZZ'(1)_{\bar k_i}\to \QZ(1)_{\bar F_\llll}$.
We put $\aaa_i=\sum_{j=1}^r[a_{ij}]
\in B_{F_\llll}(k_i)$.
It is fixed by $G_{k_i}$
and $G_{F_\llll}$.
We define
$\aaa_{D,V_\llll}
\in B_{F_\llll}(k)$
to be
$\sum_ic_i\Tr_{k_i/k}(\aaa_i)$.
It is
fixed by $G_k$ and $G_{F_{\lambda}}$.
By the first assertion of Theorem 1 [S2],
it is in
$B^0_{F_\llll}(k)$.

Now we state a variant of
Theorem 1 [S2] for relative Chow group.

\proclaim{Theorem S}
Let $U$ be a smooth scheme over a field $k$.
Let $F_\llll$ be a finite extension of $\Ql$
where $\ell\ne\chara k$
and $V_\llll$ be a smooth $F_\llll$-sheaf of rank r on $U$.
Assume the following conditions (1)-(3) are satisfied.
\roster
\item
There exists a projective smooth scheme
$X$ over $k$ including $U$ as the complement of
a divisor $D$ with relative normal crossings.
\item
The ramification of $V_\llll$ along $D$ is tame.
\item
There exists a scheme $S_0$ of finite type over $\ZZ$
and a $S_0$-scheme $U_0$ such that
$U=U_0\times_{S_0}\Spec k$ and $V_\llll$
is a pull-back of a smooth $F_\llll$-sheaf on $U_0$.
\endroster
Then the local monodromies are quasi-unipotent
and the character
$J_{D,V_\llll}:
G_k\ab
\to F_\llll^\times$
is defined.
Assume also the Gersten conjecture for
discrete valuation rings.
Let $(\det V_\llll,c_{X\bmod D})$
be an $F_{\llll}$-representation of dimension 1
of $G_k$ defined by pulling back the
$F_\llll$-representation $\det V_\llll$ of
$\PAT{U}$ by the pairing $(\ ,c_{X\bmod D})
:G_k\to \PAT{U}$.
Then
there is an isomorphism of $F_\llll$-representations
of dimension 1 of $G_k$
$$
\det R\Gamma_c(U_{\bar k},V_\llll)\otimes
(\det R\Gamma_c(U_{\bar k},F_\llll))^{\otimes-r}
\simeq
(\det V_\llll,c_{X\bmod D})\otimes
J_{D,V_\llll}^{\otimes-1}.$$
\endproclaim

\demo{Proof}
By the assumption (3),
the first assertion is an immediate consequence
of the monodromy theorem of Grothendieck.
The isomorphism is proved by the same way as Theorem 1 loc.cit.
By a specialization argument as in the proof of
Theorem 1 [S2] using Proposition 3 and 4 this time,
it is reduced to the case where $k$ is finite.
For a finite field $k$,
we prove the same formula
with $c_{X\bmod D}\in CH^n(X\bmod D)$
as in Theorem 1 [S1]
by the same argument
using the product formula by Laumon
and Lemmas in Section 2.
Then it is enough to show
$J_{D,V_\llll}(Fr_k)=\tau_{D/F}(\rho,\psi_0)$
where the right hand side is as in loc.cit..
It is easily checked by comparing the definition
of the both side in terms of Gauss sums.
We leave the detail to the reader.
\enddemo

Finally we consider an arithmetic situation.
Let $k$ and $F$ be finite extension
of the rational number field $\QQ$ and
$U$ be a smooth $k$-scheme.
We write $M'_{\QQ,F}(U)=M_F(U)$ etc. for short.
\proclaim{Lemma 5.2}
Let $X$ be a proper smooth scheme over $k$
including $U$ as the complement of a divisor
$D$ with simple normal crossings.
Let $\MM=((\EE,\nabla),V,\rho)$ be an object of $M_F(U)$
and assume that the action of the local monodromy
$T_i$ on $V$ at each component $D_i$ of
$D$ is quasi-unipotent. Then
there exists a unique extension
$(\EX,\nabla)$ of $(\EE,\nabla)$
such that the roots $\tau$
of the characteristic polynomials
$\Phi_{\EX,i}(T)$
are rational numbers satisfying
$0<\tau\le1$.
Further the $\Cal O_X $-module $\EX$
is locally free and,
if $(\EE,\nabla)$ has a multiplication by $F$,
it is extended to $\EX$.
\endproclaim
\demo{Proof}
Since the monodromy is quasi-unipotent,
and by $T_i=\exp(-2\pi\sqrt{-1}\res_i\nabla)$
(Chap. II Theorem 1.7 [D1]),
the roots are rational numbers.
By Proposition 5.4 loc.cit.,
There exists a unique extension
$(\EX,\nabla)$
satisfying the condition
after extending the scaler to $\CC$
and it is locally free.
Since the condition is invariant
by any automorphism of $\CC$ over $k$,
it is defined over $k$.
By the uniqueness, the multiplication is extended.
\enddemo
We call the extension $\EX$
defined in Lemma the canonical extension of $\EE$.
Let $(\EE,\nabla)$ be as in Lemma 5.2 with
an multiplication by $F$ and
$\EX$ be the canonical extension.
We consider a finite decreasing filtration $Fil$ on
$\EX^{(F)}$
by coherent
$\Cal O _{X\otimes_{\QQ}F}$-submodules.
We assume that the Griffiths transversality
$\nabla(Fil^p\EX^{(F)})
\subset
Fil^{p-1}\EX^{(F)}
\otimes_{\Cal O_X }\Omega^1_X(log D)$
is satisfied.
It induces a filtration
also denoted by $Fil$
on the de Rham complex
$DR(\EX^{(F)})
=(\EX^{(F)}\otimes \Omega^\bullet_X(\log D))$
by
$(Fil^pDR(\EX^{(F)}))^q
=(Fil^{p-q}\EX^{(F)}\otimes \Omega^q_X(\log D))$.

Let $\lambda$ be a finite place of $F$.
We consider a triple $\tMM=(V_\lambda,\MM,Fil)$
of smooth $F_\lambda$-sheaf $V_\llll$
on $U\et$, an object $\MM$
of $M_F(U)$ and a filtration $Fil$
on the canonical extension $\EX^{(F)}$
as above.
We assume the following compatibility.
\roster
\item
There exists an $\Cal O_k $-scheme
$U_{\Cal O_k }$ of finite type
and a smooth $F_\llll$-sheaf
$V_\ell$ on $U_{\Cal O_k }\et$
such that $U=U_{\Cal O_k }\otimes k$
and $V_\llll$ on $U$ is the pull-back of
$V_\llll$ on $U_{\Cal O_k }$.
\item
The local system $V\otimes F_\llll$
on $U\an$ is isomorphic to the pull-back of
$V_\llll$ by $U\an\to U\et$.
\endroster
Then by the assumption (1)
and by the monodromy theorem of
Grothendieck,
the local monodromy
of $V_\llll$
along each component of $D$
is quasi-unipotent.
By the assumption (2),
it is also quasi-unipotent for $V$.
Hence by Lemma 5.2,
there is a canonical extension $\EX^{(F)}$.
If $X=U=\Spec k$ and rank $\tMM=1$,
a filtration $Fil$ is determined
by a function $n:\Spec k\otimes F\to \ZZ$
by the condition
$Gr^p\EX^{(F)}(x)\ne 0$
if and only if $p=n(x)$ for $x\in
\Spec k\otimes F$.
By this correspondence,
we identify a filtration with an element in
$\ZZ^{\Spec k\otimes F}$.
The isomorphism class group of
the triples
of rank 1 on $\Spec k$
is
$$\injlim_SHom(G_{k,S}\ab,F_\llll^\times)\times
((k\otimes F)^\times\backslash
(k\otimes F\otimes \CC)^\times/
(F^\times)^{Hom(k,\CC)})\times
\ZZ^{\Spec k\otimes F}.$$
Here $S$ runs finite sets of places of $k$
and $G_{k,S}$ is the quotient of $G_k$
classifying the extensions unramified outside $S$.

For $\aaa=\sum_an_a[a]
\in B^0_F(k)$
fixed by $G_k$ an $G_F$,
we define a triple
$$\tilde J(\aaa)=
(J(\aaa),[\Gamma(c\aaa)],\langle\aaa\rangle)$$
of rank 1.
The character $J(\aaa):G_k\ab\to
F_\llll^\times$ is defined by
an algebraic Hecke character as above.
We define
$\Gamma(c\aaa)\in
(k\otimes F\otimes \CC)^\times=
(\CC^\times)^{Hom(k,\CC)\times Hom(F,\CC)}$
and
$\langle\aaa\rangle\in
\ZZ^{\Spec k\otimes F}$
as follows.
Let the fractional part
 $\langle\ \rangle:
\QZ\to\QQ\cap[0,1)$ be a section of
$\QQ\to\QZ$.
We define $\Gamma(c?):
\Bbb B=\bigoplus_{\QZ}\ZZ\to \CC^\times$
be the multiplicative map
defined by
$a\in\QZ\mapsto
\Gamma(1-\langle a\rangle)$ and
 $\langle\ \rangle:
\Bbb B\to \QQ$
be the linear map induced by
 $\langle\ \rangle$.
For $\sigma:k\to\CC$
and $\tau:F\to\CC$,
we take
isomorphisms
$\hZZ(1)_{\bar k}\simeq
\hZZ(1)_{\CC}$  and
$\hZZ(1)_{\bar F}\simeq
\hZZ(1)_{\CC}$
extending $\sigma$ and $\tau$
and let
$(\sigma,\tau):
Hom(\hZZ(1)_{\bar k},
\QZ(1)_{\bar F})\to
Hom(\hZZ(1)_{\CC},
\QZ(1)_{\CC})=\QZ$
be the induced isomorphism.
Since $\aaa$ is fixed by $G_k$ and $G_F$,
the image $(\sigma,\tau)(\aaa)\in\Bbb B$
is independent of the choice of isomorphisms.
We define
$\Gamma(c\aaa)=
\Gamma(c(\sigma,\tau)(\aaa))_{\sigma,\tau}
\in
(\CC^\times)^{Hom(k,\CC)\times Hom(F,\CC)}$.
Let $k\otimes F=
\prod_j K_j$
be the decomposition into
product of fields $K_j$.
For $K_j$, we
take isomorphisms
$\hZZ(1)_{\bar k}\simeq
\hZZ(1)_{\bar K_j}$  and
$\hZZ(1)_{\bar F}\simeq
\hZZ(1)_{\bar K_j}$
extending the inclusions and
let
$\rho_j:
Hom(\hZZ(1)_{\bar k},
\QZ(1)_{\bar F})\to
Hom(\hZZ(1)_{\bar K_j},
\QZ(1)_{\bar K_j})=\QZ$
be the induced isomorphism.
The image $\rho_j(\aaa)\in\Bbb B$
is independent of the choice of isomorphisms.
We define
$\langle\aaa\rangle=
(\langle\rho_j(\aaa)\rangle)_j\in
\ZZ^{\Spec k\otimes F}$ for $\aaa \in \Bbb B^0$.

The final result in this paper is a variant of
Theorem 1 in Section 4 and of Theorem S
for a triple $\tMM=(V_\llll,\MM,Fil)$.
We define the both sides.
The determinant of the cohomology
$\det_FR\Gamma_c(U,\tMM)$
is the triple
$(\det_{F_\llll}
R\Gamma_c(U_{\bar k},V_\llll),
\det_F
R\Gamma_c(U,\MM),
\allowmathbreak
Fil)$.
The first two components are already defined and we
define the filtration.
The filtration $Fil$ on $\EX^{(F)}$
induces a filtration $Fil$
on the de Rham complex $DR(\EX^{(F)})$.
The filtration $Fil$
on $\det_F R\Gamma_c(U,DR(\EE))$
is induced by this filtration
by using the isomorphism
$R\Gamma_c(U,DR(\EE))
\simeq R\Gamma(X^{(F)},DR(\EX^{(F)}))$
of perfect complexes of $k\otimes F$-modules.
We define a pairing $(\det_F\tMM,c)$
with an element $c\in CH^n(X\bmod D)$.
It is the triple
$((\det_{F_\llll}
V_\llll,c),
(\det_F
\MM,c),
Fil)$.
The first two components are already defined and we
define a function
$n:\Spec k\otimes F\to\ZZ$
defining a filtration.
Let $f$ denote the map $X^{(F)}\to \Spec k\otimes F$.
For a connected component $X_j$ of $X^{(F)}$,
let $n_j$
be the unique integer such that
$Gr^{n_j}(\det\EX^{(F)})\ne 0$
at the generic point of $X_j$
and let $\deg_j$ denote the map $CH^n(X_j)\to CH^0(f(X_j))=\ZZ$.
We define
$n:\Spec k\otimes F\to\ZZ$
by putting $n(y)=
\sum_{X_j\subset f^*(y)}n_j\deg_j(c|_{X_j})$.
Finally we define the Jacobi sum $J_{D,\tMM}$.
Let $\aaa_{D,V_\llll}\in B_{F_\llll}^0(k)$
be the element fixed by $G_k$ defined for the $F_\llll$-sheaf
$V_\llll$ before.
By the assumption (2),
it comes from an element
$\aaa_{D,\tMM}\in B_{F}^0(k)$
fixed also by $G_F$.
We define the triple
$J_{D,\tMM}$ to be the triple
$\tilde J(\aaa_{D,\tMM})$
for
$\aaa_{D,\tMM}\in B_{F}^0(k)$
above.

We define the notion that a triple
$\tMM=
(V_\llll,\MM,Fil)$
is determinantally motivic.
First we consider the case where
$U=\Spec k$ and $\rank \tMM=1$.
A triple $\tMM$ is called motivic
if it comes from an algebraic Hecke
character $\chi$ of $k$ with value in $F$.
More precisely,
\roster
\item
For a finite place $p$
of $k$ prime to the conductor of $\chi$
and to the characteristic of $\llll$,
the $\llll$-adic representation $V_\llll$
of the absolute Galois group $G_k$
of $k$ is unramified and
the action of the geometric Frobenius
$Fr_p$ on $V_\llll$ is the multiplication
by $\chi(p)$.
\item
The class
$[\MM]\in MPic_F(k)=
(k\otimes F)^\times\backslash
(k\otimes F\otimes \CC)^\times/
(F^\times)^{Hom(k,\CC)}$
is equal to the period
$p'(\chi)$ (8.7 [D4])
of a motive $M(\chi)$ for $\chi$
in the category of motives generated by
potentially CM abelian varieties.
\item
Decompose $k\otimes F\simeq
\prod_i K_i$
into products of fields
and let
$\chi_{alg}=
\prod_iN_{K_i/F}^{n_i}\circ can$.
Then $Gr^p(\EX^{(F)}\otimes_{k\otimes F}K_i)\ne 0$
if and only if $p=n_i$.
\endroster
In general, a triple $\tMM$ is called
determinantally motivic
if for any closed point $x\in U$,
the fiber at $x$ of the determinant
$\det_F\tMM=
(\det_{F_\llll}V_\llll,
\det_F\MM, Fil)$
is motivic in the sense just defined.

\proclaim{Theorem 3}
Let $k$ and $F$ be number fields
and $U$ be the complement of
a normal crossing divisor $D$
in a projective smooth scheme $X$ over $k$.
Let $\tMM=
(V_\llll,\MM,Fil)$
be a triple of a smooth
$F_\llll$-sheaf $V_\llll$ on $U\et$,
an object $\MM$ of $M_F(U)$ of rank $r$
satisfying conditions (1) and (2)
and a filtration $Fil$
on the canonical extension $\EX^{(F)}$.
Assume that the spectral sequence
$E_1^{p,q}=
H^{p+q}(X^{(F)},Gr^pDR(\EX^{(F)}))
\Rightarrow
H^{p+q}(X^{(F)},DR(\EX^{(F)}))
$
degenerates at $E_1$-terms
and that the Gersten conjecture holds
for discrete valuation rings
which are local rings
of a smooth model of $X$ over some
open subscheme of $\Cal O_k $.
Then

(1). There is an isomorphism
$$
{\det}_F R\Gamma_c(U,\tMM)\otimes
({\det}_F R\Gamma_c(U,1))^{\otimes-r}
\simeq
({\det}_F \tMM,c_{X\bmod D})\otimes
J_{D,\tMM}^{\otimes-1}.$$
(2). Further if $\tMM$ is determinantally motivic,
then
${\det}_F R\Gamma_c(U,\tMM)$
is also motivic.
\endproclaim
\demo{Proof}
(1).
The $F_\llll$-sheaves in the first component
are isomorphic to each other
by Theorem $S$.
By Theorem 1' in Section 4,
to prove the second components, the objects of $P_F(k)$,
are isomorphic,
it is sufficient to show that
$$\Gamma(\nabla:\EX)=
\Gamma(c\aaa_{D,\tMM})$$
$\in (k\otimes F\otimes\CC)^\times$.
Let $n_c(U,\tMM)$, $n(\det\tMM,c_{X\bmod D})$
etc. be the element
$\ZZ^{\Spec k\otimes F}$
corresponding to the filtration
on ${\det}_F R\Gamma(X,DR(\EX))$,
$(\det \EE,c_{X\bmod D})$ etc..
To prove the filtration are the same,
it is enough to show
$$
n_c(U,\tMM)
-rn_c(U,1)=
n(\det\tMM,c_{X\bmod D})-
\langle\aaa_{D,\tMM}\rangle$$
in $\ZZ^{\Spec k\otimes F}$

Let $\aaa_i=\sum_an_{i,a}[a]
\in B_{k_i}(F)$
be the element appeared in
the definition of
$\aaa_{D,V_\llll}$.
By $T_i=\exp(-2\pi\sqrt{-1}\cdot\res_i\nabla)$
and since the eigenvalues of
$\res_i\nabla$ is in $(0,1]$,
the image of
$\Phi_{\EX,i}(T)\in
k_i(T)^\times$
in
$(k_i\otimes F\otimes\CC(T))^\times
=(\CC(T)^\times)^{Hom(k_i,\CC)\times Hom(F,\CC)}$
is equal to
$(\prod_a(T-(1-\langle(\sigma,\tau)(a)\rangle))^{n_{i,a}})_{(\sigma,\tau)}$.
By $\Phi_{\EX}(T)
=\prod_iN_{k_i/k}
\Phi_{\EX,i}(T)^{c_i}$
and
$\aaa_{D,\tMM}=
\sum_ic_i\Tr_{k_i/k}\aaa_i$,
the image of
$\Phi_{\EX}(T)\in
k(T)^\times$
in
$(k\otimes F\otimes\CC(T))^\times
=(\CC(T)^\times)^{Hom(k,\CC)\times Hom(F,\CC)}$
is equal to
$(\prod_a(T-(1-\langle(\sigma,\tau)(a)\rangle))^{n_{a}})_{(\sigma,\tau)}$
for $\aaa_{D,\tMM}=\sum_an_a[a]$.
{}From this the equality
$\Gamma(\nabla:\EX)=
\Gamma(c\aaa_{D,\tMM})$
follows immediately.
Also we have
$\Tr\ \res_i\nabla=
(r-\langle(\sigma,\tau)(\aaa_i)\rangle)_{(\sigma,\tau)}$
in $k_i^\times\subset
(k_i\otimes F\otimes\CC)^\times$.
Hence we have
$-\langle\aaa_{D,\tMM}\rangle=
\sum_ic_i\cdot
\Tr_{k_i/k}(\Tr \ \res_i\nabla-r)$
in $\ZZ^{\Spec k\otimes F}
\subset k\otimes F$.
Therefore the second equality  becomes
$$
n_c(U,\tMM)
-rn_c(U,1)=
n(\det\tMM,c_{X\bmod D})+
\sum_ic_i\cdot
\Tr_{k_i/k}(\Tr \ \res_i\nabla-r).$$
Proof of this equality will complete the
proof of (1).

We compute
$n_c(U,\tMM)$.
To simplify the notation,
we write $\dim$
for the rank of
$k\otimes F$-module,
which is
a function on
$\Spec k\otimes F$
and write $\chi$ for the
alternating sum of dim.
By the assumption that the spectral sequence
degenerates at $E_1$
and by the definition of the filtration on
$\det_FR\Gamma(X^{(F)},DR(\EX^{(F)}))$,
we have
$n_c(U,\tMM)
=\sum_pp\cdot\chi
(X^{(F)},Gr^pDR(\EX^{(F)}))$.
Since
$(Gr^pDR(\EX^{(F)}))^q=
Gr^{p-q}\EX^{(F)}\otimes\Omega^q_X(\log D)$,
by putting $r=p-q$, we have
$$\gather
n_c(U,\tMM)
=\sum_{q,r}(-1)^q(q+r)\chi
(X^{(F)},
Gr^r\EX^{(F)}\otimes\Omega^q_X(\log D))\\
=\sum_r\left(r\sum_q(-1)^q\chi
(X^{(F)},
Gr^r\EX^{(F)}\otimes\Omega^q_X(\log D))\right)\\
+\sum_q(-1)^qq\chi
(X^{(F)},
\EX^{(F)}\otimes\Omega^q_X(\log D))
\endgather$$
We write
$n^1_c(U,\tMM)$ and
$n^2_c(U,\tMM)$ for the first term
and the second term respectively.
We prove equaities
$$\gather
n^1_c(U,\tMM)
=n({\det}_F\tMM,c_{X\bmod D}) \tag I\\
n^2_c(U,\tMM)
-rn^2_c(U,1)=
\sum_ic_i\cdot
\Tr_{k_i/k}(\res_i\nabla-r).\tag II
\endgather$$
Since
$n^1_c(U,1)=0$,
the equalities prove the assertion (1) of Theorem 3.
We prove the equality I.
For a connected component
$X_j$ of $X^{(F)}$,
let $n_j^r$
be the rank of $Gr^r$
at the generic point and
put $c_j=\deg_j(c_{X\bmod D}|_{X_j})$
where
$\deg_j:CH^n(X_j)\to CH^n(y)=\ZZ$
for $y=f(X_j)\in\Spec k\otimes F$.
By Riemann-Roch,
we have
$$\sum_q(-1)^q\chi_{\ky}
(X_j,
Gr^r\EX^{(F)}\otimes\Omega^q_X(\log D)|_{X_j})
=n_j^r\cdot c_j.$$
Since the integer $n_j$
appeared in the definition of
$n({\det}_F\tMM,c_{X\bmod D})$
is $n_j=\sum_rrn_j^r$,
the function
$n^1_c(U,\tMM):
\Spec k\otimes F\to\ZZ$
is $y\mapsto \sum_{X_j\subset f^*(y)}n_jc_j$.
It is
$n({\det}_F\tMM,c_{X\bmod D})$
by definition.

To prove II, we show
\proclaim{Lemma 5.3}
Let $X$ be a proper smooth scheme
over a field $k$ of dimension $n$
and $\FF$ be a locally free
$\Cal O_X $-module of rank $n$.
Let $\EE_1$ and $\EE_2$
be coherent $\Cal O_X $-modules
whose ranks are the same
on a dense open subscheme of $X$.
Then
$$\gather
\sum_q(-1)^q q
\chi(X,\EE_1\otimes\wedge^q\FF)
-\sum_q(-1)^q q
\chi(X,\EE_2\otimes\wedge^q\FF)\\
=(-1)^n(c_1(\EE_1)-c_1(\EE_2),
c_{n-1}(\FF)).
\endgather$$
\endproclaim
\demo{Proof}
By Riemann-Roch,
the left hand side is the
degree 0 part of
$$(ch(\EE_1)-
ch(\EE_2))\cdot
(\sum_q(-1)^qq\cdot
ch(\wedge^q\FF))\cdot
td(\Omega^1_X).$$
Since
$ch(\EE_1)-ch(\EE_2)$
is
$c_1(\EE_1)-c_1(\EE_2)+
\text{higher terms}$,
it is sufficient to show that
$\sum_q(-1)^qq\cdot
ch(\wedge^q\FF)$
is $(-1)^nc_{n-1}(\FF)+
\text{higher terms}$.
By splitting principle,
we put
$\sum c_i(\FF)T^{n-i}=
\prod_{i=1}^n(T+\alpha_i)$.
If we put
$a_i=\exp\alpha_i$,
the chern character
$ch(\wedge^q\FF)$
is the $q$-th fundamental symmetric polynomial
$s_q$ of $a_i$.
We consider a polynomial
$f(T)=\sum(-1)^qs_qT^{n-q}=
\prod(T-a_i)$
and compute the derivative
$f'(1)$ in two ways.
First we have
$f'(1)=
\sum(-1)^q(n-q)s_q=
n\sum(-1)^qs_q-
\sum(-1)^qqs_q$
Since
$\sum(-1)^qs_q=
f(1)=(-1)^nc_n(\FF)+
\text{higher terms}$,
we have
$\sum(-1)^qqs_q=-f'(1)$
modulo codimension $n$.
On the other hand
$f'(1)=\sum_i\prod_{j\ne i}
(1-a_j)=
(-1)^{n-1}c_{n-1}(\FF)+
\text{higher terms}$.
Hence we have
$\sum(-1)^qqs_q=
(-1)^nc_{n-1}(\FF)+
\text{higher terms}$
and Lemma is proved.
\enddemo

By Lemma applied to
$\EE_1=\EX^{(F)},
\EE_2=\Cal O_X (-D)^{(F)}$
and $\FF=\Omega^1_X(log D)$,
we have
$$
n^2_c(U,\tMM)
-rn^2_c(U,1)=
-(c_1(\EX^{(F)})-r
c_1(\Cal O_X (-D)^{(F)}))
\cdot
(-1)^{n-1}c_{n-1}(\Omega^1_X(log D)).
$$
Therefore to complete the proof
it is enough to show
$$(-c_1(\EX^{(F)}))
\cdot
(-1)^{n-1}c_{n-1}(\Omega^1_X(log D))=
\sum_ic_i\cdot
\Tr_{k_i/k}\Tr\ \res_i\nabla
$$
in $\ZZ^{\Spec k\otimes F}
\subset k\otimes F$.
We compute
$c_1(\EX^{(F)})$
by using the residue.

\proclaim{Lemma 5.4}
Let $X$ be a proper smooth scheme over
a field $k$ of characteristic 0
and $U$ be the complement of
a divisor $D$ with simple normal crossings.
Let $\nabla:\EX\to\EX\otimes\Omega^1_X(log D)$
be a logarithmic integrable connection.
Let $k_i$ be the constant field of
an irreducible component $D_i$ of $D$
and $\Tr\ \res_i\nabla\in k_i$.
Let $\partial:\bigoplus_ik_i\to
H^1(X,\Omega^1_X)$
be the boundary map of the exact sequence
$0\to\Omega^1_X\to\Omega^1_X(log D)\to\bigoplus_i\Cal O_{D_i}\to0$.
Then we have
$$
c_1(\EX)=
-\partial
(\Tr\ \res_i\nabla)_i.$$
\endproclaim
\demo{Proof}
The first chern class
$c_1(\EX)$
is the class of the
push of the $\Gm$-torsor
$\det \EX$
by $d\log:\Gm\to\Omega^1_X$.
The
$\Omega^1_X$-torser is trivialized
locally by $\nabla\log e$
for a local basis $e$ of $\det \EX$.
Hence  it is the inverse image of
$\sum_i\Tr\ \res_i\nabla\in\bigoplus_i\Cal O_{D_i}$
by $\Omega^1_X(log D)\to\bigoplus_i\Cal O_{D_i}$.
Since the boundary map is the minus of the
class of torser, Lemma is proved.
\enddemo

By Lemma 5.4, we have
$$\align
&(-c_1(\EX^{(F)}))
\cdot
(-1)^{n-1}c_{n-1}(\Omega^1_X(log D))\\
=&
\sum_i\Tr_{k_i/k}
(\res_i\nabla\cdot
\deg_{k_i}
((-1)^{n-1}c_{n-1}(\Omega^1_X(log D)|_{D_i})))\\
=&
\sum_ic_i\cdot
\Tr_{k_i/k}\Tr\ \res_i\nabla.
\endalign$$
Thus proof of Theorem 3 (1) is now complete.

(2).
By Corollary of Proposition 1 in Section 2
and by the definition of determinantally motivic,
$(\det_F\tMM,c)$
is motivic for $c\in CH^n(X\bmod D)$.
Hence it is sufficient to prove
\proclaim{Lemma 5.5}
For $\aaa\in\Bbb B^0$,
the triple $\tilde J(\aaa)$
over the abelian number field $k=F=F_{\aaa}$ is
motivic.
\endproclaim
\demo{Proof}
We may assume the coefficient
$n_a$ of $\aaa=\sum_an_a[a]$ is positive
$n_a\ge0$ for all $a\in \QZ$.
In fact we put
$\aaa^+=
\sum_{n_a>0}n_a[a]+
[-\sum_{n_a>0}n_aa]$
and
$\aaa^-=
\sum_{n_a<0}(-n_a)[a]+
[\sum_{n_a<0}n_aa]
\in \Bbb B^0$.
Then $\aaa^+$
and $\aaa^-$
satisfy the condition
$n_a\ge 0$,
fixed by $G_{\aaa}$
and $\aaa=\aaa^+-\aaa^-$.
Since the algebraic Hecke character
$J(\aaa)$ is the quotient
$J(\aaa^+)J(\aaa^-)^{-1}$,
we may assume
$n_a\ge 0$.
Further since
$J(\aaa)=J(\aaa_0)$
for $\aaa_0=\aaa-n_0[0]$,
we may assume $n_0=0$.

Let
$\aaa=\sum_{i=0}^{n+1}[a_i]$.
Clearly we may assume $n\ge 0$.
First we consider the case $n=0$.
Let $\aaa=[a/m]+[-a/m]$
for integers $(a,m)=1.$
The field $F=F_{\aaa}$
is then $\QQ(\cos\frac{2\pi}m)$.
We compute $\tilde J(\aaa)$.
For even $m$,
let $\epsilon$
be the character of order 2
corresponding to the quadratic extension
$\QQ(\cos\frac{\pi}m)$ of $F$
and $\tilde\epsilon$ be the triple
$(\epsilon,
(\tau\circ\sigma^{-1}(\exp\frac{a\pi\sqrt{-1}}{m}))_{\sigma,\tau},
0)$.
For odd $m$, we put $\tilde\epsilon=\bold1$.
We show $\tilde J(\aaa)=\tilde \epsilon\otimes \bold1(-1)$
and this proves Lemma for $n=0$.
First we compute the algebraic Hecke character
$J(\aaa)$.
Let $q$ be a finite place
of $F=F_{\aaa}$
prime to $2m$.
The norm $Nq$ satisfies $Nq=\pm1\bmod m$
since $Fr_q$ fixes $\aaa$.
If $Nq\equiv1\bmod m$, we have
$J(\aaa)(q)=
\tau(\chi,\psi)
\tau(\chi^{-1},\psi)
=
\chi(-1)Nq$
for the $m$-th power residue symbol
$\chi=
\fracwithdelims(){}q_m^a:
\kappa(q)^\times\to\mu_m$.
Further
$$\chi(-1)=
(-1)^{\frac{Nq-1}ma}=
\cases  1 & \text{if $m$ is odd or if $m$ is even and $Nq\equiv1\bmod 2m$}\\
-1 & \text{if $m$ is even and $Nq\equiv1+m\bmod 2m$.}
\endcases$$
If $Nq\equiv-1\bmod m$, we have
$J(\aaa)(q)=
\tau_{E_2}(\chi,\psi\circ\Tr_{E_2/\kappa(q)})$
where $E_2$
is the quadratic extension of $\kappa(q)$
and $\chi=
\fracwithdelims(){}{E_2}_m^a:
E_2^\times\to\mu_m$.
Using $\chi|_{\kappa(q)^\times}=1$,
an elementary computation yields
$J(\aaa)(q)=
(-1)^{\frac{Nq+1}ma}Nq$.
Hence we have
$$J(\aaa)(q)=
\cases  Nq & \text{if $m$ is odd or if $m$ is even and $Nq\equiv-1\bmod 2m$}\\
-Nq & \text{if $m$ is even and $Nq\equiv-1+m\bmod 2m$.}
\endcases$$
Thus $J(\aaa)=\epsilon\cdot N$, where $N$ is the norm, is proved.
By $\Gamma(s)\Gamma(1-s)=
2\pi\sqrt{-1}\cdot
\exp(\pi\sqrt{-1}\cdot s)
/(\exp(2\pi\sqrt{-1}\cdot s)-1)$,
we have $$\Gamma(c\aaa)=
2\pi\sqrt{-1}
\times
\cases 1 &\text{if $m$ is odd}\\
(\tau\circ\sigma^{-1}(\exp\frac{a\pi\sqrt{-1}}{m}))_{\sigma,\tau}
&\text {if $m$ is even}
\endcases$$
in $(F\otimes F)^\times\backslash
(F\otimes F\otimes \CC)^\times$.
Finally it is clear that $\langle\aaa\rangle=1$.
Thus we have proved $\tilde J(\aaa)=\tilde\epsilon\otimes \bold1(-1)$.

Assume $n\ge1$.
We show $\tilde J(\aaa)$
appears in the cohomology of a Fermat hypersurface and
hence is in the category of
motives generated by CM abelian varieties.
Let $m$ be an integer satisfying
$ma_i=0$ for $0\le i\le n+1$.
Let $U$ be the complement in
$X=\roman{Proj}\
F[T_0,\cdots,T_{n+1}]/
(\sum_i T_i)\simeq
\Bbb P^n_F$
of the divisor
$D=\bigcup_iD_i$
with simple normal crossings,
where $D_i=(T_i=0)$.
Let $X_m^n$ be the Fermat hypersurface
$\text{Proj }
F[T_0,\cdots,T_n]/(\sum_i T_i^m)$.
and
$\pi:X_m^n\to X$
be the covering $T_i\mapsto T_i^m$
etale over $U$.
Over the field $\QQ(\zeta_m)$,
it is a Galois covering with the
Galois group $G=\mu_m^{n+1}/\Delta(\mu_m)$.
The element $\aaa$
defines a character
$\aaa:G\to \QQ(\mu_m)^\times$
by $(\zeta_i)_i\mapsto\prod_i
\zeta_i^{m\langle a_i\rangle}$.
Let $p_{\aaa}=
\sum_{\zeta\in G}
\aaa(\zeta)^{-1}\zeta\in
\QQ(\mu_m)[G]$
be the projector.
Since $\aaa$ is fixed by $G_{F_{\aaa}}$,
it is in $F_{\aaa}[G]$.
Further the algebraic correspondence
$[p_{\aaa}]=
\sum_{\zeta\in G}
\aaa(\zeta)^{-1}[\Gamma_\zeta]$
where
$[\Gamma_\zeta]$
is the graph of the automorphism
$\zeta\in G$ on $X_{m,\QQ(\mu_m)}^n$
is stable by $G_{F_{\aaa}}$
and
is defined over $F_{\aaa}$.
Let $\tMM_{\aaa}$
be the $p_{\aaa}$-component of
the direct image $\pi_* 1$ on $U$
with the trivial filtration $Gr^0=\tMM_{\aaa}$.
We apply Theorem 3 (1) to
$\tMM_{\aaa}$ on $U$ and
show that
$$
H^n_c(U,\tMM_{\aaa})(-1)\simeq\tilde J(\aaa).$$
Since the left hand side is
in the category generated by
CM abelian varieties,
this complete the proof of (2).

It is well-known that
$H^q_c(U,\tMM_{\aaa})=0$
for $q\ne  n$
and
$H^n_c(U,\tMM_{\aaa})$
is of rank 1.
Also it is easy to check that
$\det_FR\Gamma_c(U,\bold1)\simeq
\bold1((-1)^{n+1})$.
Hence the left hand side of Theorem 3 (1)
for
$\tMM_{\aaa}$ on $U$
is
$(H^n_c(U,\tMM_{\aaa})(-1))^{(-1)^n}$.
We compute the right hand side.
The rank of $\tMM_{\aaa}$ over $F_{\aaa}$ is 1.
We compute $c_{X\bmod D}$.
The exact sequence
$$
H^{n-1}(X,\Kn)\to
\bigoplus_iH^{n-1}(D_i,\Kn)\to
CH^n(X\bmod D)\to
H^n(X,\Kn)\to0$$
is
$$F^\times\to
\bigoplus_iF^\times\to
CH^n(X\bmod D)\to
\ZZ\to 0$$
and a section
$\ZZ\to CH^n(X\bmod D):1\mapsto [x]$
is defined for an
$F$-rational point
$x=(t_0:\cdots:t_{n+1})
\in U(F)$.
We show that
$c_n(\Omega^1_X(log D),\res)=
[x]+(t_i)_i$.
Let $\omega=
\sum_{i=1}^{n+1}
t_id\log(T_i/T_0)$.
The zero locus of $\omega$
is $[x]$
and
$\res_{D_i}\omega=
t_i$ for $0\le i\le n+1$.
Using Lemma 2.1,
we have
$$
c_n(\Omega^1_X(log D),\res)=
[x]+\sum_i\{t_i\}\cdot
c_{n-1}(\Omega_{D_i}(\log D|_{D_i}))
=[x]+(t_i)_i.
$$
{}From this, by an elementary computation,
we have
$(\tMM_{\aaa},c_{X\bmod D})=1$.
We also see that
the formation of $c_{X\bmod D}$
commutes with base change
and we do not need assume
Gersten resolution in the proof of (1) for this case.
The assumption that the
spectral sequence degenerates at
$E_1$ is satisfied by
the degeneracy of the Hodge spectral sequence
for $X_m^n$.
Finally by
$c_i=(-1)^{n-1}$,
we see
$\aaa_{D,\tMM_{\aaa}}= (-1)^{n-1}\aaa$.
Therefore the right hand side of (1)
is $\tilde J(\aaa)^{(-1)^n}$.
Thus $H^n_c(U,\tMM_{\aaa})(-1)
\simeq\tilde J(\aaa)$
is proved and
proof of Theorem 3 is completed.
\enddemo
\enddemo

\enddocument